\magnification 1200
\input amssym
\input miniltx.tex
\input graphicx.sty
\input pictex
\input miniltx
\input color.sty

  \makeatletter
  \def \Gin @driver{pdftex.def} \resetatcatcode \font \bbfive = bbm5 \font \bbseven = bbm7 \font \bbten = bbm10  
  
  \font \eightbf = cmbx8
  \font \eighti = cmmi8 \skewchar \eighti = '177
  \font \eightit = cmti8
  \font \eightrm = cmr8
  \font \eightsl = cmsl8
  \font \eightsy = cmsy8 \skewchar \eightsy = '60
  \font \eighttt = cmtt8 \hyphenchar \eighttt = -1

  \font \sixi = cmmi6 \skewchar \sixi = '177
  \font \sixrm = cmr6
  \font \sixsy = cmsy6 \skewchar \sixsy = '60
  \font \tensc = cmcsc10
  
   \scriptfont \bffam = \bbseven \scriptscriptfont \bffam = \bbfive \textfont \bffam = \bbten
\newskip \ttglue
  \def \eightpoint {\def \rm {\fam 0 \eightrm }\relax
    \textfont 0= \eightrm
    \scriptfont 0 = \sixrm
    \scriptscriptfont 0 = \fiverm
    \textfont 1 = \eighti
    \scriptfont 1 = \sixi
    \scriptscriptfont 1 = \fivei
    \textfont 2 = \eightsy
    \scriptfont 2 = \sixsy
    \scriptscriptfont 2 = \fivesy
    \textfont 3 = \tenex
    \scriptfont 3 = \tenex
    \scriptscriptfont 3 = \tenex
    \def \it {\fam \itfam \eightit }\relax
    \textfont \itfam = \eightit
    \def \sl {\fam \slfam \eightsl }\relax
    \textfont \slfam = \eightsl
    \def \bf {\fam \bffam \eightbf }\relax
    \textfont \bffam = \bbseven
    \scriptfont \bffam = \bbfive
    \scriptscriptfont \bffam = \bbfive
    \def \tt {\fam \ttfam \eighttt }\relax
    \textfont \ttfam = \eighttt
    \tt
    \ttglue = .5em plus.25em minus.15em
    \normalbaselineskip = 9pt
    \def \MF {{\manual opqr}\-{\manual stuq}}\relax
    \let \sc = \sixrm
    \let \big = \eightbig
    \setbox \strutbox = \hbox {\vrule height7pt depth2pt width0pt}\relax \normalbaselines \rm }

  \def \setfont #1{\font \auxfont =#1 \auxfont }
  \long \def \withfont #1#2{{\setfont {#1}#2}}
  \def \TRUE {Y}
  \def \FALSE {N}
  \def \EMPTY {}
  \def \ifundef #1{\expandafter \ifx \csname #1\endcsname \relax }
  \def \undefrule {\kern 2pt \vrule width 5pt height 5pt depth 0pt \kern 2pt}
  \def \UndefLabels {}
  \def \possundef #1{\ifundef {#1}\undefrule {\eighttt #1}\undefrule \global \edef \UndefLabels {\UndefLabels #1\par
}\else \csname #1\endcsname \fi }
  \newcount \secno \secno = 0
  \newcount \stno \stno = 0
  \newcount \eqcntr \eqcntr = 0
  \ifundef {showlabel} \global \def \showlabel {\FALSE } \fi
  \ifundef {auxwrite} \global \def \auxwrite {\TRUE } \fi
  \ifundef {auxread} \global \def \auxread {\TRUE } \fi
  \def \define #1#2{\global \expandafter \edef \csname #1\endcsname {#2}} \long
  \def \error #1{\medskip \noindent {\bf ******* #1}}
  \def \fatal #1{\error {#1\par Exiting...}\end }
  \def \advseqnumbering {\global \advance \stno by 1 \global \eqcntr =0}
  \def \current {\ifnum \secno = 0 \number \stno \else \number \secno \ifnum \stno = 0 \else .\number \stno \fi \fi }
  \begingroup
  \catcode `\@=0
  \catcode `\\=11
  @global @def @textbackslash {\}
  @endgroup
  \def \space { }
  \def \deflabel #1#2{\if \TRUE \showlabel \hbox {\sixrm [[ #1 ]]} \fi \ifundef {#1PrimarilyDefined}\define
{#1}{#2}\define {#1PrimarilyDefined}{#2}\if \TRUE \auxwrite \immediate \write 1 {\textbackslash newlabel {#1}{#2}}\fi
\else \edef \old {\csname #1\endcsname }\edef \new {#2}\if \old \new \else \fatal {Duplicate definition for label ``{\tt
#1}'', already defined as ``{\tt \old }''.}\fi \fi }
  
  \def \label #1 {\deflabel {#1}{\current }}
  \def \equationmark #1 {\ifundef {InsideBlock} \advseqnumbering \eqno {(\current )} \deflabel {#1}{\current } \else
\global \advance \eqcntr by 1 \edef \subeqmarkaux {\current .\number \eqcntr } \eqno {(\subeqmarkaux )} \deflabel
{#1}{\subeqmarkaux } \fi }
  \def \split #1.#2.#3.#4;{\global
  \def \parone {#1}\global
  \def \partwo {#2}\global
  \def \parthree {#3}\global
  \def \parfour {#4}}
  \def \NA {NA}
  \def \ref #1{\split #1.NA.NA.NA;(\possundef {\parone }\ifx \partwo \NA \else .\partwo \fi )}
   \newcount \bibno \bibno = 0
  
  \def \Bibitem #1 #2; #3; #4 \par {\smallbreak \global \advance \bibno by 1 \item {[\possundef
{#1}]}#2, {``#3''},
#4.{\ifundef {USED#1}\color {red}\ (Not cited)\fi }\par \ifundef {#1PrimarilyDefined}\else \fatal {Duplicate definition for bibliography item ``{\tt #1}'', already
defined in ``{\tt [\csname #1\endcsname ]}''.}  \fi \ifundef {#1}\else \edef \prevNum {\csname #1\endcsname } \ifnum
\bibno =\prevNum \else \error {Mismatch bibliography item ``{\tt #1}'', defined earlier (in aux file ?) as ``{\tt
\prevNum }'' but should be ``{\tt \number \bibno }''.  Running again should fix this.}  \fi \fi \define
{#1PrimarilyDefined}{#2}\if \TRUE \auxwrite \immediate \write 1 {\textbackslash newbib {#1}{\number \bibno }}\fi }
  \def \jrn #1, #2 (#3), #4-#5;{{\sl #1}, {\bf #2} (#3), #4--#5}
  \def \Article #1 #2; #3; #4 \par {\Bibitem #1 #2; #3; \jrn #4; \par }
  \def \references {\begingroup \bigbreak \eightpoint \centerline {\tensc References} \nobreak \medskip \frenchspacing }
  \catcode `\@=11
  \def \c@itrk #1{{\bf \possundef {#1}}}
  \def \c@ite #1{{\rm [\c@itrk{#1}]}\define {USED#1}{NA}}
  \def \sc@ite [#1]#2{{\rm [\c@itrk{#2}\hskip 0.7pt:\hskip 2pt #1]}\define {USED#2}{NA}}
  \def \du@lcite {\if \pe@k [\expandafter \sc@ite \else \expandafter \c@ite \fi }
  \def \cite {\futurelet \pe@k \du@lcite } \catcode `\@=12
  \def \Headlines #1#2{\nopagenumbers \headline {\ifnum \pageno = 1 \hfil \else \ifodd \pageno \tensc \hfil \lcase {#1}
\hfil \folio \else \tensc \folio \hfil \lcase {#2} \hfil \fi \fi }}
  \def \title #1{\medskip \centerline {\withfont {cmbx12}{\ucase {#1}}}}
  \def \Subjclass #1#2{\footnote {\null }{\eightrm #1 \eightsl Mathematics Subject Classification: \eightrm #2.}}  \long
  \def \Quote #1\endQuote {\begingroup \leftskip 35pt \rightskip 35pt \parindent 17pt \eightpoint #1\par \endgroup }
  \long \def \Abstract #1\endAbstract {\vskip 1cm \Quote \noindent #1\endQuote }

  \def \Note #1{\footnote {}{\eightpoint #1}}
  \def \Date #1 {\Note {\it Date: #1.}}
  \def \today {\ifcase \month \or January\or February\or March\or April\or May\or June\or July\or August\or September\or
October\or November\or December\fi \space \number \day , \number \year }
   \newdimen \boxht
  \def \ifht #1#2#3{\ifdim \boxht <11pt #1\else \ifdim \boxht <13pt #2\else #3\fi \fi }
  \def \hparen #1{\setbox 0\hbox {$#1$}\boxht =\ht 0 \advance \boxht by \dp 0 \ifht {(}{\big (}{\Big (}#1\ifht {)}{\big
)}{\Big )}}

  \def \fix {\smallskip \noindent $\blacktriangleright $\kern 12pt}
  \def \iskip {\medskip \noindent }
  
  \def \newpage {\vfill \eject }
  \def \ucase #1{\edef \auxvar {\uppercase {#1}}\auxvar }
  \def \lcase #1{\edef \auxvar {\lowercase {#1}}\auxvar }
  \def \emph #1{{\it #1}}
  \def \section #1 \par {\global \advance \secno by 1 \stno = 0 \goodbreak \bigbreak \noindent {\bf \number \secno
.\enspace #1.}  \nobreak \medskip \noindent }
  \def \state #1 #2\par {\begingroup \def \InsideBlock {} \medbreak \noindent \advseqnumbering {\bf \current .\enspace #1.\enspace \sl #2\par }\medbreak
\endgroup }
  \def \definition #1\par {\state Definition \rm #1\par } \newcount \CloseProofFlag
  \def \closeProof {\eqno \endproofmarker \global \CloseProofFlag =1}
  
  \long \def \Proof #1\endProof {\begingroup
    \def \InsideBlock {} \global \CloseProofFlag =0 \medbreak \noindent {\it Proof.\enspace }#1 \ifnum \CloseProofFlag =0
\hfill $\endproofmarker $ \looseness = -1 \fi \medbreak \endgroup }
  \def \quebra #1{#1 $$$$ #1}
  \def \explica #1#2{\mathrel {\buildrel \hbox {\sixrm #1} \over #2}}

  \def \explain #1#2{\explica {\ref {#1}}{#2}}
  
  \def \=#1{\explain {#1}{=}}
  \def \pilar #1{\vrule height #1 width 0pt}

   \newcount \fnctr \fnctr = 0
  \def \fn #1{\global \advance \fnctr by 1 \edef \footnumb {$^{\number \fnctr }$}\footnote {\footnumb }{\eightpoint
#1\par \vskip -10pt}}
  \def \text #1{\hbox {#1}}
  \def \bool #1{[{#1}]}
  
  \def \Item #1{\smallskip \item {{\rm #1}}} \newcount \zitemno \zitemno = 0
  \def \iItemize {\global \zitemno = 0}
  \def \zitemplus {\global \advance \zitemno by 1 \relax }
  \def \rzitem {\romannumeral \zitemno }
  \def \rzitemplus {\zitemplus \rzitem }
  \def \iItem {\Item {{\rm (\rzitemplus )}}}
  
  \def \iItemmark #1 {\deflabel {#1}{\current .\rzitem } \deflabel {Local#1}{\rzitem }}
  \def \ds {\displaystyle }
  \def \and {\hbox {\quad and \quad }}

  \def \imply {\mathrel {\Rightarrow }}
  \def \IFF {\kern 7pt\Leftrightarrow \kern 7pt}

  \def \IMPLY {\kern 7pt \Rightarrow \kern 7pt}
  \def \for #1{\quad \forall \,#1}

  \def \endproofmarker {\square }
  \def \"#1{{\it #1}\/}
  \def \umlaut #1{{\accent "7F #1}}
  \def \inv {^{-1}}
  \def \*{\otimes }
  \def \caldef #1{\global \expandafter \edef \csname #1\endcsname {{\cal #1}}}
  \def \mathcal #1{{\cal #1}}

  \def \bfdef #1{\global \expandafter \edef \csname #1\endcsname {{\bf #1}}} \bfdef N \bfdef Z \bfdef C \bfdef R
  
  \if \TRUE \auxread \IfFileExists {\jobname .aux}{\input \jobname .aux}{\null } \fi
  \if \TRUE \auxwrite \immediate \openout 1 \jobname .aux \fi
  \def \close {\if \EMPTY \UndefLabels \else \message {*** There were undefined labels ***} \iskip ****************** \
Undefined Labels: \tt \par \UndefLabels \fi \if \TRUE \auxwrite \closeout 1 \fi \par \vfill \supereject \end }

  \def \prod {\mathop {\mathchoice {\hbox {$\mathchar "1351$}}{\mathchar "1351}{\mathchar "1351}{\mathchar "1351}}}
  \def \bigvee {\mathop {\mathchoice {\hbox {$\mathchar "1357$}}{\mathchar "1357}{\mathchar "1357}{\mathchar "1357}}}
  \def \bigwedge {\mathop {\mathchoice {\hbox {$\mathchar "1356$}}{\mathchar "1356}{\mathchar "1356}{\mathchar "1356}}}
  \def \sum {\mathop {\mathchoice {\hbox {$\mathchar "1350$}}{\mathchar "1350}{\mathchar "1350}{\mathchar "1350}}}
  \def \bigcup {\mathop {\mathchoice {\raise 1pt \hbox {$\mathchar "1353$}}{\mathchar "1353}{\mathchar "1353}{\mathchar "1353}}}
  \def \bigcap {\mathop {\mathchoice {\raise 1.5pt \hbox {$\mathchar "1354$}}{\mathchar "1354}{\mathchar "1354}{\mathchar "1354}}}

  \def \part #1#2{
    % \vskip0pt plus.3\vsize \penalty-500 \vskip0pt plus-.3\vsize \vskip 2cm \hrule    \vskip 0.5cm
    \vfill\eject
    \centerline {\withfont {cmr17 }{PART #1}}
    \bigskip
    \centerline {\withfont {cmr12}{#2}}
    \vskip 1cm}

  \font \gothicfont =eufb10
  \font \smallgothicfont =eufb10 scaled 833
  \font \vtt =cmvtt10

  \def \Imply {\ \mathrel {\Rightarrow }\ }
  \def \Iff {\ \mathrel {\Leftrightarrow }\ }
  \def \|{\mathrel {|}}
  \def \arw #1{\ {\buildrel #1 \over \longrightarrow }\ }
  \def \tabrule {\noalign {\hrule }}

  \def \I {{\cal I}}
  \def \P {{\cal P}}
  \def \F {{\cal F}}
  \def \E {{\cal E}}
  \def \S {{\cal S}}
  \def \B {{\cal B}}

  \def \goth #1{\hbox {\gothicfont #1}}
  \newdimen \largura \largura =5pt
  \def \gothEhat {\goth E\kern -\largura \widehat {\vrule height 6.5pt width \largura depth -10pt}}
  \def \hull {\goth H(S)}
  \def \ehull {\goth E(S)}

  \def \tight {{\hbox {\sixrm tight}}}
  \def \sixrmbox #1{{\hbox {\sixrm #1}}}
  \def \ess {\sharp }

  \def \X {{\cal X}}
  \def \subX {{\kern -1pt\scriptscriptstyle \X }}
  \def \SX {S_\subX }

  \def \K {{\bf K}}
  
  \def \BA {{\cal B}_A}
  \def \eone {\goth E_1(S)}

  \def \supop {^{_\sixrmbox {op}}}
  \def \subinf {_{^\infty }}

  \def \spec {\gothEhat (S)}
  \def \varspec #1{\gothEhat #1(S)}
  \def \seone {\varspec {_1}}
  \def \sinf {\varspec {_{\scriptscriptstyle \infty }}}

  \def \sinfop {\varspec {\subinf \supop }}
  \def \sop {\varspec {_\sixrmbox {op}}}

  \relax

  \def \id {\hbox {id}}
  \def \tS {\tilde S}
  \def \sla {semilattice}
  \def \uf {ultrafilter}
  \def \itmproof #1 {\medskip \noindent #1\enspace }
  \def \interior #1{\mathaccent '27#1}
  \definecolor {pink}{rgb}{1.0,0.21,0.7}
  \definecolor {grey}{rgb}{0.6,0.6,0.6}
  \definecolor {orange}{rgb}{1.0,0.5,0.0}
  \def \exists {\mathchar "0239\kern 1pt }
  
  \def \aster 0 {\ast _0}
  \def \Bool #1{\left [#1\right ]}
  \def \I {{\cal I}}
  \def \itemproof #1{\bigskip \noindent (#1)\ }

  \newcount \aitemno \aitemno = -1
  \def \boxlet #1{\hbox to 6.5pt{\hfill #1\hfill }}
  
  \def \aitemconv {\ifcase \aitemno a\or b\or c\or d\or e\or f\or g\or
h\or i\or j\or k\or l\or m\or n\or o\or p\or q\or r\or s\or t\or u\or
v\or w\or x\or y\or z\else zzz\fi }
  \def \aitem {\global \advance \aitemno by 1\Item {(\boxlet \aitemconv )}}
  \def \aitemmark #1 {\deflabel {#1}{\aitemconv }}

\def \writeIndex #1#2#3#4{\write #1 {\textbackslash #2 {#3}{#4}{\folio }}}

\def \printIndex #1#2#3#4{\noindent \hbox to #4{\hfill #1\hfill }\quad \vrule width 0pt depth 7pt #2 \dotfill \quad #3 \par }

\def \# #1; #2; {\writeIndex 3{symbol}{#1}{#2}}
\def \section #1 \par {\color {black}\global \advance \secno by 1 \stno = 0
  \goodbreak \bigbreak
  \noindent {\bf \number \secno .\enspace #1.}  \nobreak \medskip \noindent
  \writeIndex 2{contents}{\number \secno .}{#1}}

\openout 3 symbols.aux

\def \# #1; #2; {\writeIndex 3{symbol}{#1}{#2}}

\openout 3 symbols.aux

\centerline {\bf REPRESENTATIONS OF THE INVERSE HULL OF}
\smallskip
\centerline {\bf A $0$\kern 0.6pt-LEFT CANCELLATIVE SEMIGROUP}
\medskip
\centerline {\tensc R. Exel and B. Steinberg}
%\footnote {}{\eightrm Typeset on \today .}
\Subjclass {2010}{20M18, 46L55}

\Abstract
  \def\smallhull{\hbox{\smallgothicfont H}(S)}%
  \def\smallehull{\hbox{\smallgothicfont E}(S)}%
  A semigroup $S$ containing a zero element is said to be $0$-left cancellative if $st =sr\neq0$ implies that $t =r$.
  Given such an $S$ we build an inverse semigroup $\smallhull$, called the inverse hull of $S$.
  Motivated by the study of certain C*-algebras associated to $\smallhull $ (a task that we will address in a subsequent
article) we carry out a detailed analysis of the spectrum of the idempotent semilattice $\smallehull$  of
$\smallhull$ with a special interest in identifying the ultra-characters.
  In order to produce examples of characters on $\smallehull$, we introduce the notion of \emph{strings} in a semigroup,
atempting to make sense of the \emph{infinite paths} which are of great importance in the study of graph C*-algebras.
  Our strongest results are obtained under the assumption that $S$ admits \emph{least common multiples}, but we also touch upon
the notion of \emph{finite alignment}, motivated by the corresponding notion from the theory of higher rank graphs, and
which has also appeared in recent papers by Spielberg and collaborators.
  \endAbstract

\vskip 1cm

  \IfFileExists {./contents.aux}{
    \bigbreak
    \centerline {\tensc Contents}
    \nobreak \bigskip
    \catcode `\@=11
    \input contents.aux
    \catcode `\@=12
    \vskip 1cm
  }{\bigskip \noindent \vtt *** File ``contents.aux" is missing.  It might be generated after a rerun. *** \bigskip }

\openout 2 contents.aux

\newpage

\section Introduction

The theory of semigroup C*-algebras has a long history, beginning with Coburn's work \cite {CoOne} and \cite {CoTwo}
where the C*-algebra of the additive semigroup of the natural numbers is studied in connection to Toeplitz operators.
  In \cite {MurOne} Murphy generalized this construction to the positive cone of an ordered group, and later to left
cancellative semigroups (\cite {MurTwo}, \cite {MurThree}).
  The C*-algebras studied by Murphy turned out to be too wild, even for nice looking semigroups such as ${\bf N}\times
{\bf N}$, and this prompted Li \cite {Li} to introduce an alternative C*-algebra for a left cancellative semigroup.  By
definition a semigroup $S$ is said to be left cancellative provided, for every $r,s,t\in S$, one has that
  $$
  s t =s r \Imply t =r.
  \equationmark Cancellative
  $$

Many interesting semigroups possess a zero element, namely an element $0$ such that
  $$
  s0=0s=0,
  $$
  for every $s$, and it is obvious that the presence of a zero prevents a semigroup from being left cancellative.  In
this work we focus on $0$-left cancellative semigroups, meaning that \ref {Cancellative} is required to hold only when
the terms in its antecedent are supposed to be nonzero.  This dramatically opens up the scope of applications including
a wealth of interesting semigroups, such as those arising from subshifts and, more generally, languages over a fixed
alphabet.  This also allows for the inclusion of the semigroupoids of \cite [Section 14]{actions} and left cancellative
categories, once the multiplication is extended to all pairs of elements by setting it to be zero whenever it is not
already defined.  See Section \ref {ExamplesSection} for more examples.

Starting with a $0$-left cancellative semigroup $S$, the crucial point is first to build an inverse semigroup $\hull $,
which we call \emph {the inverse hull} of $S$, in analogy with the notion of the inverse hull of a left cancellative
semigroup, cf.\ \cite [Section~1.9]{CP} and \cite {Cherubini}.  Once in possession of the inverse hull, one may invoke
any of the now standard constructions of C*-algebras from inverse semigroups, such as
  the tight C*-algebra \cite {actions} or
  Paterson's \cite {Paterson} universal C*-algebra.  Indeed
  our initial motivation was to study such C*-algebras, but the present work is instead focused on the passage from the
original semigroup to its inverse hull, rather than the much better understood passage from there to the associated
C*-algebras.  Particularly demanding is the work geared towards understanding the idempotent semilattice of $\hull $,
which we denote by $\ehull $, as well as its spectrum.  By a standard gadget $\ehull $ is put in correspondence with a
subsemilattice of the power set of $S\setminus \{0\}$, whose members we call the \emph {constructible sets}, by analogy
with a similar concept relevant to Li's work in \cite {Li}.

Our proposal, to be further developed in the second part of this article (curently in preparation), is to consider the
tight groupoid of $\hull $.  The unit space of this groupoid is well know to be the tight spectrum of the semilattice
$\ehull $, so it is crucial to understand the tight characters and, in view of \cite [12.9]{actions}, also the
ultra-characters.

One of the main examples motivating our pursuit of the present line of ideas is the path semigroup associated to a graph
$E$.  See our discussion at the end of section \ref {ExamplesSection} for a definition of this semigroup.  According to
\cite {KPRR}, one may associate to any locally finite directed graph $E$ (without sinks for simplicity), a C*-algebra
$C^*(E)$, which happens to coincide with the C*-algebra of a canonically associated \'etale groupoid ${\cal G}_E$, whose
unit space turns out to be the infinite path space of $E$.  We therefore set out to look for ways of producing a
canonically defined \'etale groupoid ${\cal G}(S)$ from any $0$-cancellative semigroup $S$, generalizing the
construction of ${\cal G}_E$ from the path semigroup of $E$.

In the case of the path semigroup of a graph, the ultra-characters correspond to infinite paths (see also \cite
[19.11]{actions}) so, attempting to make sense of infinite paths on an arbitrary semigroup, we introduce the concept of
strings \ref {DefineStrin}, an idea already present in \cite [19.10]{actions}, and which is one of our fundamental tools
when studying the spectrum of the {\sla } of constructible sets.  However there are $0$-left cancellative semigroups in
the wild where the nice relationship between ultra-characters and maximal paths (as observed in the graph case) is all
but lost, requiring a much more detailed analysis, which is carried out in section \ref {StringSection}.  See in
particular the example presented after \ref {BigNewTightResult}.

Regarding the problem of fully understanding the spectrum of $\ehull $, including the identification of the tight and
ultra-characters, we believe the present work represents only a modest beginning in a mammoth task lying ahead.  This
impression comes from situations in which similar spectra have been more or less understood, such as in \cite {infinoa}
and in \cite {DokuchaExel}, illustrating the high degree of complexity one should expect.

One of the main working hypothesis adopted in this work is the existence of \emph {least common multiples}: if $s$ and
$t$ are elements of the semigroup $S$, we say that $r$ is a least common multiple of $s$ and $t$, provided $r$ is a
common multiple of $s$ and $t$, and $sS\cap tS=rS$.  See Definition \ref {DefLCM} for more details.  While the exact
form of this notion does not seem to be present in the literature, it is of course motivated by the usual notion from
arithmetic as well as similar notions extensively employed in the literature of semigroup C*-algebras, such as in \cite
{AielloContiRossiStammeier}, \cite {BrownloweLarsenStammeier}, \cite {KwasniewskiLarsen}, \cite {Stammeier} and \cite
{Starling}.

Under the assumption that the semigroup admits least common multiples we are able to prove some of our strongest
results, beginning with the description in Corollary \ref {FormOfHull} of a \emph {normal form} for elements of the
inverse hull, including of course constructible sets.

Roughly speaking, the difficulty in proving major results in any area of mathematics is inversely proportional to the
strength of the chosen set of axioms.  The huge generality of semigroups allowed by our limited set of conditions
certainly makes that task very difficult but still we believe we have managed to prove a result we think will find
interesting applications, namely Theorem \ref {NonOpenUltra}, which essentially characterizes the set of all
ultra-characters (and hence the tight spectrum, by taking closures), although it depends on the ad-hoc knowledge of
ground ultra-characters \ref {DefineGroundChars}.

Besides semigroups admitting least common multiples we also study \emph {finitely aligned semigroups}, namely semigroups
whose finitely generated right ideals are close under intersection (sometimes called Howson semigroups in analogy with
\cite {Howson}), and which holds true for the semigroup of finite paths on a finitely aligned higher rank graph \cite
{KumjianPask}.

In a sense the present work should be thought of as a continuation of the study of semigroupoids started in \cite
[Section 14]{actions} and \cite {semigpds}.  In fact, as discussed in \ref {Semigroupoids}, given a semigroupoid, one
may set the undefined products to be zero and thus obtain a semigroup with many of the relavant properties studied here.
However, the specific associativity axiom assumed at the beginning of \cite [Section 14]{actions} is too strong and
excludes many intersting examples which we are now able to treat.

A semigroup with zero $S$ is called \emph {categorical at zero} (see Definition \ref {CatAtZero}), provided for every
$r,s,t\in S$, one has that
  $$
  rs\neq 0, \hbox { and } st\neq 0 \Imply rst\neq 0.
  $$
  Semigroups arising from semigroupoids, as in \ref {Semigroupoids}, are easily seen to be categorical at zero, but
semigroups defined from subshifts, such as in Example \ref {DefineSubshSgrp} do not share this property unless the
subshift is Markov.

Another well known condition usually considered in the study of semigroups is the existence of \emph {right local
units}.  By this we mean that, for every $s$ in $S$, there exists an idempotent element $e$ in $S$ such that $s = se$.
In the specific case of $0$-left cancellative semigroups this is in fact equivalent to saying that $s\in sS$, for every
$s$ in $S$ (see \ref {LocalUnits}).  Unital semigroups of course have right local units as do the semigroups arising
from categories.  However, once more the semigroups coming from subshifts are excluded.

Since one of our main motivations is to be able to apply our theory to subshift semigroups, in many of our general
results we have strived to avoid assuming strong hypotheses such as being categorical at zero or the existence of right
local units.

Last but not least we should mention the work by Spielberg and collaborators on left cancellative small categories
(\cite {SpielbergA}, \cite {SpielbergB} and \cite {BedosSpielberg}), which goes very much in the same direction we are
heading, with some significant differences in hypothesis.  On the one hand the papers mentioned above only deal with
categories, which may be viewed as special cases of the semigroupoids of \cite {actions} and \cite {semigpds}, but on
the other hand they transcend the singly aligned assumption of \cite [20.1]{actions} by thoroughly exploring the
finitely aligned situation (and to a certain extent the infinitely aligned case as well).  However, as already
mentioned, semigroups from categories are categorical at zero, and hence they exclude the main examples we have in mind,
namely subshift semigroups.

One of the advantages of our theory is that it can work with quotients of left cancellative categories by ideals.  If
one quotients a left cancellative category by an ideal (e.g., quotient the free monoid or semigroup by the ideal of
non-factors of a subshift), the resulting category is no longer left cancellative but the corresponding semigroup is
$0$-left cancellative.

The results in this paper have already been announced in \cite {announce}.  In addition, in a forthcoming paper we will
apply the results obtained here to study $0$-left cancellative semigroups arising from subshifts and their relationship
to various C*-algebras that have appeared in the literature motivated by Matsumoto's original work \cite {MatsuOri},
such as the Carlsen-Matsumoto C*-algebras of \cite {MatsuCarl}.  See also \cite {CarlsenSilvestrov}.

\part {ONE}{Semigroups}

\section Representations of semigroups

Let $S$ be a semigroup, namely a nonempty set equipped with an associative operation.

A \emph {zero element}\/ for $S$ is a (necessarily unique) element $0\in S$, satisfying
  $$
  s0 = 0s = 0, \for s\in S.
  $$

In what follows we will fix a semigroup $S$ possessing a zero element. Note that one can always adjoin a zero element to
a semigroup.  The set of idempotent elementss of $S$ will be denoted by $E(S)$.

\definition \label DefRepre
  Let $\Omega $ be any set.  By a \emph {representation of $S$ on $\Omega $} we shall mean any map
  $$
  \pi \colon S\to \I (\Omega ),
  $$
  \# $\textbackslash I (\Omega )$; Symmetric inverse semigroup on the set $\Omega $;
  \# $\pi $; Representation of a semigroup;
  where $\I (\Omega )$ is the symmetric inverse
  semigroup\fn {The symmetric inverse semigroup on a set $\Omega $ is the inverse semigroup formed by all partially
defined bijections on $\Omega $.}  on $\Omega $, such that
  \iItemize
  \iItem $\pi _0$ is the empty map on $\Omega $, and
  \iItem $\pi _s\circ \pi _t=\pi _{st}$, for all $s$ and $t$ in $S$.

Given a set $\Omega $, and any subset $X\subseteq \Omega $, let
  $\id _X$
  \# $\id _X$; Identity function on the set $X$;
  denote the identity function on $X$, so that $\id _X$ is
  an element of $E\big (\I (\Omega )\big )$, the \emph {idempotent {\sla }} of $\I (\Omega )$.  One in fact has that
  $$
  E\big (\I (\Omega )\big ) = \{\id _X: X\subseteq \Omega \},
  $$
  \# $\textbackslash P(\Omega )$; The set of all subsets of $\Omega $;
  so we may identify $E\big (\I (\Omega )\big )$ with the meet {\sla } $\P (\Omega )$ formed by all subsets of $\Omega
$.

\definition \label IntroEsFs
  Given a representation $\pi $ of $S$, for every $s$ in $S$ we will denote the domain of $\pi _s$ by $F^\pi _s$, and
the range of $\pi _s$ by $E^\pi _s$, so
  \# $F^\pi _s$; Domain of $\pi _s$;
  \# $E^\pi _s$; Range of $\pi _s$;
  \# $F _s$; Simplified notation for domain of $\pi _s$;
  \# $E _s$; Simplified notation for range of $\pi _s$;
  that $\pi _s$ is a bijective mapping
  $$
  \pi _s\colon F^\pi _s \to E^\pi _s.
  $$
  We will moreover let
  \# $f^\pi _s$; Identity function on $F^\pi _s$;
  \# $e^\pi _s$; Identity function on $E^\pi _s$;
  $$
  f^\pi _s: = \pi _s\inv \pi _s = \id _{F^\pi _s}
  \and
  e^\pi _s:= \pi _s\pi _s\inv = \id _{E^\pi _s}.
  $$

If $\pi $ is a representation of $S$ on a set $\Omega $, and if $\Omega '$ is a subset of $\Omega $ such that
  $$
  F^\pi _s \subseteq \Omega '\ \and E^\pi _s \subseteq \Omega ', \for s\in S,
  $$
  then evidently $\pi $ may be considered as a representation on $\Omega '$.  Moreover any point of $\Omega \setminus
\Omega '$ will have little relevance for $\pi $.

An example of a subset of $\Omega $ satisfying the above is clearly obtained by taking
  \# $\Omega _\ess $; Essential subset;
  $$
  \Omega _\ess = \big (\bigcup _{s\in S}F^\pi _s \big ) \cup \big (\bigcup _{s\in S} E^\pi _s\big ),
  \equationmark EssSubset
  $$
  which we will henceforth refer to as the \emph {essential subset} for $\pi $.

\definition \label DefEssRep
  A representation $\pi $ of $S$ is said to be \emph {essential}\/ provided $\Omega _\ess = \Omega $.

\state Proposition \label InvertEssential
  Let $\pi $ be a representation of a unital semigroup $S$ on a set $\Omega $, and let $u$ be an invertible element of
$s$.  Then
  $$
  \Omega _\ess = F^\pi _u = E^\pi _u.
  $$

\Proof
  First note that $F^\pi _1=E^\pi _1$ as $\pi _1$ is idempotent.  Since $\pi _1\pi _s=\pi _s=\pi _s\pi _1$, it follows
that
  $$
  F^\pi _s,E^\pi _s\subseteq F^\pi _1=E^\pi _1
  $$
  for all $s\in S$ and hence $\Omega _\ess =F^\pi _1=E^\pi _1$.  If $u$ is invertible with inverse $v$, then $\pi _u\pi
_v=\pi _1=\pi _v\pi _u$ shows that $F^\pi _1\subseteq F^\pi _u$ and $E^\pi _1\subseteq E^\pi _u$.  We deduce that
$\Omega _\ess =F^\pi _u=E^\pi _u$ as required.  \endProof

Let us fix, for the time being, a representation $\pi $ of $S$ on $\Omega $.  Whenever there is only one representation
in sight we will drop the superscripts in $F^\pi _s$, $E^\pi _s$, $f^\pi _s$, and $e^\pi _s$, and adopt the simplified
notations $F_s$, $E_s$, $f_s$, and $e_s$.

\state Proposition \label Covar
  Given $s$ and $t$ in $S$, one has that
  \iItemize
  \iItem $\pi _se_t=e_{st}\pi _s$, and
  \iItem $f_t\pi _s=\pi _sf_{ts}$.

\Proof
  We have
  $$
  \pi _s e_t =
  \pi _s\pi _s\inv \pi _s e_t =
  \pi _sf_s e_t =
  \pi _se_t f_s =
  \pi _s\pi _t \pi _t\inv \pi _s\inv \pi _s =
  \pi _{st}\pi _{st}\inv \pi _s =
  e_{st}\pi _s.
  $$
  As for (ii), we have
  $$
  f_t \pi _s =
  f_t \pi _s \pi _s\inv \pi _s =
  f_t e_s \pi _s =
  e_s f_t \pi _s =
  \pi _s\pi _s\inv \pi _t\inv \pi _t\pi _s =
  \pi _s\pi _{ts}\inv \pi _{ts} =
  \pi _sf_{ts}. \closeProof
  $$
  \endProof

\definition \label GenSgANdConstr
  \iItemize
  \iItem The inverse subsemigroup of $\I (\Omega )$ generated by the set $\{\pi _s : s \in S\}$ will be denoted by
  $\I (\Omega ,\pi )$.
  \# $\textbackslash I (\Omega ,\pi )$; The inverse semigroup generated by the range of a representation of $\pi $;
  \iItem Given any $X\in \P (\Omega )$ such that $\id _X$ belongs to $E\big (\I (\Omega ,\pi )\big )$, we will say $X$
is a \emph {$\pi $-constructible subset}.
  \iItem The collection of all $\pi $-constructible subsets of $\Omega $ will be denoted by $\P (\Omega ,\pi )$.  In
symbols
  $$
  \P (\Omega ,\pi )= \big \{X\in \P (\Omega ) : \id _X \in E\big (\I (\Omega ,\pi )\big )\big \}.
  $$
  \# $\textbackslash P (\Omega ,\pi )$; The collection of all $\pi $-constructible sets;
  Observe that by \ref {IntroEsFs}, one has that $E_s$ and $F_s$ are $\pi $-constructible sets.  For the special case of
$s=0$, we have $E_s=F_s=\emptyset $, so the empty set is $\pi $-constructible as well.

Since $\P (\Omega ,\pi )$ corresponds to the idempotent {\sla } of $\I (\Omega ,\pi )$ by definition, it is clear that
$\P (\Omega ,\pi )$ is a {\sla }, and in particular the intersection of two $\pi $-constructible sets is again $\pi
$-constructible.  In what follows we would like to characterize the $\pi $-constructible sets.

\state Lemma \label Conjuga
  For every $s$ in $S$, and every $X\in \P (\Omega )$, let
  $$
  s[X] := \pi _s(F_s\cap X), \and s\inv [X] := \pi _s\inv (E_s\cap X).
  $$
  One then has that
  \iItemize
  \iItem $\pi _s\,\id _X\,\pi _s\inv = \id _{s[X]}$, and
  \iItem $\pi _s\inv \,\id _X\,\pi _s= \id _{s\inv [X]}$.

\Proof
  We have
  $$
  \pi _s\,\id _X\,\pi _s\inv =
  \pi _sf_s\,\id _X\,\pi _s\inv =
  \pi _s\,\id _{F_s}\,\id _X\,\pi _s\inv \quebra =
  \pi _s\,\id _{F_s\cap X}\,\pi _s\inv =
  \id _{\pi _s(F_s\cap X)} =
  \id _{s[X]},
  $$
  proving (i).  A similar argument proves (ii).
  \endProof

\state Proposition \label CharactSL
  The family $\P (\Omega ,\pi )$ of $\pi $-constructible sets is the smallest subset of $\P (\Omega )$ containing every
$E_s$, and which is invariant under the maps
  $$
  X\mapsto s\inv [X], \and X\mapsto s[X],
  $$
  for all $s$ in $S$.

\Proof
  Given any $s$ in $S$, we have already seen that $E_s\in \P (\Omega ,\pi )$.  Furthermore, given $X$ in $\P (\Omega
,\pi )$, one has that $\ \pi _s\inv \,\id _X\,\pi _s\ $ and $\ \pi _s\,\id _X\,\pi _s\inv \ $ are both idempotent
elements of $\I (\Omega ,\pi )$, so we may deduce from \ref {Conjuga} that $s\inv [X]$ and $s[X]$ belong to $\P (\Omega
,\pi )$.

This proves that $\P (\Omega ,\pi )$ satisfies the conditions mentioned in the statement, and it therefore remains to
prove that $\P (\Omega ,\pi )$ is the smallest such collection.  In other words, given any collection $\F $ of subsets
of $\Omega $ satisfying the conditions in the statement, we must show that $\P (\Omega ,\pi )\subseteq \F $.  In order
to do this, pick any $X$ in $\P (\Omega ,\pi )$, so there exists some $\alpha $ in $\I (\Omega ,\pi )$ such that
  $$
  \id _X=\alpha \alpha \inv .
  $$

By definition of $\I (\Omega ,\pi )$, we have that $\alpha $ may be written as a product $\alpha =\alpha _1\alpha
_2\ldots \alpha _n$, where, for each $i$, there is an ${s_i}\in S$, such that either $\alpha _i=\pi _{s_i}$, or $\alpha
_i=\pi _{s_i}\inv $.

We will accomplish our goal of showing that $X\in \F $ by induction on $n$.  If $n=1$, and if $\alpha _1=\pi _{s_1}$,
then
  $$
  \id _X = \alpha _1\alpha _1\inv = \pi _{s_1}\pi _{s_1}\inv = \id _{E_{s_1}},
  $$
  whence $X=E_{s_1}$, so it lies in $\F $, by hypothesis.  Still under the assumption that $n=1$, but supposing now that
$\alpha _1=\pi _{s_1}\inv $, we have
  $$
  \id _X =
  \alpha _1\inv \alpha _1 =
  \pi _{s_1}\inv \pi _{s_1} =
  \id _{F_{s_1}},
  $$
  so
  $$
  X=F_{s_1} = s_1\inv [E_{s_1}] \in \F .
  $$

  Next assume that $n>1$, and let $\beta =\alpha _2\ldots \alpha _n$, so that
  $\beta \beta \inv = \id _Y$, where $Y$ lies in $\F $ by the induction
hypothesis.  Moreover
  $$
  \id _X = \alpha \alpha \inv = \alpha _1\beta \beta \inv \alpha _1\inv =
  \alpha _1\id _Y\alpha _1\inv .
  $$
  It therefore follows from \ref {Conjuga} that $X$ is either equal to $s_1[Y]$
or to $s_1\inv [Y]$, according to whether $\alpha _1=\pi _{s_1}$ or $\alpha
_1=\pi _{s_1}\inv $.  In any case we conclude that $X\in \F $, completing the
proof.
  \endProof

\section Cancellative semigroups

Beginning with this section we will restrict our attention to semigroups possessing certain special properties regarding
cancellation.

\definition
  Let $S$ be a semigroup containing a zero element.  We will say that $S$ is
\emph {$0$-left cancellative}, or \emph {left cancellative away from zero} if,
for every $r,s,t\in S$,
  $$
  s t =s r \neq 0 \Imply t =r,
  $$
  and \emph {$0$-right cancellative} if
  $$
  t s =r s \neq 0 \Imply t =r.
  $$
  If $S$ is both $0$-left cancellative and $0$-right cancellative, we will say
that $S$ is \emph {$0$-cancellative}.  Adjoining a zero to a left cancellative semigroup will result in a $0$-left cancellative semigroup and
so our study will subsume the classical case.

\fix In what follows we will fix a $0$-left cancellative semigroup $S$.  Occasionally, we will also assume that $S$ is
$0$-right cancellative.

If $X\subseteq S$ and $s\in S$, let us write $s\inv X$
  for the preimage of $X$ under left multiplication by $s$, namely
  $$
  s\inv X:=\{t\in S: st\in X\}.
  $$

For any $s$ in $S$ we will let
  $$
  F_s = \{x\in S: sx\neq 0\}=s\inv (S\setminus \{0\}),
  $$
  \#   $F_s$; Domains for the regular representation;
  and
  $$
  E_s = \{y\in S: y=sx\neq 0, \hbox { for some } x\in S\}=sS\setminus \{0\}.
  $$
  \#   $E_s$; Ranges for the regular representation;

Observe that the correspondence ``$x\to sx$'' gives a map from $F_s$ onto $E_s$,
which is onto by definition of $E_s$ and
one-to-one by virtue of $0$-left cancellativity.

\definition \label DefineThetaS
  For every $s$ in $S$ we will denote by $\theta _s$ the bijective mapping given
by
  $$
  \theta _s: x\in F_s \mapsto sx\in E_s.
  $$
  \# $\theta $; Regular representation;

Observing that $0$ is neither in $F_s$, nor in $E_s$, we see that these are both
subsets of
  $$
  S':= S\setminus \{0\},
  \equationmark IntroSPrime
  $$
  \# $S'$; Set of all nonzero elements of $S$;
  so we may view $\theta _s$ as a partially defined bijection on $S'$, which is
to say that $\theta _s\in \I (S')$.  We also notice that when $s=0$, both
$F_s$ and $E_s$ are empty, so $\theta _s$ is the empty map.

\state Proposition \label DefineRegRep
  The correspondence
  $$
  s \in S \mapsto \theta _s\in \I (S')
  $$
  is a representation of $S$ on $S'$, henceforth called the \emph {regular
representation of} $S$.

\Proof
  As already seen, $\theta _0$ is the empty map on $\I (S')$, so it suffices to
check \ref {DefRepre.ii}.
  Notice that a given element $x$ in $S'$ lies in the domain of $\theta _s \circ
\theta _t$ if and only if $tx\neq 0$, and $s(tx)\neq 0$.  These two conditions are
obviously equivalent to $(st)x\neq 0$, which is to say that $x$ lies in the domain
of $\theta _{st}$.  Moreover, for any $x$ in this common domain we have
  $$
  \theta _s\big (\theta _t(x)\big ) = s(tx) = (st)x = \theta _{st}(x),
  $$
  so $\theta _s \circ \theta _t=\theta _{st}$.
  \endProof

Regarding the notations introduced in \ref {IntroEsFs} in relation to the
regular representation, notice that
  $$
  F_s=F^\theta _s, \and E_s=E^\theta _s.
  $$

So far nothing guarantees that the left regular representation is essential, so let us now study its essential subset,
beginning with the following trivial fact whose easy proof is left for the reader.

\state Lemma
  Given an element $s$ in $S$, one has that
  \medskip
  $$
  \matrix {
  \ds s\in \bigcup _{t\in S}E^\theta _t & \iff &s\in S^2, & \quad \hbox {and}\cr \pilar {16pt}
  \ds s\in \bigcup _{t\in S}F^\theta _t & \iff & Ss\neq \{0\}.}
  $$

  \medskip
We thus see that $s$ fails to be in the essential subset of $\theta $ if and only if $s$ possesses the property
defined below:

\definition \label DefDegenElt
  A nonzero element $s$ in $S$ is said to be \emph {degenerate} if
  $$
  s\notin S^2, \and Ss=\{0\}.
  $$

The following is thus a simple interpretation of the terms involved:

\state Proposition \label RegRepEssential
  Denoting the essential subset for $\theta $ by $S'_\ess $, one has that
  $$
  S'\setminus S'_\ess = \{s\in S': s \hbox { is degenerate}\}.
  $$
  Therefore $\theta $ is essential if and only if $S$ possesses no degenerate elements.

So far nothing guarantees that the left regular representation is injective, but in case injectivity of $\theta $ is desired,
let us now discuss the appropriate conditions for this.

\state Definition \label RightReductive \rm A semigroup $S$ is called \emph {right reductive} if it acts faithfully on
the left of itself, that is, $sx=tx$ for all $x\in S$ implies $s=t$.

Of course every unital semigroup is right reductive.  If $S$ is a right reductive $0$-left cancellative semigroup, then
it embeds in $\I (S')$ via $s\mapsto \theta _s$.

Observe that if $S$ is $0$-right cancellative, then a single $x$ for which $sx=rx$, as long as this is nonzero, is
enough to imply that $s=t$.  So, in a sense, right reductivity is a weaker version of $0$-right cancellativity.

\state Proposition \label sSEqualSImpliesInvert
  Suppose that, besides being $0$-left cancellative, $S$ is also right reductive.  If $s$ is an element of $S$ such that
$sS=S$, then $S$ is a unital semigroup and $s$ is invertible.

\Proof
  Choosing $u$ in $S$ such that $su=s$, we will prove that $u$ is an identity for $S$.
  In order to do so, first notice that
  $s^2S=sS=S$ so, excluding the elementary case in which $S=\{0\}$, we have that $s^2\neq 0$.  Consequently from
  $
  sus=ss\neq 0,
  $
  we deduce that $us=s$.  Since any element $t$ in $S$ may be written in the form $t=sx$, for some $x$ in $S$, we have
  $$
  ut=usx = sx = t,
  $$
  so $u$ is a left identity for $S$.  In addition, given $t$ in $S$, we have that
  $$
  tux = tx, \for x\in S,
  $$
  so $tu=t$ by right reductivity.  Therefore $u$ is also a right identity, hence a (two-sided) identity and we see that
$S$ is unital.

In order to prove that $s$ is invertible, let $t$ be such that $st=u$.  Then
  $$
  sts = us = s = su,
  $$
  and by $0$-left cancellativity we get that $ts=u$, so $s$ is invertible and $s\inv =t$.  \endProof

The following definition introduces one of the main concepts studied in this work.

\definition
  The \emph {inverse hull} of a $0$-left cancellative semigroup $S$, henceforth denoted by $\hull $, is the inverse
subsemigroup of $\I (S')$ generated by the set $\{\theta _s : s \in S\}$. Thus, in the terminology of \ref
{GenSgANdConstr.i} we have
  \# $\hull $; Inverse hull of $S$;
  $$
  \hull = \I (S',\theta ).
  $$

The reader should compare the above with the notion of \emph {inverse hull} considered in \cite [Section~1.9]{CP} and
\cite {Cherubini}.

The collection of $\theta $-constructible subsets of $S'$ is of special importance to us, so we would like to give it a
special notation:

\definition \label DefIdempSLA
  The idempotent {\sla } of $\hull $, which we will tacitly identify with the {\sla } of $\theta $-constructible subsets
of $S'$, will be denoted by $\ehull $.
  Thus, in the terminology of \ref {GenSgANdConstr.iii} we have
  $$
  \ehull = \P (S',\theta ).
  $$
  \# $\ehull $; Set of all $\theta $-constructible subsets of $S'$;

Since the $\theta $-constructible sets may be described by \ref {CharactSL}, it is interesting to have a concrete
description for the maps mentioned there in the special case of the regular representation.

\state Lemma \label InterpPibPif
  Regarding the regular representation $\theta $ of $S$ on $S'$, for every $s$ in $S$, and every $X\subseteq S'$, one
has that
  \iItemize
  \iItem $s\inv [X] = \{y\in S': sy\in X\}$, and
  \iItem $\phantom {\inv }s[X] = \{y\in S': y=sx, \hbox { for some } x\in X\} = sX\setminus \{0\}$.

\Proof
  Left for the reader.
  \endProof

It will be of importance to identify some properties of $0$-left cancellative semigroups that will play a role later.

\state Proposition \label SomeCancelProps Let $S$ be a $0$-left cancellative semigroup.
    \iItemize
    \iItem If $e\in E(S)$ and $s\in S\setminus \{0\}$, then $es\neq 0$, if and only if $es=s$, that is, $s\in
eS\setminus \{0\}$.
    \iItem If $s\in S\setminus \{0\}$, then $s\in sS$ if and only if $se=s$ for a necessarily unique idempotent $e$.
    \iItem If $sS=S$ and $S$ is right reductive, then $S$ is unital and $s$ is invertible.

\Proof For the non-trivial ``only if'' direction of the first item, assume that $es\neq 0$. Then $es=ees\neq 0$ implies
$s=es$ by $0$-left cancellativity.  For the second item, if $sx=s$, then $sxx=sx=s\neq 0$ implies that $x^2=x$ by
$0$-left cancellativity.  Also $sx=s=sy$ implies $x=y$ by $0$-left cancellativity.  This establishes the second item.
The third item is the content of Proposition \ref {sSEqualSImpliesInvert}.
    \endProof

A semigroup $S$ is said to have \emph {right local units} if $S=SE(S)$, that is, for all $s\in S$, there exists $e\in
E(S)$ with $se=s$.  A unital semigroup has right local units for trivial reasons.  If $S$ has right local units, then
$sS=0$ implies that $s=0$.  From Proposition~\ref {SomeCancelProps} we obtain the following corollary.

\state Corollary \label LocalUnits Let $S$ be a $0$-left cancellative semigroup.  Then $S$ has right local units if and
only if $s\in sS$ for all $s\in S$.

\state Proposition \label OrthogIdems Let $S$ be a right reductive, $0$-left cancellative semigroup and suppose that $e$
and $f$ are idempotent elements of $S$ with $e\neq f$.  Then $ef=0$.

\Proof This is obvious if $e$ or $f$ is $0$.  So assume that $e\neq 0\neq f$.  If $ef\neq 0$, then $ef=f$ by
Proposition~\ref {SomeCancelProps.i}.  But then $fef=ff=f\neq 0$ and so $fe\neq 0$.  Therefore, $fe=e$ by
Proposition~\ref {SomeCancelProps.i}.  But then $eS=fS$ and so by Proposition~\ref {SomeCancelProps.i} we obtain that
$ex=x=fx$ for all $x\in eS=fS$ and $ex=0=fx$ for all $x\notin eS=fS$.  Thus $e=f$ as $S$ is right reductive.
    \endProof

\definition If $S$ is a $0$-left cancellative, right reductive semigroup with right local units, then for $s\in
S\setminus \{0\}$, we denote by
  $s^+$
  \# $s^+$; Right local unit for $s$;
  the unique idempotent with $ss^+=s$.  If $S$ is unital, then $s^+=1$.

   We can associate to a left cancellative category $C$ (i.e., a category of monics) a semigroup $S(C)$ by letting
$S(C)$ consist of the arrows of $C$ together with a zero element $0$.  Products that are undefined in $C$ are made zero
in $S(C)$ and the remaining products are as in $C$. It is straightforward to check that $S(C)$ is $0$-left cancellative,
right reductive and has right local units. If $f\colon c\to d$ is an arrow of $C$, then $f^+=1_c$.  The case when $C$ is
the category associated to a higher rank graph \cite {KumjianPask}, will be considered later in this paper.

The special case of $0$-left cancellative semigroups of the form $S(C)$, with $C$ a left cancellative category, is
considered in detail by Spielberg in~\cite {SpielbergB} (which was posted after our announcement of the results of this
paper~\cite {announce}, but was done independently of our work).  An advantage of our more general framework is that if
one has an ideal in a left cancellative category, then the quotient category by this ideal need not be left
cancellative, but factoring $S(C)$ by the ideal will result in a $0$-left cancellative semigroup.  This type of ideal
construction, applied to categories of paths, lets one go from shifts of finite type to arbitrary shifts by forbidding
patterns.

An ideal in a semigroup $S$ is a non-empty subset $I$ such that $SI\cup IS\subseteq I$.  The Rees quotient $S/I$ is the
quotient of $S$ by the congruence identifying $I$ to a single element (which will be the zero element of $S/I$).  Each
element of $S\setminus I$ forms its own equivalence class.  The class of $0$-left cancellative semigroups is evidently
closed under Rees quotients and taking subsemigroups (containing $0$).

\state Proposition \label LocalUnitsRightRed Suppose that $S$ is $0$-left cancellative, right reductive and has right
local units.
    \iItemize
    \iItem If $s\in S\setminus \{0\}$, then $sx\neq 0$ implies $x\in s^+S$.
    \iItem If $s,t\in S\setminus \{0\}$ and $sx\neq 0\neq tx$, then $s^+=t^+$.
    \iItem If $T$ is a subsemigroup of $S$ containing $0$ that is closed under the unary operation $s\mapsto s^+$, then
$T$ is $0$-left cancellative, right reducitve and has right local units.
    \iItem If $I$ is a proper ideal of $S$, then the Rees quotient $S/I$ is $0$-left cancellative, right reductive and
has right local units.

\Proof For the first item, $0\neq sx=ss^+x$ implies $s^+x\neq 0$ and hence $x\in s^+S$ by Proposition~\ref
{SomeCancelProps}.  For the second item, we have that $x\in s^+S\cap t^+S\subseteq s^+t^+S$.  It follows that
$s^+t^+\neq 0$ and so $s^+=t^+$ by Proposition~\ref {OrthogIdems}.  Suppose that $T$ is a subsemigroup closed under the
unary operation.  Then it is obviously $0$-left cancellative and has right local units.  Suppose that $s\in T\setminus
\{0\}$ and $t\in T$ with $sx=tx$ for all $x\in T$.  Then $s=ss^+=ts^+$.  We conclude that $t\neq 0$ and so $s^+=t^+$ by
the second item. Therefore, $s=ts^+=tt^+=t$.  Thus $T$ is right reductive.  The final item is proved in the same way as
the previous one; we omit the details.
    \endProof

Note that in general a subsemigroup or Rees quotient of a right reductive semigroup need not be right reductive so the
right local units play a key role in the last two items.

\section Categorical at zero semigroups

At this point we would like to remark that Definition \ref {GenSgANdConstr.ii}, especially when applied to the regular
representation $\theta $, is motivated by Li's use of the term \emph {constructible} in \cite {Li}.  However we should
notice that, contrary to the situation treated in \cite {Li}, our $\theta $-constructible subsets are not necessarily
related to right ideals.

On the positive side, under special conditions on $S$ we shall soon prove that any $\theta $-constructible set is the
nonzero part of a right ideal in $S$.

\definition \label CatAtZero
  (\cite {Munn}) Let $S$ be a semigroup with zero.
  We will say that $S$ is \emph {categorical at zero} if, for every $r,s,t\in S$, one has that
  $$
  rs\neq 0, \hbox { and } st\neq 0 \Imply rst\neq 0.
  $$

When applied to unital semigroups, the above concept is not very interesting since it reduces to the absence of zero
divisors.  In fact, if one is allowed to take $s=1$, then the above condition would read
  $$
  r\neq 0, \hbox { and } t\neq 0 \Imply rt\neq 0.
  \equationmark NoZeroDivisors
  $$

The reason for the terminology is that if $C$ is a category, then the semigroup $S(C)$ constructed above is categorical
at zero.  However not all categorical at zero semigroups have the form $S(C)$, as illustrated by semigroups arising from
Markov subshifts to be introduced below.

Recall that a subset $R\subseteq S$ is said to be a \emph {right ideal}\/ when $Rs\subseteq R$, for every $s$ in $S$.
Notice that right ideals always contain the zero element.

Generalizing our notation $S'$ introduced in \ref {IntroSPrime}, for each $X\subseteq S$, let us put
  $$
  X'=X\setminus \{0\}.
  $$

The following result employs the square bracket notation defined in \ref {Conjuga}.

\state Lemma \label StrAssocResult
  Given any $s$ in $S$, and any right ideal $R\subseteq S$, one has that
  \iItemize
  \iItem $s[R']\cup \,\{0\}$ is a right ideal in $S$, and
  \iItem $s\inv [R']\cup \,\{0\}$ is a right ideal in $S$, provided $S$ is categorical at zero.

\Proof
  Observing that
  $$
  s[R']\cup \,\{0\}\={InterpPibPif}
  (s{R'}\setminus \{0\})\cup \,\{0\}=
  s{R'}\cup \,\{0\}=
  sR,
  $$
  the proof of (i) is clear.

Regarding (ii), pick $x$ in $s\inv [R']\cup \,\{0\}$, and $y$ in $S$.  We must then prove that
  $$
  xy\in s\inv [R']\cup \,\{0\}.
  $$
  If $xy=0$, there is nothing to be done, so we suppose that $xy\neq 0$. Consequently $x\neq 0$, and then we see that
$x\in s\inv [R']$, so $sx\in R'$, by \ref {InterpPibPif}.

In particular $sx\neq 0$, so $sxy\neq 0$ because $S$ is categorical at zero.  Moreover, since $sx\in R$, we also have
that $sxy\in R$, and consequently $sxy\in R'$, which implies that $xy\in s\inv [R']$, as desired.
  \endProof

To see that being categorical at zero is important in \ref {StrAssocResult.ii}, let $S$ be a semigroup not possessing
this property, and take $s,x,y\in S$ with $sx$ and $xy$ nonzero, but $sxy=0$.

Considering $S$ as a right ideal in itself, notice that $sx\in S'$, so $x\in s\inv [S']$.  However $xy$ is not in $s\inv
[S']\cup \,\{0\}$, because neither is $xy=0$, nor is $sxy$ in $S'$.  So $s\inv [S']\cup \,\{0\}$ is not a right ideal.

\state Proposition \label ConstrIsIdeal
  If $S$ is a $0$-left cancellative semigroup which is categorical at zero, then every $\theta $-constructible subset of
$S'$ coincides with the set of nonzero elements of some right ideal of $S$.

\Proof
  Letting
  \def \Rr {\goth R}
  $$
  \Rr = \{R': R\hbox { is a right ideal of } S\},
  $$
  we claim that $\P (S',\theta )\subseteq \Rr $.  In view of \ref {CharactSL}, in order to prove this claim all we need
to do is to show that $\Rr $ contains every $E^\theta _s$, and that $\Rr $ is invariant under the two maps referred to
in \ref {CharactSL}.

Since
  $$
  E^\theta _s= (sS)',
  $$
  we see that $E^\theta _s$ lies in $\Rr $.  Moreover, given any $s$ in $S$, and any right ideal $R\subseteq S$, we have
by \ref {StrAssocResult.ii} that $T:= s\inv [R']\cup \,\{0\}$ is a right ideal.  Observing that $0\not \in s\inv [R']$,
we have that
  $$
  s\inv [R'] = T\setminus \{0\}= T'\in \Rr ,
  $$
  so we see that $\Rr $ is invariant under the first map referred to in \ref {CharactSL}. Regarding the second one, let
us again be given a right ideal $R\subseteq S$.  We then have by \ref {StrAssocResult.i} that
  $T:= s[R']\cup \,\{0\}$ is a right ideal.  Therefore
  $$
  s[R']\setminus \{0\}= T\setminus \{0\}= T'\in \Rr .
  $$
  This proves that $\Rr $ is invariant under the second map referred to in \ref {CharactSL}, and hence our claim that
$\P (S',\theta )\subseteq \Rr $ is verified, from where the statement follows.
  \endProof

Recalling our discussion right before the statement of \ref {ConstrIsIdeal}, when we observed that $s\inv [S']\cup
\,\{0\}$ is not a right ideal, we see that the $\theta $-constructible set
  $
  F^\theta _s = s\inv [S']
  $
  is not the nonzero part of a right ideal.  This says that the hypothesis that $S$ is categorical at zero in \ref
{ConstrIsIdeal} cannot be removed.

Since many of the examples we have in mind involve semigroups which are not categorical at zero, we will unfortunately
not be in a position to benefit from Proposition \ref {ConstrIsIdeal}.  It is given above mostly for the purpose of
comparing our work with Li's \cite {Li} study of C*-algebras of semigroups.  However, in the case of the semigroup
associated to a left cancellative category with least common multiples, it will be useful.

\state Proposition \label CategoricalLocalUnits Let $S$ be a categorical at zero, $0$-left cancellative, right reductive
semigroup with right local units.  Let $s\in S\setminus \{0\}$. Then $sx\neq 0$ if and only if $x\in s^+S\setminus
\{0\}$.  In other words, $F^{\theta }_s=s^+S\setminus \{0\}$ and $E^{\theta }_s=sS\setminus \{0\}$.

\Proof Note that $ss^+=s\neq 0$ and so $sx=ss^+x$ is non-zero if and only $s^+x\neq 0$.
    \endProof

\section Least common multiples

We now wish to introduce a class of semigroups possessing a property inspired by the notion of least common multiples
from arithmetic.  In order to do so we need to consider the question of divisibility.

\definition
  Given $s$ and $t$ in a semigroup $S$, we will say that $s$
  \emph {divides}
  $t$, in symbols
  $$
  s\| t,
  $$
  \# $s\| t$; $s$ divides $t$;
  or that $t$ is a \emph {multiple} of $s$, when either $s=t$, or there is some $u$ in $S$ such that
  $su=t$.  In other words, $s$ divides $t$ if and only if
  $$
  t\in \{s\}\cup sS.
  $$

This should actually be called \emph {left-division}, since one could alternatively define \emph {right-divisibility}
upon replacing the above expression ``$su=t$" with ``$us=t$".  However we will not have any use for right-division, and
hence we may safely use the term \emph {division} to mean \emph {left-division}.

We observe that division is a reflexive and transitive relation, so it may be seen as a (not necessarily anti-symmetric)
order relation upon defining ``$\leq $" by
  $$
  s\leq t\ \Leftrightarrow \ s\|t.
  \equationmark OrderDivision
  $$ This is the dual of Green's quasi-order $\leq _{\mathcal R}$, which is usually considered in semigroup theory.

  Since any $s$ in $S$ divides $0$, one has that $0$ is the maximum element of $S$.  Should $S$ be a unital semigroup,
with unit denoted $1_S$, then $1_S$ is a minimum element, a property shared by all other invertible elements of $S$.

Incidentally, when $S$ is unital, or more generally has right local units, we may define division is a slightly simpler
way since
  $$
  s\|t \iff t\in sS.
  $$

For the strict purpose of simplifying the description of the division relation, regardless of whether or not $S$ is
unital, we shall sometimes employ the unitized semigroup
  $$
  \tS := S\cup \{1\},
  $$
  \# $\textbackslash tS$; Unitized semigroup $S\cup \{1\}$;
  where $1$ is any element not belonging to $S$, made to act like a unit for $S$.  For every $s$ and $t$ in $S$ we
therefore have that
  $$
  s\|t \iff \exists u\in \tS , \ su=t.
  \equationmark UnitDivision
  $$

Having enlarged our semigroup, we might as well extend the notion of divisibility:

\definition Given $v$ and $w$ in $\tS $, we will say that $v\|w$ when there exists some $u$ in $\tS $, such that $vu=w$,
i.e., $w\in v\tS $.

Notice that if $v$ and $w$ are in $S$, then the above notion of divisibility coincides with the previous one by \ref
{UnitDivision}.  Analysing the new cases where this extended divisibility may or may not apply, notice that:
  $$
  \matrix {
  \forall w\in \tS , &
  1\|w, \hfill \cr \pilar {20pt}
  \forall v\in \tS , &
  v\|1 \iff v=1.}
  \equationmark NewDivision
  $$

  The introduction of $\tS $ brings with it several pitfalls, not least because $\tS $ might not be $0$-left
cancellative: when $S$ already has a unit, say $1_S$, then in the identity ``$s1_S=s1$'', we are not allowed to left
cancel $s$, since $1_S\neq 1$.  One should therefore exercise extra care when working with $\tS $.

The notion of least common multiples is well studied for unital semigroups.  We could not find much in the literature in
the non-unital setting and a number of subtleties arise.

\definition \label DefLCM
  Let $S$ be a semigroup and let $s,t\in S$.  We will say that an element $r\in S$ is a \emph {least common multiple}
for $s$ and $t$ when
  \iItemize
  \iItem $sS\cap tS = rS$,
  \iItem both $s$ and $t$ divide $r$.

Observe that when $S$ has right local units then $r\in rS$, by \ref {LocalUnits}, and hence condition \ref {DefLCM.i}
trivially implies \ref {DefLCM.ii}, so the former condition alone suffices to define least common multiples.
  However, in a semigroup without right local units it is not true that condition \ref {DefLCM.i} implies \ref
{DefLCM.ii}.
  For example, if $S^2=\{0\}$, then \ref {DefLCM.i} is always satisfied by any $s,t,r$, but \ref {DefLCM.ii} is only
satisfied for $s\neq t$ when $r=0$.
  Notice that if $S^2=\{0\}$ then $s$ and $0$ are least common multiples of $s$ with itself.

  The above example indicates some of the strange things that can happen when $S$ lacks local units and is not right
reductive. Nevertheless, some of the main examples we have in mind, such as \ref {LanguageSemigroup} below, behave very
well with respect to least common multiples even though they do not admit local units.

Regardless of the existence of right local units, observe that when $sS\cap tS=\{0\}$, then $0$ is always a least common
multiple for $s$ and $t$ because $s$ and $t$ always divide $0$.  In addition, notice that condition \ref {DefLCM.ii}
holds if and only if $r\tS \subseteq s\tS \cap t\tS $, and therefore one has that $r$ is a least common multiple for $s$
and $t$ if and only if
  $$
  sS\cap tS = rS \subseteq r\tS \subseteq s\tS \cap t\tS .
  \equationmark sStStilde
  $$

\definition \label DefLCMSgrp
  We shall say that a semigroup $S$ \emph {admits least common multiples} if there exists a least common multiple for
each pair of elements of $S$.

\section Examples

\label ExamplesSection Even though we believe examples are of fundamental importance in any mathematical work, we have
hitherto postponed their presentation to give us time to build the necessary terminology needed to highlight their
relevant properties.

Our first class of examples comes from Language Theory.
  Let $\Lambda $ be any finite or infinite set, henceforth called \emph {the alphabet}, and let $\Lambda ^+$ be the free
semigroup generated by $\Lambda $, namely the set of all finite words in $\Lambda $ of positive length (and hence
excluding the empty word), equipped with the multiplication operation given by concatenation. Incidentally recall that
the free monoid on $\Lambda $ is customarily denoted $\Lambda ^*$; it includes the empty string.

Let $L$ be a \emph {language} on $\Lambda $, namely any nonempty subset of $\Lambda ^+$.  We will furthermore assume
that $L$ is \emph {closed under prefixes and suffixes}, that is, for every $\alpha $ and $\beta $ in $\Lambda ^+$, one
has
  $$
  \alpha \beta \in L \Imply \alpha \in L, \hbox { and } \beta \in L .
  $$
  This is equivalent to $L$ being closed under factors: $\alpha \beta \gamma \in L$ implies $\beta \in L$ for all $\beta
\in \Lambda ^+$ and $\alpha ,\gamma \in \Lambda ^*$.

  Define a multiplication operation on
  $$
  S:=L\cup \,\{0\},
  $$
  where $0$ is any element not belonging to $\Lambda ^+$, by
  $$
  \def \quad { }
  \alpha \cdot \beta =\left \{\matrix {
  \alpha \beta , & \hbox { if $\alpha ,\beta \neq 0$, and $\alpha \beta \in L$}, \cr \pilar {12pt}
  0, & \hbox { otherwise.}\hfill }\right .
  $$

The reader will have no difficulty in proving the following:

\state Proposition \label LanguageSemigroup
  Given any language $L\subseteq \Lambda ^+$, closed under prefixes and suffixes, the above multiplication operation is
associative thus making $S$ a semigroup with zero. Moreover $S$ is $0$-cancellative and admits least common multiples.
There are no idempotent elements in $S$, whence $S$ lacks right local units.

One may also see $S$ as the Rees quotient $\Lambda ^+/I$, where $I=\Lambda ^+\setminus L$.  Because $L$ is closed under
factors, $I$ is obviously an ideal of $\Lambda ^+$ and hence $S$ is a semigroup.  The fact that $S$ is $0$-left
cancellative may also be deduced from the fact that any Rees quotient of a $0$-left cancellative semigroup shares this
property.

Notice that $S$ has no nonzero idempotent elements and hence it cannot have right local units.  In extreme cases $S$
could also fail to be right reductive such as when all words in $L$ have length one.

We may now easily give an example where $S$ is not categorical at zero: take any nonempty alphabet $\Lambda $, let $L$
be the language consisting of all words of length at most two, and let $S=L\cup \{0\}$, as above.  If $a$, $b$ and $c$,
are members of $\Lambda $, we have that $abc=0$, but $ab$ and $bc$ are nonzero, so $S$ is not categorical at zero.
Regarding our discussion after the proof of \ref {ConstrIsIdeal}, observe that $F^\theta _a$ consists of all elements
$x$ of $L$ such that $ax\neq 0$, so that $F^\theta _a$ is precisely the set of all words of length $1$, which is
certainly not the nonzero part of a right ideal in $S$.

One important special case of the above is based on \emph {subshifts}.  Given an alphabet $\Lambda $, as above, consider
the \emph {left shift}, namely the mapping $\sigma \colon \Lambda ^{\bf N}\to \Lambda ^{\bf N}$ given by
  $$
  \sigma (x_1x_2x_3\ldots )=x_2x_3x_4\ldots .
  $$

  A nonempty subset $\X \subseteq \Lambda ^{\bf N}$ is called a
    \emph {subshift}\fn {The term \emph {subshift} is often applied to the map obtained by restricting $\sigma $ to $\X
$.
  Moreover, in the field of symbolic dynamics it is also required that $\X $ be closed in the product topology of
$\Lambda ^{\bf N}$.}
  when it is invariant under $\sigma $ in the sense that $\sigma (\X )\subseteq \X $.

  Given a subshift $\X $, let $L_\X \subseteq \Lambda ^+$ be the \emph {language of} $\X $, namely the set of all finite
words occuring in some infinite word belonging to $\X $. Then $L_\X $ is clearly closed under prefixes and suffixes, and
hence we are back in the conditions of example \ref {LanguageSemigroup}.

The fact that $\X $ is invariant under the left shift is indeed superfluous, as any nonempty subset $\X \subseteq
\Lambda ^{\bf N}$ would lead to the same conclusion. However, languages arising from subshifts have been intensively
studied in the literature, hence the motivation for considering this situation.

\definition \label DefineSubshSgrp
  Given a subshift $\X $, we will denote by $\SX $ the semigroup built from $L_\X $ as in \ref {LanguageSemigroup}.

Given an alphabet $\Lambda $, let us be given a matrix
  $$
  A= \{A_{x,y}\}_{x,y\in \Lambda }
  $$
  such that $A_{x,y}\in \{0,1\}$, for all $x,y\in \Lambda $.  Such a matrix is sometimes called a \emph {transition
matrix}.  Define $\X _A$ to be the set of all infinite words
  $$
  x_1x_2x_3\ldots \in \Lambda ^{\bf N},
  $$
  such that
  $$
  A_{x_i,x_{i+1}}=1, \for i\in {\bf N}.
  $$
  It is easy to see that $\X _A$ is invariant under the left shift hence a subshift.
  This is usually referred to as the \emph {Markov subshift} associated to the transition matrix $A$.

The reader will have no difficulty in checking that $S_{\X _A}$ is categorical at zero for every transition matrix $A$,
but it is easy to exhibit subshifts $\X $ for which $S_\X $ does not share this property.

Markov subshifts may be used to exhibit a semigroup which is categorical at zero but is not of the form $S(C)$, as
hinted at in the paragraph following \ref {NoZeroDivisors}.  In fact a semigroup arising from \ref {LanguageSemigroup}
is never isomorphic to some $S(C)$ since the former has no nonzero idempotent elements while the latter has many, namely
the identity morphism of each object of $C$.

Markov subshifts may indeed be used to produce a semigroup which is categorical at zero and yet is not isomorphic to any
\emph {subsemigroup} of $S(C)$, no matter which category $C$ one takes.  To see this, consider the alphabet $\Lambda
=\{x_1,x_2\}$ and let $A$ be transition matrix
  $$
  A=\pmatrix {1 & 1\cr 1 & 0}.
  $$

Notice that the words
  $x_1x_1$, $x_1x_2$, and $x_2x_1$, belong to the language of $\X _A$, but $x_2x_2$ is forbidden, precisely because $A_
{x_2, x_2}=0$.

Should there exist a category $C$ such that $S$ is a subsemigroup of $S(C)$, the fact that, say, $x_1x_2\neq 0$ would
lead one to believe that $d(x_1)$, namely the \emph {domain} of $x_1$, coincides with $r(x_2)$, the \emph {range} of
$x_2$.
  But then for similar reasons one would have
  $$
  d(x_2) = r(x_1) = d(x_1) = r(x_2),
  $$
  which would imply that $x_2x_2\neq 0$, a contradiction.

Another interesting class of examples is obtained from the quasi-lattice ordered groups of \cite {Nica}, which we would
now like to briefly describe.

Given a group $G$ and a unital subsemigroup $P\subseteq G$, one defines a partial order on $G$ via
  $$
  x\leq y \iff x\inv y \in P.
  $$
  The quasi-lattice condition says that, whenever elements $x$ and $y$ in $G$ admit a common upper bound, namely an
element $z$ in $G$ such that $z\geq x$ and $z\geq y$, then there exists a least common upper bound, usually denoted
$x\vee y$.

Under this situation, consider the semigroup $S = P\cup \{0\}$, obtained by adjoining a zero to $P$.  Then, for every
nonzero $s$ in $S$, i.e., for $s$ in $P$, one has that
  $$
  sS=\{x\in P: x\geq s\}\cup \{0\},
  $$
  so that the multiples of $s$ are precisely the upper bounds of $s$ in $P$, including zero.

  If $t$ is another nonzero element in $S$, one therefore has that
  $s$ and $t$ admit a nonzero common multiple if and only if $s$ and $t$ admit a common upper bound in $P$, in which
case $s\vee t$ is a least common multiple of $s$ and $t$.

On the other hand, when $s$ and $t$ admit no common upper bound, then obviously $s\vee t$ does not exist, but still $s$
and $t$ admit a least common multiple in $S$, namely $0$.

Summarizing our discussion so far we have the following:

\state Proposition
  Let $(G,P)$ be a quasi-lattice ordered group.  Then the $0$-left cancellative semigroup $S:=P\cup \{0\}$ admits least
common multiples.

We give now an example of a $0$-left cancellative semigroup $S$ where, for all $s,t\in S$, there exists $r\in S$ with
$sS\cap tS=rS$, but $S$ fails to admit least common multiples in our sense. Let $S=\{a,b,c,ab,ba,c^2,0\}$ where
$ac=bc=c^2$ and all other non-obvious products are $0$.  In particular, any product of three elements of $S$ is $0$ and
so $S$ is associative.  Then $aS\cap bS=aS\cap cS=bS\cap cS=cS=\{c^2,0\}$ but $c$ is not a common multiple of $a$ and
$b$.  For all other $x\in \{ab,ba,c^2,0\}$, we have $xS=0$.

Another class of examples may be obtained from semigroupoids, as defined in \cite [Section 14]{actions}.  Given a
semigroupoid $\Lambda $, consider the semigroup $S = \Lambda \cup \{0\}$, where $0$ is any element not belonging to
$\Lambda $, with multiplication defined by $f\cdot 0 = 0\cdot f = 0$, for all $f$ in $S$, while for $f$ and $g$ in $S$,
we put
  $$
  f\cdot g = \left \{ \matrix {
  fg, & \hbox {if } (f,g)\in \Lambda ^{(2)}, \cr \pilar {12pt}
  0, & \hbox {otherwise}.\hfill
  }\right .
  $$

\state Proposition \label Semigroupoids
  Given a semigroupoid $\Lambda $, let $S$ be the semigroup constructed above.  Then
  \iItemize
  \iItem $S$ is categorical at zero,
  \iItem if every element of $\Lambda $ is monic \cite [Definition 14.5]{actions}, then $S$ is $0$-left cancellative.

Of course \ref {Semigroupoids.i} is a consequence of the choice of the strong associativity property in \cite {actions}.
Conversely, given a semigroup $S$ with zero, one could let $\Lambda =S\setminus \{0\}$, with partial multiplication
defined on
  $$
  \Lambda ^{(2)} = \{(s, t)\in \Lambda : st \neq 0\}.
  $$
  If $S$ is categorical at zero one may prove that $\Lambda $ satisfies the associativity property of \cite [Section
14]{actions}, but one could alternatively generalize the notion of semigroupoid by assuming a less stringent
associativity axiom.

We should also mention a few other classes of examples which are in fact special cases of some of the above examples
but, given their role in the modern literature, it is perhaps worth singling them out.

Given a small category $\cal C$ in which all arrows are monomorphisms, as already mentioned we may associate to $\cal C$
a semigroup $S({\cal C})$ consisting of the arrows of ${\cal C}$ together with a zero element, where the product is
extended by making it zero whenever not already defined.  Then $S({\cal C})$ is a $0$-left cancellative semigroup with
right local units.  This is of course a special case of example \ref {Semigroupoids}.

Given any directed graph $E$, one may view the collection of all finite paths in $E$ (with or without the vertices
which, if included, could be viewed as paths of length zero) as a semigroupoid (also as a category if the vertices are
included) in which every element is monic and hence one may again build a semigroup as in \ref {Semigroupoids}.  Should
one prefer not to include the vertices, this may also be though of as a special case of example \ref
{LanguageSemigroup}, where the alphabet is taken to be the set of all edges, the language consisting of all finite paths
in $E$.

The situation in the above paragraph could also be applied to any higher rank graph with similar conclusions.

\section Normal Form

\fix Throughout this section we will fix a $0$-left cancellative semigroup $S$ admitting least common multiples.

Given a representation $\pi $ of $S$ on a set $\Omega $, we will now concentrate our attention in giving a concrete
description for the elements of $\I (\Omega ,\pi )$ (see Definition \ref {GenSgANdConstr}), provided $\pi $ satisfies
certain special properties, which we will now describe.

Initially notice that if $s\|r$, then the range of $\pi _r$ is contained in the range of $\pi _s$ because either $r=s$,
or $r=su$, for some $u$ in $S$, in which case
  $
  \pi _r = \pi _s\pi _u.
  $
  \ So, using the notation introduced in \ref {IntroEsFs},
  $$
  E^\pi _r \subseteq E^\pi _s.
  $$ When $r$ is a least common multiple of $s$ and $t$, it then follows that
  $$
  E^\pi _r \subseteq E^\pi _s \cap E^\pi _t.
  $$

\definition \label DefRespect
  A representation $\pi $ of $S$ is said to \emph {respect least common multiples} if, whenever $r$ is a least common
multiple of elements $s$ and $t$ in $S$, one has that
  $
  E^\pi _r = E^\pi _s \cap E^\pi _t.
  $

As an example, notice that the regular representation of $S$, defined in \ref {DefineRegRep}, satisfies the above
condition since the fact that $rS=sS\cap tS$ implies that
  $$
  E^\theta _r =
  rS\setminus \{0\} =
  (sS \cap tS)\setminus \{0\} =
  (sS\setminus \{0\}) \cap (tS\setminus \{0\}) =
  E^\theta _s \cap E^\theta _t.
  \equationmark RanLCM
  $$

\fix From now on we will moreover fix a representation $\pi $ of $S$ on a set $\Omega $, assumed to respect least common
multiples.

Since $\pi $ will be the only representation considered for a while, we will
  use the simplified notations
  $F_s$, $E_s$, $f_s$, and $e_s$.

There is a cannonical way to extend $\pi $ to $\tS $ by setting
  $$
  F_1=E_1=\Omega , \and \pi _1=\id _\Omega .
  $$
  It is evident that $\pi $ remains a multiplicative map after this extension.  Whenever we find it convenient we will
therefore think of $\pi $ as defined on $\tS $ as above.  We will accordingly extend the notations
  $f_s$ and $e_s$ to allow for any $s$ in $\tS $, in the obvious way.

\state Proposition \label ExistLCMinSTilde
  Given $u$ and $v$ in $\tS $, there exists $w$ in $\tS $ such that
  \iItemize
  \iItem $uS\cap vS = wS$,
  \iItem both $u$ and $v$ divide $w$.

\Proof
  When $u$ and $v$ lie in $S$, it is enough to take $w$ to be a (usual) least common multiple of $u$ and $v$.  On the
other hand, if $u=1$, one takes $w=v$, and if $v=1$, one takes $w=u$.  \endProof

Based on the above we may extend the notion of least common multiples to $\tS $, as follows:

\definition \label NewDefLCM
  Given $u$ and $v$ in $\tS $, we will say that an element $w$ in $\tS $ is a least common multiple of $u$ and $v$,
provided \ref {ExistLCMinSTilde.i-ii} hold.  In the exceptional case that $u=v=1$, only $w=1$ will be considered to be a
least common multiple of $u$ and $v$, even though there might be another $w$ in $S$ satisfying \ref
{ExistLCMinSTilde.i-ii}.

It is perhaps interesting to describe the exceptional situation above, where we are arbitrarily prohibiting by hand that
an element of $S$ be considered as a least common multiple of $1$ and itself, even though it would otherwise satisfy all
of the required properties.
  If $w\in S$ is such an element, then
  $$
  wS = 1S\cap 1S=S,
  $$
  so, in case we throw in the assumption that $S$ is right-reductive, we deduce from \ref {sSEqualSImpliesInvert} that
$S$ is unital and $w$ is invertible.  Thus, in hindsight it might not have been such a good idea to add an external unit
to $S$ after all!

On the other hand, when $s$ and $t$ lie in $S$, it is not hard to see that any least common multiple of $s$ and $t$ in
the new sense of \ref {NewDefLCM} must belong to $S$, and hence it must also be a least common multiple in the old sense
of \ref {DefLCM}.

\state Proposition \label ExtensionRespLCM
  Let $\pi $ be a representation of\/ $S$ on a set\/ $\Omega $.  If $\pi $ respects least common multiples then so does
its natural extension to $\tS $.  Precisely, if $u$ and $v$ are elements of $\tS $, and if $w\in \tS $ is a least common
multiple of $u$ and $v$, then
  $
  E _w = E _u \cap E _v.
  $

\Proof
  If $u$ and $v$ lie in $S$, then $w$ is necessarily a least common multiple of $u$ and $v$ in the old sense of \ref
{DefLCM}, so the result follows by hypothesis.

If $u=v=1$, then $w=1$ by
  default\fn {Should we have allowed in \ref {NewDefLCM} that another element $w$ of $S$ be considered a least common
multiple of $1$ and itself, at this point we would be required to prove that $E_w=\Omega $.  This would still be within
reach, as long as we loaded up on our hypotheses, requiring $S$ to be right-reductive and $\pi $ to be essential.  With
all of this we could invoke \ref {sSEqualSImpliesInvert} to deduce that $S$ is unital and $w$ is invertible, and then by
\ref {InvertEssential} we would obtain the desired equality.  In conclusion we believe that adding a little exception to
Definition \ref {NewDefLCM} is a small price to pay for a result with fewer hypotheses and hence wider applicability.},
  and the result follows trivially.

Up to interchanging $u$ and $v$, the last case to be considered is when $u=1$ and $v\in S$.  In this case notice that
$v\|w$, hence $w$ must be in $S$.  Moreover,
  $$
  wS = uS\cap vS = S\cap vS = vS.
  $$

  In an unforeseen twist of fate,
  if follows from the above that $w$ is a least common multiple of $v$ and itself.  So, by hypothesis
  $$
  E _w = E _v \cap E _v = E _v,
  $$
  whence
  $$
  E _w = E _v = \Omega \cap E _v = E_1 \cap E _v = E_u \cap E _v,
  $$
  concluding the proof.  \endProof

\state Lemma \label TudoTil
  Let $\pi $ be a representation of $S$ respecting least common multiples, and let $w\in \tS $ be a least common
multiple for given elements $u$ and $v$ in $\tS $.  Using \ref {UnitDivision} to write $w = ux = vy$, with $x,y\in \tS
$, one has that
  \iItemize
  \iItem $e _ue _v=e _w$,
  \iItem $\pi _u\inv \pi _v= \pi _xf _w\pi _y\inv $.

\Proof
  Since $\pi $ respects least common multiples, even at the level of $\tS $ by \ref {ExtensionRespLCM}, the hypothesis
gives $E _u\cap E _v=E _w$, from which (i) follows.  Regarding (ii) we have
  $$
  \pi _u\inv \pi _v=
  \pi _u\inv \pi _u\pi _u\inv \pi _v\pi _v\inv \pi _v=
  \pi _u\inv e _ue _v\pi _v\explica {(i)}=
  \pi _u\inv e _w\pi _v=
  \pi _u\inv \pi _w\pi _w\inv \pi _v\quebra =
  \pi _u\inv \pi _u\pi _x\pi _y\inv \pi _v\inv \pi _v=
  f _{u} \pi _x\pi _y\inv f _{v} \={Covar.ii}
  \pi _xf _{ux} f _{vy} \pi _y\inv =
  \pi _xf _{w} \pi _y\inv .
  $$

As the careful reader may have noticed, we are using \ref {Covar.ii} above for the extended representation, a result
that can be proved without any difficulty since no assumption was made regarding faithfulness of the representation.
  \endProof

\definition
  Given a representation $\pi $ of $S$, and given any nonempty finite subset $\Lambda \subseteq \tS $, we will let
  $$
  F^\pi _\Lambda = \bigcap _{u\in \Lambda } F^\pi _u, \and
  f^\pi _\Lambda = \prod _{u\in \Lambda }f^\pi _u.
  $$
  \# $F^\pi _\Lambda $; Intersection of $F^\pi _u$, for $u$ in $\Lambda $;
  \# $f^\pi _\Lambda $; Product of $f^\pi _u$, for $u$ in $\Lambda $;

When there is only one representation of $S$ in sight, as in the present moment, we will drop the superscripts and use
the simplified notations
  $F_\Lambda $ and
  $f_\Lambda $.

We should remark that, since each $f _s$ is the identity map on $F _s$, one has that $f _\Lambda $ is the identity map
on $F _\Lambda $.

Also notice that, since $f_1=\id _\Omega $, the presence of $1$ in $\Lambda $ has no effect in the sense that $f_\Lambda
= f_{\Lambda \cup \{1\}}$, for every $\Lambda $.  Thus, whenever convenient we may assume that $1\in \Lambda $.

As already indicated we are interested in obtaining a description of the inverse semigroup $\I (\Omega ,\pi )$.  In that
respect it is interesting to observe that most elements of the form
  $f_\Lambda $ belong to $\I (\Omega , \pi )$, but there is one exception, namely when $\Lambda =\{1\}$.  In this case
we have
  $$
  f_{\{1\}} = \id _\Omega ,
  $$
  which may or may not lie in $\I (\Omega , \pi )$.  However, when $\Lambda \cap S\neq \emptyset $, then surely
  $$
  f_\Lambda \in \I (\Omega , \pi ).
  \equationmark fLambdaInHull
  $$

\state Lemma \label AContadoProduto
  Let $u_1,v_1, u_2,v_2\in S$, and let $\Lambda _1$ and $\Lambda _2$ be nonempty finite subsets of $\tS $.  Let $w$ be a
least common multiple of $v_1$ and $u_2$, and write
  $
  w=v_1x=u_2y,
  $
  for suitable $x,y\in \tS $.  Then
  $$
  (\pi _{u_1} f_{\Lambda _1} \pi _{v_1}\inv )(\pi _{u_2} f_{\Lambda _2} \pi _{v_2}\inv ) = \pi _uf_\Lambda \pi _v\inv ,
  $$
  where
  $u= u_1x$, $v=v_2y$, and $\Lambda = \Lambda _1x\cup \{w\}\cup \Lambda _2y$.

\Proof
  We have
  $$
  \pi _{u_1} f_{\Lambda _1} \pi _{v_1}\inv \pi _{u_2} f_{\Lambda _2} \pi _{v_2}\inv \={TudoTil}
  \pi _{u_1} f_{\Lambda _1} \pi _{x}f_w\pi _{y}\inv f_{\Lambda _2} \pi _{v_2}\inv \={Covar.ii}
  $$
  $$
  =
  \pi _{u_1} \pi _{x}f_{\Lambda _1x}f_wf_{\Lambda _2y} \pi _{y}\inv \pi _{v_2}\inv =
  \pi _{u_1x}f_{\Lambda _1x\cup \{w\}\cup \Lambda _2y} \pi _{v_2y}\inv .
  \closeProof
  $$
  \endProof

We should remark that, whenever we are looking at a term of the form
  $\pi _uf_\Lambda \pi _v\inv $, we may assume that $u,v\in \Lambda $, because
  $$
  \pi _uf_\Lambda \pi _v\inv =
  \pi _u\pi _u\inv \pi _uf_\Lambda \pi _v\inv \pi _v\pi _v\inv \quebra =
  \pi _uf_u f_\Lambda f_v \pi _v\inv =
  \pi _uf_{\{u\}\cup \Lambda \cup \{v\}} \pi _v\inv ,
  \equationmark LambdaHasu
  $$
  so $\Lambda $ may be replaced by $\{u\}\cup \Lambda \cup \{v\}$ without altering the above term.  Moreover, as in \ref
{fLambdaInHull}, observe that if $\Lambda \cap S\neq \emptyset $, then
  $$
  \pi _uf_\Lambda \pi _v\inv \in \I (\Omega , \pi ).
  $$

The following is the promissed concrete description of the elements of $\I (\Omega ,\pi )$.

\state Theorem \label GenFormHull
  Let $S$ be a $0$-left cancellative semigroup admitting least common multiples.
  Also let $\pi $ be a representation of\/ $S$ on a set $\Omega $, assumed to respect least common multiples.
  Then
  $$
  \I (\Omega ,\pi ) =
  \big \{\pi _uf_\Lambda \pi _v\inv :
  \Lambda \subseteq \tS \hbox { is finite, } \Lambda \cap S\neq \emptyset , \hbox { and } u,v\in \Lambda
  \big \}.
  $$

\Proof
  Let us temporarily denote the set appearing in the right hand side above by $\mathcal J$, observing that $\mathcal J
\subseteq \I (\Omega ,\pi )$, as already noted.

We next claim that $\mathcal J$ is an inverse subsemigroup of $\I (\Omega ,\pi )$.  In order to prove it, observe first
that $\mathcal J$ clearly contains the inverse of its elements, so we just need to check that $\mathcal J$ is closed
under multiplication.  Given two elements of $\mathcal J$, say
  $$
  \pi _{u_1} f_{\Lambda _1} \pi _{v_1}\inv \and \pi _{u_2} f_{\Lambda _2} \pi _{v_2}\inv ,
  $$
  we have by \ref {ExistLCMinSTilde} that there exists a least common multiple for $v_1$ and $u_2$, say $w$.  We may
then write
  $
  w=v_1 x=u_2 y,
  $
  with $x,y\in \tS $, and then by \ref {AContadoProduto} we have
  $$
  (\pi _{u_1} f_{\Lambda _1} \pi _{v_1}\inv )(\pi _{u_2} f_{\Lambda _2} \pi _{v_2}\inv ) =
  \pi _{u} f_{\Lambda } \pi _{v}\inv ,
  \equationmark ProdBoldT
  $$
  where $u=u_1x$, $v=v_2y$, and
  $
  \Lambda =\Lambda _1x\cup \{w\}\cup \Lambda _2y.
  $
  So \ref {ProdBoldT} indeed represents an element in $\mathcal J$, thus proving that $\mathcal J$ is an inverse
semigroup as claimed.

Given $s$ in $S$, we have
  $$
  \pi _s = \pi _sf_{\{s, 1\}}\pi _1\inv \in \mathcal J,
  $$
  whence $\mathcal J$ contains the inverse semigroup generated by the $\pi _s$, namely $\I (\Omega ,\pi )$.
  \endProof

With this we may describe the constructible sets in a more concrete way than done in \ref {CharactSL}.

\state Proposition
  Under the assumptions of \ref {GenFormHull}, the $\pi $-constructible subsets of\/ $\Omega $ are precisely the sets of
the form
  $$
  X= \pi _u(F_\Lambda ) ,
  $$
  where $\Lambda \subseteq \tS $ is a finite subset, $\Lambda \cap S\neq \emptyset $, and $u\in \Lambda $.

\Proof
  We leave it for the reader to check that all sets of the above form are $\pi $-constructible, and let us instead show
that every $\pi $-constructible set $X$ is of the above form.

Given any such $X$, we have that $\id _X$ is an idempotent element of $\I (\Omega ,\pi )$, so there exists some $\alpha
$ in $\I (\Omega ,\pi )$ such that $\id _X=\alpha \alpha \inv $.  By \ref {GenFormHull} we may write $\alpha =\pi
_uf_\Lambda \pi _v\inv $, where
  $\Lambda $ is a finite subset of $\tS $, with $\Lambda \cap S\neq \emptyset $, and
  $u,v\in \Lambda $.
  It follows that
  $$
  \id _X =
  \alpha \alpha \inv =
  (\pi _uf_\Lambda \pi _v\inv )(\pi _vf_\Lambda \pi _u\inv ) =
  \pi _uf_\Lambda f_vf_\Lambda \pi _u\inv \quebra =
  \pi _uf_\Lambda \pi _u\inv =
  \pi _u\id _{F_\Lambda } \pi _u\inv = \id _{\pi _u(F_\Lambda )},
  $$
  so
  $
  X= \pi _u(F_\Lambda ),
  $
  as desired.
  \endProof

Recalling that the regular representation of $S$ respects least common multiples, our last two results apply to give:

\state Corollary \label FormOfHull
  Let $S$ be a $0$-left cancellative semigroup admitting least common multiples.  Then
  $$
  \hull =
  \big \{\theta _uf_\Lambda \theta _v\inv :
  \Lambda \subseteq \tS \hbox { is finite, } \Lambda \cap S\neq \emptyset , \hbox { and } u,v\in \Lambda
  \big \},
  $$
  and
  $$
  \ehull =
  \big \{uF_{\Lambda } ,\
  \Lambda \subseteq \tS \hbox { is finite, } \Lambda \cap S\neq \emptyset , \hbox { and } u\in \Lambda
  \big \}.
  $$

We now plan to use the above result to describe the order relation of $\ehull $ in an especially useful way.  Attempting
to motivate what is to come, let $\Lambda $ be a finite subset of $\tS $, with $\Lambda \cap S\neq \emptyset $, and let
$u\in \Lambda $, so that $uF_\Lambda $ is a general element of $\ehull $ by \ref {FormOfHull}.  Given any $x$ in $S$,
the reader is invited to check that
  $$
  uxF_{\Lambda x} \subseteq uF_\Lambda .
  $$
  In addition, if $\Delta $ is any finite subset of $\tS $ such that $\Lambda x\subseteq \Delta $, then clearly
$F_{\Delta } \subseteq F_{\Lambda x}$, so
  $$
  uxF_\Delta \subseteq uF_\Lambda .
  \equationmark StandardSubset
  $$

  We will next prove that the above example is the most general situation in which a member of $\ehull $ is contained in
another one.

\state Proposition \label FormOfSubSets
  Given a finite subset $\Lambda \subseteq \tS $, with $\Lambda \cap S\neq \emptyset $, and given $u$ in $\Lambda $,
suppose that $X$ is a $\theta $-constructible set such that
  $$
  X \subseteq uF_{\Lambda }.
  $$
  Then there is some $x$ in $\tS $, and a finite subset $\Delta \subseteq \tS $, with $\Lambda x\subseteq \Delta $, such
that
  $$
  X = uxF_{\Delta }.
  $$

\Proof
  Using \ref {FormOfHull}, write $X=vF_\Gamma $, with $v\in \Gamma \subseteq \tS $, and $\Gamma \cap S\neq \emptyset $.
  Observe that the sets $vF_{\Gamma }$ and $uF_{\Lambda }$ are respectively the ranges of the idempotent elements
  $
  \theta _vf_\Gamma \theta _v\inv
  $
  and
  $
  \theta _uf_\Lambda \theta _u\inv .
  $ By hypothesis we then have that
  $$
  \id _X =
  \id _{uF_{\Lambda } \cap X} =
  \id _{uF_{\Lambda }} \ \id _X =
  (\theta _uf_\Lambda \theta _u\inv )(\theta _vf_\Gamma \theta _v\inv ).
  $$

  Let $w$ be a least common multiple of $u$ and $v$, and write $w=ux=vy$, for suitable elements $x$ and $y$ in $\tS $.
By \ref {AContadoProduto} the above product turns out to be
  $$
  \theta _{ux}f_{\Lambda x\cup \{w\}\cup \Gamma y}\theta _{vy}\inv =
  \theta _wf_{\Lambda x\cup \{w\}\cup \Gamma y}\theta _w\inv =
  \theta _wf_{\Delta }\theta _w\inv ,
  $$
  where $\Delta =\Lambda x\cup \{w\}\cup \Gamma y$.  We then conclude that
  $$
  X = wF_{\Delta } = uxF_{\Delta }.
  \closeProof
  $$
  \endProof

In the case of semigroups having right local units we may give a slightly more precise description for $\I (\Omega ,\pi
)$.

\state Proposition \label LocUnitsFormOfHull
  Let $S$ be a $0$-left cancellative semigroup admitting least common multiples and right local units.
  Also let $\pi $ be a representation of\/ $S$ on a set $\Omega $, assumed to respect least common multiples.  Then any
nonzero element in $\I (\Omega , \pi )$ may be written as
  $$
  \pi _sf_\Lambda \pi _t\inv ,
  $$
  where $\Lambda $ is a nonempty finite subset of $S$, and $s,t\in \Lambda $.  If moreover $S$ is right-reductive, one
may also assume that $s^+=t^+$, and that $\Lambda \subseteq Ss^+$.

\Proof
  Given any nonzero element $g\in \hull $, use \ref {GenFormHull} to write
  $$
  g=\pi _uf_\Lambda \pi _v\inv ,
  $$
  where
  $\Lambda \subseteq \tS $ is finite, $\Lambda \cap S\neq \emptyset $, and $u,v\in \Lambda $.

We then claim that we may assume that $u$ lies in $S$.  In order to see this, pick any $s$ in $\Lambda \cap S$, and
recall that there is an idempotent element $s^+$ of $S$ such that $s=ss^+$.  Then $\pi _{s^+}$ is an idempotent element
in $\I (\Omega )$, and in particular $\pi (s^+)=\pi (s^+)\inv $.  So
  $$
  f_s =
  \pi _s\inv \pi _s =
  \pi _{ss^+}\inv \pi _s =
  \pi _{s^+}\inv \pi _s\inv \pi _s = \pi _{s^+}f_s,
  \equationmark FsHasPisPlus
  $$
  and consequently
  $$
  g=\pi _uf_\Lambda \pi _v\inv =
  \pi _uf_sf_\Lambda \pi _v\inv =
  \pi _u\pi _{s^+}f_sf_\Lambda \pi _v\inv =
  \pi _{us^+}f_\Lambda \pi _v\inv .
  $$

  Noticing that $us^+\in S$, and that the argument presented in \ref {LambdaHasu} allows us to assume that $us^+\in
\Lambda $, the claim is proven. In an entirely similar way one checks that $v$ may also be taken in $S$.

We will therefore assume that $u, v\in S$, so that
  $$
  u,v\in \Lambda \cap S=:\Lambda '.
  $$
  The only difference between $\Lambda $ and $\Lambda '$, if any, is that $1$ might be in the former but not in the
latter.  In any case it is clear that $f_\Lambda ={f_\Lambda '}$, so
  $$
  g=\pi _uf_\Lambda \pi _v\inv = \pi _uf_{\Lambda '}\pi _v\inv ,
  $$
  proving the first part of the statement.  To address the last part let us suppose from now on that $S$ is
right-reductive.

Observing that $u=uu^+$, and that $\pi _{u^+}$ is and idempotent element, and hence commutes with $f_\Lambda $, we have
  $$
  g = \pi _uf_\Lambda \pi _v\inv =
  \pi _u\pi _{u^+}f_\Lambda \pi _{v^+}\pi _v\inv =
  \pi _uf_\Lambda \pi _{u^+}\pi _{v^+}\pi _v\inv =
  \pi _uf_\Lambda \pi _{(u^+v^+)}\pi _v\inv .
  $$

We then must have that $u^+=v^+$ since otherwise \ref {OrthogIdems} gives $u^+v^+=0$, and we would deduce from the above
that $g=0$.

Given any $t\in \Lambda $, notice that
  $$
  \pi _uf_\Lambda = \pi _u\pi _{u^+}f_tf_\Lambda \={FsHasPisPlus}
  \pi _u\pi _{u^+}\pi _{t^+}f_tf_\Lambda =
  \pi _u\pi _{(u^+t^+)}f_tf_\Lambda .
  $$
  Should $t^+$ not coincide with $u^+$, the above would again imply that $g=0$, so necessarily $t^+=u^+$, and hence
  $$
  t=tt^+=tu^+\in Su^+,
  $$
  thus proving that $\Lambda \subseteq Su^+$.  This concludes the proof.  \endProof

As before we may also describe $\pi $-constructible sets based on \ref {LocUnitsFormOfHull}.

\state Corollary \label LocUnitsFormOfSlatice
  Under the conditions of \ref {LocUnitsFormOfHull}, any $\pi $ constructible subset of\/ $\Omega $ may be written as
  $$
  \pi _s(F_\Lambda ),
  $$
  where $\Lambda $ is a finite subset of $S$, and $s\in \Lambda $.  If moreover $S$ is right-reductive, one may also
assume that that $\Lambda \subseteq Ss^+$.

\Proof
  Left for the reader.  \endProof

Let us say that $g$ dominates $f$ in an inverse semigroup if $g\geq f$.  In the special case in which $S$ is also
$0$-right cancellative we have:

\state Proposition \label ZeroEUnitary
  Let $S$ be a
  $0$-cancellative\fn
    {Recall that a semigroup is $0$-cancellative when it is both $0$-left cancellative and $0$-right cancellative.}
  semigroup admitting least common multiples. Then $\hull $ is a
  $0$-$E$-unitary\fn
    {An inverse semigroup is called $0$-$E$-unitary, or $E^*$-unitary, if whenever an element $g$ dominates a nonzero
idempotent, then $g$ itself is idempotent.}
  inverse semigroup.  Conversely, if $S$ is $0$-left cancellative and has right local units, then $\hull $
$0$-$E$-unitary implies that $S$ is $0$-cancellative.

\Proof
  Pick any $g\in \hull $, and use \ref {FormOfHull} to write
  $$
  g=\theta _uf_\Lambda \theta _v\inv ,
  $$
  where $\Lambda \subseteq \tS $ is finite, $\Lambda \cap S\neq \emptyset $, and $u,v\in \Lambda $, and suppose that $g$
dominates a nonzero idempotent.  Since the idempotent elements in $\I (S')$ are identity functions on their domains, it
follows that $g$ admits a fixed point, say $p\in S'$.

Notice that $p$ lies in the domain of $g$, which is a subset of $E_v$, so we may write $p=vs$, for some $s$ in $S$.  We
then have
  $$
  vs=p=g(p) = us \neq 0.
  $$

We now wish to conclude from the above that $u=v$, which will in turn imply that $g$ is idempotent and the proof will be
finished.  In case both $u$ and $v$ lie in $S$, the desired conclusion that $u=v$ clearly follows from the fact that $S$
is $0$-right cancellative.  If both $u$ and $v$ lie in $\tS \setminus S = \{1\}$, then $u=v$ for obvious reasons.  We
must therefore deal with the remaining situation in which one of $u$ and $v$ lie in $S$, while the other coincides $1$.
By symmetry we will suppose, without loss of generality, that $u$ is in $S$ and $v=1$.  It then follows that $s=us$,
whence also $us=u^2s$, so $0$-right cancellativity implies that $u^2=u$.  Therefore $\theta _u$ is idempotent, and hence
so is
  $\theta _uf_\Lambda =g$.

  Suppose that $S$ is $0$-left cancellative and has right local units.  Assume that $\hull $ is $0$-$E$-unitary and that
$su=tu\neq 0$ with $s,t,u\in S$.  We compute that
  $$
  \theta _s\theta _t\inv \theta _{tu}\theta _{tu}\inv =
  \theta _s\theta _t\inv \theta _t\theta _u\theta _u\inv \theta _t\inv =
  \theta _{su}\theta _{tu}\inv ,
  $$
 which is idempotent and non-zero as $su=tu\neq 0$. Thus $\theta _s\theta _t\inv $ is idempotent.  Let $e$ be an
idempotent with $te=t$.  From $0\neq tu=teu$, we obtain $u=eu$ and hence $0\neq su=seu$ implies that $se\neq 0$ and so
$se=se^2$ implies $s=se$.  Thus $\theta _s\theta _t\inv (te)=se=s$.  But by idempotence, $\theta _s\theta _t\inv
(te)=te=t$.  Thus $s=t$ and so $S$ is right $0$-cancellative.
  \endProof

The description of a given element of $\hull $ in the form $\theta _uf_\Lambda \theta _v\inv $, as in \ref {FormOfHull},
is far from unique and, in fact, it might not be easy to find a unique representation.  However, should an element of
$\hull $ posses two distinct representations of the above form, certain relations between these may be identified.  In
order to carry out this analysis, we will first develop a few techical tools.

The first such tool is intended to point out a situation in which uniqueness does fail.

\state Lemma \label RangeInsv
  Let $S$ be a
  $0$-left cancellative
  semigroup admitting least common multiples.
  Let $\Lambda $ be a finite subset of $\tS $, with $\Lambda \cap S\neq \emptyset $, and let $u,v\in \Lambda $.  Given
$w$ in $S$, suppose that the range of
  $\theta _uf_\Lambda \theta _v\inv $
  is contained in $E_{uw}$.  Then
  \iItemize
  \iItem $F_\Lambda \subseteq E_w$,
  \iItem $\theta _uf_\Lambda \theta _v\inv = \theta _{uw}f_{\Lambda w} \theta _{vw}\inv $.

\Proof
  The hypothesis about the range of $\theta _uf_\Lambda \theta _v\inv $ implies that
  $$
  \theta _uf_\Lambda \theta _v\inv =
  e_{uw} \theta _uf_\Lambda \theta _v\inv =
  \theta _{uw}\theta _{uw}\inv \theta _uf_\Lambda \theta _v\inv =
  \theta _{uw}\theta _w\inv \theta _u\inv \theta _uf_\Lambda \theta _v\inv \quebra =
  \theta _{uw}\theta _w\inv f_uf_\Lambda \theta _v\inv =
  \theta _{uw}\theta _w\inv f_\Lambda \theta _v\inv \={Covar.ii}
  \theta _{uw}f_{\Lambda w} \theta _w\inv \theta _v\inv =
  \theta _{uw}f_{\Lambda w} \theta _{vw}\inv ,
  $$
  proving (ii).  Given that $u,v\in \Lambda $, we have
  $$
  f_\Lambda =
  f_uf_\Lambda f_v =
  \theta _u\inv \theta _uf_\Lambda \theta _v\inv \theta _v =
  \theta _u\inv e_{uw}\theta _uf_\Lambda \theta _v\inv \theta _v \={Covar.i}
  $$
  $$
  =
  \theta _u\inv \theta _u e_wf_\Lambda \theta _v\inv \theta _v =
  f_ue_wf_\Lambda f_v =
  e_wf_\Lambda .
  $$
  Therefore
  $$
  \id _{F_\Lambda } =
  \id _{E_w}\id _{F_\Lambda } =
  \id _{E_w\cap F_\Lambda },
  $$
  so
  $
  F_\Lambda = E_w\cap F_\Lambda ,
  $
  and then
  $
  F_\Lambda \subseteq E_w.
  $
  \endProof

We next prove our first uniqueness result, assuming two given representations of the same element of $\hull $ already
share some ingredients.  From now on we will have to rely on $0$-right cancellativity.

\state Lemma \label TeLabmdaIgual
  Let $S$ be a
  $0$-cancellative
  semigroup admitting least common multiples.
  For each $i=1, 2$, let $\Lambda _i$ be a finite subset of $\tS $ having a nonempty intersection with $S$, and let
$u_i, v_i\in \Lambda _i$ be such that
  $$
  \theta _{u_1}f_{\Lambda _1} \theta _{v_1}\inv =
  \theta _{u_2}f_{\Lambda _2} \theta _{v_2}\inv \neq 0, \and
  u_1=u_2.
  $$
  Then $f_{\Lambda _1} = f_{\Lambda _2}$.  In addition,
  \iItemize
  \iItem if either both $v_1$ and $v_2$ lie in $S$, or both $v_1$ and $v_2$ lie in $\tS \setminus S = \{1\}$, then
$v_1=v_2$.
  \iItem if $v_1\in S$ and $v_2=1$, then $v_1$ is an idempotent element of $S$, and $\theta _{v_1}\geq f_{\Lambda _1}$.
  \iItem same as in (ii) with subscripts ``$\, 1$'' and ``$\, 2$'' interchanged.

\Proof
  Based on (i) we will simply write $u$ for $u_1$ or $u_2$.

  Let $z$ be any nonzero element in the common domain of $\theta _uf_{\Lambda _1} \theta _{v_1}\inv $ and $\theta
_uf_{\Lambda _2} \theta _{v_2}\inv $, so that
  $z\in E_{v_1}\cap E_{v_2}$.

  By definition we have that $E_{v_i} = v_iS\setminus \{0\}$ (in fact when $v_i=1$, this is not quite the definition of
$E_{v_i}$, although it is still obviously true)
  so we may write
  $z=v_1x_1=v_2x_2$, with $x_1,x_2\in S$.  It is perhaps worth insisting that $x_1$ and $x_2$ indeed lie in $S$, as
opposed to $\tS $.

We then have
  $$
  ux_1= \theta _uf_{\Lambda _1} \theta _{v_1}\inv (z)=\theta _uf_{\Lambda _2} \theta _{v_2}\inv (z)=ux_2 \neq 0,
  $$
  so $x_1=x_2$, by $0$-left cancellativity (regardless of the fact that $\tS $ might not satisfy this property), and
hence
  $$
  v_1x_1=v_2x_2=v_2x_1\neq 0.
  \equationmark vxOneTwo
  $$

  In order to prove (i) we must check that $v_1=v_2$.  Under the first alternative of (i), this follows from $0$-right
cancellativity, while it is plain obvious under the second alternative. Still under the conditions of (i) we then have
that
  $$
  f_{\Lambda _1} =
  f_uf_{\Lambda _1}f_{v_1} =
  \theta _u\inv \theta _uf_{\Lambda _1}\theta _{v_1}\inv \theta _{v_1} \quebra =
  \theta _u\inv \theta _uf_{\Lambda _2}\theta _{v_2}\inv \theta _{v_2} =
  f_uf_{\Lambda _2}f_{v_2} = f_{\Lambda _2},
  $$
  completing the proof under (i).  So now let us assume that $v_1\in S$ and $v_2=1$.  We then have from \ref {vxOneTwo}
that $v_1x_1=x_1$, so also $v_1^2x_1=v_1x_1$, hence by $0$-right cancellativity we deduce that $v_1$ is idempotent.
Therefore $\theta _{v_1}=\theta _{v_1}\inv =f_{v_1}$, and then
  $$
  f_{\Lambda _1} =
  f_uf_{\Lambda _1}f_{v_1} =
  \theta _u\inv \theta _uf_{\Lambda _1}\theta _{v_1}\inv =
  \theta _u\inv \theta _uf_{\Lambda _2}\theta _{v_2}\inv =
  f_{\Lambda _2}.
  $$

  We finally have
  $$
  f_{\Lambda _1}\theta _{v_1} = f_{\Lambda _1}f_{v_1} = f_{\Lambda _1},
  $$
  so $\theta _{v_1}\geq f_{\Lambda _1}$, proving (ii), while (iii) is proved in a similar way.
  \endProof

The following result is our best shot at identifying relations between two descriptions of a single element of $\hull $
when uniqueness fails.

\state Theorem \label Equality
  Let $S$ be a $0$-cancellative semigroup admitting least common multiples.
  For $i=1,2$, let $\Lambda _i$ be a finite subset of $\tS $ intersecting $S$, and let $u_i,v_i\in \Lambda _i$ be such
that
  $$
  \theta _{u_1}f_{\Lambda _1} \theta _{v_1}\inv = \theta _{u_2}f_{\Lambda _2} \theta _{v_2}\inv \neq 0.
  $$
  Then there are $x_1,x_2\in \tS $, such that
  \iItemize
  \iItem $\theta _{u_i}f_{\Lambda _i} \theta _{v_i}\inv = \theta _{u_ix_i}f_{\Lambda _ix_i} \theta _{v_ix_i}\inv $, for
$i=1,2$,
  \iItem $u_1x_1=u_2x_2$, and $F_{\Lambda _1x_1}=F_{\Lambda _2x_2}$.
  \medskip \noindent Moreover at least one of the following three properties hold:
  \smallskip \item {\rm (a)} $v_1x_1=v_2x_2$, or
  \smallskip \item {\rm (b)} $v_1x_1$ is an idempotent element in $S$, and $\theta _{v_1x_1}\geq f_{\Lambda _1x_1}$, and
$v_2=x_2=1$, or
  \smallskip \item {\rm (c)} same as in (b) with subscripts ``$\, 1$'' and ``$\, 2$'' interchanged.

\Proof
  Notice that the range of the nonempty map mentioned in the hypothesis is contained in
  $E_{u_1} \cap E_{u_2}$.
  Letting $w$ be a least common multiple of $u_1$ and $u_2$, write
  $w=u_1x_1=u_2x_2$, with $x_1,x_2\in \tS $.

Since $u_1S\cap u_2S=wS$, we have that $E_{u_1}\cap E_{u_2} = E_w$, so for every $i=1,2$, the range of $\theta
_{u_i}f_{\Lambda _i} \theta _{v_i}\inv $ is contained in $E_w=E_{u_ix_i}$.  By \ref {RangeInsv} we then conclude that
  $$
  \theta _{u_i}f_{\Lambda _i} \theta _{v_i}\inv = \theta _{u_ix_i}f_{\Lambda _ix_i} \theta _{v_ix_i}\inv ,
  $$
  thus proving (i) and the first part of (ii).  Having already seen that $u_1x_1=u_2x_2=w$, notice that
  $$
  \theta _wf_{\Lambda _1x_1} \theta _{v_1x_1}\inv = \theta _wf_{\Lambda _2x_2} \theta _{v_2x_2}\inv \neq 0,
  $$
  so the conclusion follows from \ref {TeLabmdaIgual}.
  \endProof

Finally, we handle the case of a categorical at zero semigroup, generalizing several known results in the literature.
This theorem applies, in particular, to the inverse hull of a left cancellative category.

\state Theorem \label CatatZeroHull Let $S$ be a categorical at zero semigroup that is $0$-left cancellative, right
reductive, has right local units and least common multiples.  Then the non-zero elements of the inverse hull $\hull $
are precisely those elements of the form $\theta _s\theta _t\inv $ with $s^+=t^+$.  Moreover, if $s_1^+=t_1^+$ and
$s_2^+=t_2^+$, then $\theta _{s_1}\theta _{t_1}\inv =\theta _{s_2}\theta _{t_2}^{\inv }$ if and only if there exist
$x,y$ with $xy=s_1^+$, $yx=s_2^+$, $s_1x=s_2$, $t_1x=t_2$, $s_2y=s_1$ and $t_2y=t_1$.

\Proof Since $S$ is categorical at zero, we have that $f_s=\theta _{s^+}$ for all $s\in S$ by Proposition~\ref
{CategoricalLocalUnits}.  Thus if $s^+=t^+$ and $\Lambda \subseteq Ss^+$, then
  $$
  \theta _sf_{\Lambda }\theta _t\inv =\theta _s\theta _{s^+}\theta _t\inv =\theta _s\theta _t\inv .
  $$
  Also note that
  $$
  F_s=s^+S\setminus \{0\}=t^+S\setminus \{0\}\supseteq tS\setminus \{0\}=E_t,
  $$
  and so $\theta _s\theta _t\inv \neq 0$ if $s^+=t^+$. Observe that $\theta _s\theta _t\inv \colon tS\setminus \{0\}\to
sS\setminus \{0\}$.

Assume that $s_1,s_2,t_1,t_2,x,y$ are as above.  Then $xS=s_1^+S=t_1^S$ since $xy=s_1^+$ and $s_1x=s_2\neq 0$ implies
$s_1^+x\neq 0$ and so $x\in s_1^+S$.  Thus $\theta _x\theta _x\inv =\theta _{s_1^+}$ and hence $\theta _{s_2}\theta
_{t_2}\inv = \theta _{s_1}\theta _x\theta _x\inv \theta _{t_1}\inv =\theta _{s_1}\theta _{t_1}\inv $.

Conversely, if $\theta _{s_1}\theta _{t_1}\inv = \theta _{s_2}\theta _{t_2}\inv $, then $s_1S=s_2S$ and so $s_1x=s_2$
and $s_2y=s_1$ for some $x,y\in S$.  Then $s_1xy=s_1$ and $s_2yx=s_2$.  Therefore, $xy=s_1^+$ and $yx=s_2^+$ by $0$-left
cancellativity.
  Now $s_1=\theta _{s_1}\theta _{t_1}\inv (t_1)=\theta _{s_2}\theta _{t_2}\inv (t_1)$
   and so $t_1=t_2z$ with $s_2z=s_1=s_2y$, whence $y=z$.
  Therefore, $t_1=t_2y$. Similarly, $t_2=t_1w$ with $s_1w=s_2=s_1x$ and hence $w=x$, whence $t_2=t_1x$.  This completes
the proof.  \endProof

Notice that in Theorem \ref {CatatZeroHull} if $S=S(C)$ where $C$ is a left cancellative category with least common
multiples, then the elements $x,y$ above will be isomorphisms in $C$.

\section Finitely aligned semigroups

We consider here a generalization of the lcm property.

\definition \label FinitelyALigned
  A $0$-left cancellative semigroup $S$ is said to be \emph {finitely aligned}, or to have the \emph {(right) Howson
property}, provided that, for every $s$ and $t$ in $S$, there is a finite sequence $\{r_j\}_{j=1}^n$ of elements of $S$
such that
  \iItemize
  \iItem
  $
  sS\cap tS=\bigcup _{j=1}^n r_jS,
  $
  \iItem both $s$ and $t$ divide $r_j$, for every $j=1,\ldots ,n$.

Notice that when $S$ has right local units one has that $r\in rS$, for all $r\in S$, so condition \ref
{FinitelyALigned.ii} above follows from \ref {FinitelyALigned.i}.

\fix For simplicity, we will stick to the case that $S$ has right local units so for this section, let us assume that
$S$ is a $0$-left cancellative semigroup admitting right local units.

Under the present hypothesis we then have that the intersection of finitely generated right ideals is finitely
generated, hence the motivation for the terminology adopted above \cite {Howson}.

\definition We shall say that $S$ is \emph {strongly finitely aligned} if, for all $s,t\in S$, there exists a finite set
$B\subseteq S\setminus \{0\}$ (possibly empty) such that $sS\cap tS=BS$ and $bS\cap b'S=\{0\}$ for $b\neq b'\in B$.
Here we interpret $\emptyset S=0$.  We call $B$ a basis for $sS\cap tS$.

For example, any $0$-left cancellative right lcm semigroup with right local units is strongly finitely aligned.  We
shall get more examples from higher rank graphs.

 Our main example of a strongly finitely aligned semigroup comes from a finitely aligned higher rank graph \cite
{KumjianPask}; the reader is referred to \cite {KumjianPask} for all undefined notions.  Let $\Lambda $ be a $k$-graph
with degree functor $d$.  We write use $\vee $ for the pointwise maximum on $\Bbb N^k$.  Put $S=\Lambda \cup \{0\}$
where all undefined products in $\Lambda $ are made $0$.  Then $S$ is a $0$-cancellative, categorical at zero, right and
left reductive semigroup with local units.  We claim that $S$ is strongly finitely aligned in our sense if and only if
$\Lambda $ is finitely aligned in the usual sense; note that it is singlely aligned when $S$ has lcms.

 First assume that $S$ is strongly finitely aligned in our sense and let $\lambda ,\mu \in \Lambda $.  Let $B$ be a
basis for $\lambda S\cap \mu S$ and assume that $\lambda S\cap \mu S\neq \{0\}$.  We claim that if $\gamma \in B$, then
$d(\gamma )=d(\lambda )\vee d(\mu )$.  Indeed, by the unique factorization property, we must have $\gamma =\rho \eta $
with $d(\rho )=d(\lambda )\vee d(\mu )$ and $\rho \in \lambda S\cap \mu S$ by the unique factorization property.  Thus
$\rho =\gamma '\tau $ with $\gamma '\in B$.  From $\gamma =\gamma '\tau \eta $, we deduce that $\gamma =\gamma '$ and
$d(\tau \eta )=0$, i.e., $\tau $ and $\eta $ are identities.  Thus $d(\gamma )=d(\rho ) = d(\lambda )\vee d(\mu )$.
Next observe that if $\rho \in \lambda S\cap \mu S$ and $d(\rho )=d(\lambda )\vee d(\mu )$, then $\rho \in B$.  Indeed,
$\rho =\gamma \eta $ with $\gamma \in B$.  Hence $d(\rho )=d(\gamma )$ by what we just observed and so $\rho =\gamma $.
It now follows that there is a bijection $\Lambda ^{min}(\lambda ,\mu )\to B$ given by $(\alpha ,\beta )\mapsto \lambda
\alpha =\mu \beta $ and so $\Lambda $ is finitely aligned.

 Next assume that $\Lambda $ is finitely aligned in the usual sense and suppose that $\lambda S\cap \mu S\neq \{0\}$.
Let $B=\{\lambda \alpha \mid (\alpha ,\beta )\in \Lambda ^{min}(\lambda ,\mu )\}$.  We claim that $B$ is a basis for
$\lambda S\cap \mu S$.  Clearly, $B\subseteq \lambda S\cap \mu S$. If $0\neq \gamma \in \lambda S\cap \mu S$, then by
the unique factorization property we must have that $\gamma =\rho \eta $ with $d(\rho )=d(\lambda )\vee d(\mu )$. Then
if $\rho = \lambda \alpha =\mu \beta $, we have that $(\alpha ,\beta )\in \Lambda ^{min}(\lambda ,\mu )$ and so $\rho
=\lambda \alpha \in B$.  Also if $\gamma ,\gamma '\in B$ and $0\neq \tau \in \gamma S\cap \gamma 'S$, then from
$d(\gamma )=d(\gamma ')$ we must have $\gamma =\gamma '$ by the unique factorization property.  Thus $B$ is a basis.  We
conclude that $S$ is strongly finitely aligned.  We shall show in our sequel paper that the tight $C^*$-algebra of the
strongly finitely aligned $0$-cancellative semigroup $S(\Lambda )$ associated to a finitely aligned $\Lambda $ is the
higher rank graph $C^*$-algebra.

 For the remainder of this work, we will focus on the case of right lcm semigroups, but future work will consider
further the finitely aligned case.

\section Free product

In this section, we study free products of $0$-left cancellative monoids in order to produce new examples lcm monoids
and finitely aligned monoids.

If $M$ and $N$ are monoids, their free product $M* N$ is their coproduct in the category of monoids.  It is a standard
fact that $M$ and $N$ embed in their free product and each element of $M* N$ can be uniquely expressed as a product of
the form $m_1n_1m_2n_2\cdots m_kn_k$ with $m_i\in M\setminus \{1\}$, for $2\leq i\leq k$, $n_i\in N\setminus \{1\}$ for
$1\leq i\leq k-1$ and $m_1\in M$, $n_k\in N$.

Assume now that $M$ and $N$ are monoids with zero.  Let us denote by $M\aster 0 N$ their coproduct in the category of
monoids with zero and call it the $0$-free product of $M$ and $N$.  In other words $M\aster 0 N$ is a monoid with zero
equipped zero-preserving homomorphisms $M\to M\aster 0 N$ and $N\to M\aster 0 N$ such that any zero-preserving
homomorphisms $M\to T$ and $N\to T$ to a monoid with zero $T$ `extend' uniquely to $M\aster 0 N$.

\state Proposition \label freeprodasreesquotient Let $M$ and $N$ be monoids with zero.  Let $I$ be the ideal of $M*N$
generated by the respective zeroes $0_M, 0_N$ of $M$ and $N$.  Then $M\aster 0 N\cong (M*N)/I$.

\Proof It is clear that if $M\to T$ and $N\to T$ are zero-preserving maps, then their extension to $M*N$ maps $I$ to $0$
and hence factors uniquely through $(M*N)/I$.  It follows that $(M*N)/I$ (equipped with the canonical maps $M\to
(M*N)/I$ and $N\to (M*N)/I$) has the correct universal property to be $M\aster 0 N$.
    \endProof

\state Corollary \label normalformfree Suppose that $M$ and $N$ are non-trivial monoids with zero.  Then $M$ and $N$
embed into $M\aster 0 N$ and each non-zero element of $M\aster 0 N$ can be uniquely written in the form
$m_1n_1m_2n_2\cdots m_kn_k$ with $m_i\in M\setminus \{0,1\}$, for $2\leq i\leq k$, $n_i\in N\setminus \{0,1\}$ for
$1\leq i\leq k-1$ and $m_1\in M\setminus \{0\}$, $n_k\in N\setminus \{0\}$.

\Proof This is immediate from the normal form theorem for free products of monoids and the observation that if $M$ and
$N$ are non-trivial, then an element of $M*N$ belongs to the ideal generated by the zeroes of $M$ and $N$ if and only if
its normal form contains a $0$ syllable.
    \endProof

For a non-zero, non-identity element $u$ of $M\aster 0 N$, we say that the normal form of $u=m_1n_1m_2n_2\cdots m_kn_k$
ends in an $M$-syllable if $n_k=1$, and that $m_k$ is the last syllable of $u$, and otherwise we say that it ends in an
$N$-syllable and $n_k$ is the last syllable of $u$.  The total number of non-identity syllables in the normal form of
$u$ is called the syllable length of $u$.  We take $0$ and $1$ to have syllable length $0$.

As a consequence, one can show that the $0$-free product of non-trivial $0$-left ($0$-right) cancellative monoids is
$0$-left ($0$-right) cancellative.

\state Theorem \label freeiscanc Let $M$ and $N$ be non-trivial $0$-left cancellative monoids.  Then $M\aster 0 N$ is
$0$-left cancellative.  The dual result holds for $0$-right cancellative monoids.

\Proof By induction on syllable length, it is enough to show that if $x\in M\cup N$ and $xu=xv\neq 0$, then $u=v$.  By
symmetry, we may assume that $x\in M\setminus \{0\}$.  Let $u=m_1n_1\cdots m_kn_k$ and $v=m_1'n_1'\cdots m_r'n_r'$ be
the normal forms as per \ref {normalformfree}.  Then $xu$ has normal form $(xm_1)n_1\cdots m_kn_k$ and $xv$ has normal
form $(xm_1')n_1'\cdots m_r'n_r'$ and so $xm_1=xm_1'\neq 0$, $k=r$ and $m_i=m_i'$, $n_j=n_j'$ for $2\leq i\leq k$ and
$1\leq j\leq k$.  As $M$ is $0$-left cancellative $m_1=m_1'$ and so $u=v$.
    \endProof

Next we want to prove that being an lcm monoid or a (strongly) finitely aligned monoid is closed under $0$-free product.
We begin by describing a set of representatives of the principal right ideals of a $0$-free product.  Recall that the
$\mathcal R$-class of an element of a monoid is the set of all elements which generate its principal right ideal.

\state Proposition \label proprightideals Let $M$ and $N$ be non-trivial monoids with $0$.  Let $T_M$ and $T_N$ be a
complete set of representatives of the non-zero $\mathcal R$-classes of $M$ and $N$, respectively, with $1\in T_M$ and
$1\in T_N$.  Then a complete set of representatives of the $\mathcal R$-classes of $M\aster 0 N$ consists of $0,1$ and
all elements whose normal forms end in a syllable from $(T_M\cup T_N)\setminus \{1\}$.

\Proof We prove that the $\mathcal R$-class of $w\in M\aster 0 N$ appears in the list above by induction on the syllable
length of $w$.  If the syllable length of $w$ is zero, there is nothing to prove.  Suppose that $w\in M\aster 0
N\setminus \{0,1\}$ has normal form $ux$ with $x\in M\setminus \{0,1\}\cup N\setminus \{0,1\}$ and $u$ has syllable
length one less than $x$.  Without loss of generality, assume that $x\in M$.  If $x$ is a right invertible element of
$M$, then $ux$ generates the same right ideal as $u$ and the result follows by induction.  If $x$ is not right
invertible, then $xM=x'M$ for a unique $x'\in T_M\setminus \{1\}$.  Then $ux'$ generates the same right ideal as $w=ux$
and $ux'$ belongs to our list of representatives.

First note that if $u$ ends in an $N$-syllable and $x\in T_M\setminus \{1\}$, then $1\notin xM$ and so $u$ is a prefix
of the normal form of any element of $ux(M\aster 0 N)$.  A similar observation holds for $ux$ if $u$ ends in an
$M$-syllable and $x\in T_N\setminus \{1\}$.  It now follows that the elements on our list generate distinct principal
right ideals.
    \endProof

\state Theorem \label closelcmunderfree Let $M$ and $N$ be non-trivial $0$-left cancellative lcm ((strongly) finitely
aligned) monoids.  Then $M\aster 0 N$ is also a $0$-left cancellative lcm ((strongly) finitely aligned) monoid.

\Proof We know $M\aster 0 N$ is $0$-left cancellative by \ref {freeiscanc}.  Let $T_M$ and $T_N$ be a complete set of
representatives of the non-zero $\mathcal R$-classes of $M$ and $N$, respectively, with $1\in T_M$ and $1\in T_N$.  Then
a complete set of representatives of the $\mathcal R$-classes of $M\aster 0 N$ consists of $0,1$ and all elements whose
normal forms end in a syllable from $(T_M\cup T_N)\setminus \{1\}$ by \ref {proprightideals}.  If $ux$ is a normal form
with $x\in T_M\setminus \{1\}$ and $u$ empty or ending in an $N$-syllable, then the non-zero right multiples of $ux$
have normal form $uyz$ where $y$ is a right multiple of $x$ in $M$ and $z$ is empty or has first syllable from $N$.  The
situation is dual if $x\in T_N\setminus \{1\}$ and $u$ is empty or ends in an $M$-syllable.  It follows that if $R$ and
$R'$ are principal right ideals with $R\cap R'\neq \{0\}$, then either $R\subseteq R'$, $R'\subseteq R$ or $R$ and $R'$
are generated by $ux$ and $ux'$ (written in normal form) with $x,x'\in (T_M\cup T_N)\setminus \{1\}$.  If $M$ and $N$
are lcm monoids, then the least common multiple of $ux,ux'$ is $uy$ where $y$ is a least common multiple of $x,x'$ in
the respective factor $M$ or $N$.  If $M$ and $N$ are (strongly) finitely aligned and $B$ is a finite generating set
(basis) for $xM\cap x'M$, if $x,x'\in T_M$, or for $xN\cap x'N$, if $x,x'\in T_N$, then $uB$ is a generating set (basis)
for $ux(M\aster 0 N)\cap ux'(M\aster 0 N)$.
    \endProof

Note that if $G$ is a group, then $G\cup \{0\}$ is always a $0$-left cancellative lcm monoid.  So we can build lots of
strongly finitely aligned $0$-cancellative monoids by taking free products of groups with adjoined zeroes and higher
rank $k$-graph on one vertex.

\section Strings

Regarding the semigroup $\SX $ constructed from a subshift $\X $, as in example \ref {DefineSubshSgrp}, suppose we want
to recover $\X $ from the algebraic structure of $\SX $.  Given a generic element of $\X $, namely an \emph {infinite}
word
  $$
  x = x_1 x_2 x_3 \ldots ,
  $$
  we may \emph {approximate} $x$ by members of $\SX $ by considering the sequence of \emph {finite} words
  $\{s_n\}_{n\in {\bf N}}$, given by
  $$
  s_n = x_1 x_2 \ldots x_n, \for n\in {\bf N}.
  $$
  We will not attempt to give a precise definition for the meaning of the word \emph {approximate} in this context, but
we will instead introduce a general concept which is expected to play the role of the above heuristic method for an
arbitrary semigroup.

\fix Throughout this section $S$ will be a fixed $0$-left cancellative semigroup.

\definition \label DefineStrin
  A nonempty subset
  $\sigma \subseteq S$
  \# $\sigma $; A string;
  is said to be a \emph {string} in $S$, if
  \iItemize
  \iItem $0\notin \sigma $,
  \iItem for every $s$ and $t$ in $S$, if $s\|t$, and $t\in \sigma $, then $s\in \sigma $, \iItemmark StringHereditary
  \iItem for every $s_1$ and $s_2$ in $\sigma $, there is some $s$ in $\sigma $ such that $s_1\| s$, and $s_2\| s$.
\iItemmark StringDirected

\bigskip An elementary example of a string is the set of divisors of any nonzero element $s$ in $S$, namely,
  $$
  \delta _s= \{t\in S: t\| s\}.
  \equationmark StringDivisors
  $$
  \# $\delta _s$; The string of divisors of $s$;

Considering the partially ordered set $\mathcal P$ formed by all principal right ideals $s\tS $ of $S$, whose smallest
element is $\{0\}$, notice that for every string $\sigma $, one has that
  $$
  \mathcal F:= \{s\tS : s\in \sigma \}
  $$
  is a proper filter \cite [12.1]{actions}
  on $\P $.  Conversely, if $\mathcal F$ is a proper filter on $\P $, then
  $$
  \sigma := \big \{s\in S: s\tS \in \mathcal F\big \}
  $$
  is a string.  Therefore, strings are essentially the same as proper filters on $\P $.

Strings often contain many elements, but there are some exceptional strings consisting of a single semigroup element. To
better study these it is useful to introduce some terminology.

\definition \label DefinePrime
  Given a nonzero $s$ in $S$ we will say that $s$ is:
  \iItemize
  \iItem \emph {prime}, if the only divisor of $s$ is $s$, itself, or, equivalently, if $\delta _s=\{s\}$,
  \iItem \emph {irreducible}, if there are no two elements $x$ and $y$ in $S$ such that $s=xy$, or, equivalently, if
$s\notin S^2$.

It is evident that any irreducible element is prime, but there might be prime elements which are not irreducible.  For
example, in the semigroup $S=\{0, s, e\}$, with multiplication table given by

  \bigskip \hfill \vbox {\begingroup
  \def \tabrule {\noalign {\hrule }}
  \def \vr {\vrule height 12pt}
  \def \bx #1{\hbox to 20pt {\ \hfill #1\hfill }}
  \offinterlineskip
  \halign {
    \strut
    \vr #&\vr #&\vr #&\vr #\vr \cr
    \tabrule
   \bx {$\times $} & \bx {$0$} & \bx {$e$} & \bx {$s$}\cr \tabrule
   \bx {$0$} & \bx {$0$} & \bx {$0$} & \bx {$0$}\cr \tabrule
   \bx {$e$} & \bx {$0$} & \bx {$e$} & \bx {$0$}\cr \tabrule
   \bx {$s$} & \bx {$0$} & \bx {$s$} & \bx {$0$}\cr \tabrule
  }
  \endgroup }\hfill \null

  \bigskip \noindent one has that $s$ is prime but not irreducible because $s=se\in S^2$.

\state Proposition \label PrimeStrings A singleton $\{s\}$ is a string if and only if $s$ is prime.

\Proof If $s$ is prime then the singleton $\{s\}$ coincides with $\delta _s$, and hence it is a string.  Conversely,
supposing that $\{s\}$ is a string, we have by \ref {StringHereditary} that $\delta _s\subseteq \{s\}$, from where it
follows that $s$ is prime.
  \endProof

\definition \label DefineSStar
  The set of all strings in $S$ will be denoted by $S^\star $.
  \# $S^\star $; The set of all strings in $S$;

From now on our goal will be to define an action of $S$ on $S^\star $.

\state Proposition \label MappingStrings
  Let $\sigma $ be a string in $S$, and let $r\in S$.  Then
  \iItemize
  \iItem if\/ $0$ is not in $r\sigma $, one has that
  \setbox 1\hbox {$\{t\in S: t\| rs, \hbox { for some } s\in \sigma \}$}
  \setbox 0\hbox {$r\inv *\sigma $}
  $$
  \hbox to \wd 0{\hfill $r*\sigma $}:= \copy 1
  $$
  \# $r*\sigma $; Product of the semigroup element $r$ by the string $\sigma $;
  is a string whose intersection with $rS$ is nonempty.
  \iItem If $\sigma $ is a string whose intersection with $rS$ is nonempty, then
  $$
  \copy 0:= \hbox to \wd 1 {$\{t\in S: rt\in \sigma \}$\hfill }
  $$
  \# $r\inv *\sigma $; Inverse product of the semigroup element $r$ by the string $\sigma $;
  is a string, and $0$ is not in $r(r\inv *\sigma )$.

\Proof
  In order to prove that $r*\sigma $ satisfies \ref {DefineStrin.i} we
  argue by contradiction: if $0$ is in $r*\sigma $, then there exists some $s$ in $\sigma $ such that $0\| rs$, whence
$rs=0$, which is ruled out by hypotheses.

Since division is a transitive relation, it is clear that $r*\sigma $ satisfies \ref {StringHereditary}.

Regarding \ref {StringDirected}, for each $i=1, 2$, let $t_i\in r*\sigma $, and pick $s_i\in \sigma $, such that $t_i\|
rs_i$.  We may then choose $u_i\in \tS $ so that
  $$
  t_iu_i=rs_i.
  $$

  Since $s_1$ and $s_2$ lie in the $\sigma $, we may furthermore choose $v_1,v_2\in \tS $, such that
  $$
  s_1v_1=s_2v_2\in \sigma .
  $$
  Setting $w_i=u_iv_i$, we then have that
  $$
  t_1w_1 =
  t_1u_1v_1 =
  rs_1v_1 =
  rs_2v_2 =
  t_2u_2v_2 =
  t_2w_2.
  $$ Given that $s_1v_1\in \sigma $, it is clear that $rs_1v_1\in r*\sigma $, so we may complement our findings above by
writing
  $$
  t_1w_1 = t_2w_2\in r*\sigma ,
  $$
  hence proving \ref {StringDirected} for $r*\sigma $.
  The final requirement of (i) is easily checked by noticing that
  $$
  \emptyset \neq r\sigma \subseteq (r*\sigma ) \cap rS.
  $$ This also proves the required condition that $r*\sigma $ be nonempty.

Addressing (ii), notice that since $\sigma \cap rS \neq \emptyset $, we may pick some $t$ in $S$ such that $rt\in \sigma
$, so we see that $t$ lies in $r\inv *\sigma $, proving the latter to be a nonempty set.  Also, $0$ is not in $r\inv
*\sigma $, since otherwise $r\cdot 0$ would be in $\sigma $, hence $r\inv *\sigma $ satisfies \ref {DefineStrin.i}.
Assuming that
  $$
  s\| t\in r\inv *\sigma ,
  $$
  we have that
  $$
  rs\| rt\in \sigma ,
  $$
  and we deduce that $rs\in \sigma $, hence $s \in r\inv *\sigma $, thus proving that $r\inv *\sigma $ satisfies \ref
{StringHereditary}.

 Let us now prove that $r\inv *\sigma $ satisfies \ref {StringDirected}.  For this let us pick $t_1$ and $t_2$ in $r\inv
*\sigma $, so that $rt_1,rt_2\in \sigma $, and then we may find $u_1,u_2\in \tS $, such that
  $$
  rt_1u_1 = rt_2u_2 \in \sigma .
  $$

  By $0$-left cancellativity we deduce that $t_1u_1 = t_2u_2$, and it is clear that the element represented by either
side of this equality lies in $r\inv *\sigma $.

To finish we observe that if $t\in r\inv *\sigma $, then $rt$ lies in $\sigma $, so $rt\neq 0$.
  \endProof

It should be noted that, under the assumptions of \ref {MappingStrings.i}, one has that
  $$
  r\sigma \subseteq r*\sigma ,
  \equationmark ProdInForw
  $$
  and in fact $r*\sigma $ is the hereditary closure of $r\sigma $ relative to the order relation \ref
{OrderDivision}. In addition we have:

\state Proposition \label StringContained
  If $\sigma $ and $\tau $ are strings with $r\sigma \subseteq \tau $, then $\sigma $ satisfies the assumption of \ref
{MappingStrings.i}, and $r*\sigma \subseteq \tau $.

\Proof
  If $r\sigma \subseteq \tau $, then $0$ is clearly not in $r\sigma $, while the inclusion $r*\sigma \subseteq \tau $
follows from the fact that $\tau $ is hereditary, and the observation already made that $r*\sigma $ is the hereditary
closure of $r\sigma $.  \endProof

We will now define a representation of $S$ on the set $S^\star $ of all strings in $S$, as follows.

\state Proposition \label IntroStarAction
  For each $r$ in $S$, put
  \# $\theta ^\star $; Representation of $S$ on $S^\star $;
  \# $F^\star _r$; Domain of $\theta ^\star _r$;
  \# $E^\star _r$; Range of $\theta ^\star _r$;
  $$
  F^\star _r=\{\sigma \in S^\star : r\sigma \not \ni 0\}, \and
  E^\star _r=\{\sigma \in S^\star : \sigma \cap rS \neq \emptyset \}.
  $$
  Also let
  $$
  \theta ^\star _r\colon F^\star _r\to E^\star _r
  $$
  be defined by $\theta ^\star _r(\sigma ) = r*\sigma $, for every $\sigma \in F^\star _r$.  Then:
  \iItemize
  \iItem $\theta ^\star _r$ is bijective, and its inverse is the mapping defined by
  $$
  \sigma \in E^\star _r\mapsto r\inv *\sigma \in F^\star _r.
  $$
  \iItem Viewing $\theta ^\star $ as a map from $S$ to $\I (S^\star )$, one has that $\theta ^\star $ is a
representation of $S$ on $S^\star $.

\Proof
  (i)\enspace For $\sigma \in F^\star _r$, and $t\in S$, one has that
  $$
  t \in r\inv *(r*\sigma ) \iff
  rt \in r*\sigma \iff
  \exists s\in \sigma ,\ rt \| rs \quebra \iff
  \exists s\in \sigma ,\ \exists u\in \tS ,\ rtu= rs
  $$
  Observing that $rs\neq 0$, by hypothesis, we may use $0$-left cancellativity to conclude that the above is equivalent
to
  $$
  \exists s\in \sigma ,\ \exists u\in \tS ,\ tu= s \iff
  \exists u\in \tS ,\ tu\in \sigma \iff t\in \sigma .
  $$
  This shows that $r\inv *(r*\sigma ) =\sigma $, and we will next prove that $r*(r\inv *\sigma ) =\sigma $.  For this,
pick $\sigma $ in $E^\star _r$, and notice that for any given $t\in S$, one has that
  $$
  t\in r*(r\inv *\sigma ) \iff
  \exists s\in r\inv *\sigma ,\ t \| rs \quebra \iff
  \exists s\in S,\ (rs\in \sigma ) \wedge (t \| rs).
  \equationmark CondForVaiVolta
  $$
  The last sentence above implies that $t$ divides an element of $\sigma $, and hence that $t\in \sigma $, therefore
showing that
  $r*(r\inv *\sigma ) \subseteq \sigma $.

To prove the reverse inclusion, pick $t$ in $\sigma $.  Observing that $\sigma $ is in $E^\star _r$, we may find $x$ in
$S$ such that $rx\in \sigma $.  Using \ref {StringDirected}, pick $u$ and $v$ in $\tS $, such that $tu=rxv\in \sigma $,
so $t\|rxv$ and, upon choosing $s=xv$, we see that \ref {CondForVaiVolta} holds.  So $t\in r*(r\inv *\sigma )$,
concluding the proof of (i).

\itmproof (ii)
  It is clear that $F^\star _0 = E^\star _0 = \emptyset $, so $\theta ^\star _0$ is the empty map. Given $r_1$ and $r_2$
in $S$, we must now prove that the domain of $\theta ^\star _{r_1}\circ \theta ^\star _{r_2}$ coincides with the domain
of $\theta ^\star _{r_1r_2}$, namely $F^\star _{r_1r_2}$, and that
  $$
  r_1*(r_2*\sigma ) = (r_1r_2)*\sigma ,
  \equationmark CondForHomo
  $$
  for every $\sigma $ in the above common domain.

For this, notice that a given string $\sigma $ lies in the domain of $\theta ^\star _{r_1}\circ \theta ^\star _{r_2}$ if
and only if
  $$
  (\sigma \in F^\star _{r_2}) \wedge (r_2*\sigma \in F^\star _{r_1}) \iff
  (0\notin r_2\sigma ) \wedge \big (0\notin r_1(r_2*\sigma )\big ).
  \equationmark CondForDomCompos
  $$

Suppose by way of contradiction that a string $\sigma $ satisfying the above equivalent condition fails to be in
$F^\star _{r_1r_2}$.  Then
  $$
  0 \in r_1r_2\sigma \explain {ProdInForw}\subseteq r_1(r_2*\sigma ),
  $$
  which contradicts \ref {CondForDomCompos}.  Therefore we see that the domain of $\theta ^\star _{r_1}\circ \theta
^\star _{r_2}$ is contained in $F^\star _{r_1r_2}$.

  In order to prove the reverse inclusion, pick $\sigma $ in $F^\star _{r_1r_2}$.  Then $0\notin r_1r_2\sigma $, from
where one deduces that $0\notin r_2\sigma $, thus verifying the first condition in the right hand side of \ref
{CondForDomCompos}, and we claim that the second condition also holds.  Arguing by contradiction, suppose that $0\in
r_1(r_2*\sigma )$, which is to say that $r_1t=0$, for some $t$ in $r_2*\sigma $.  It follows that there are elements $u$
in $\tS $, and $s$ in $\sigma $, such that $tu=r_2s$, and then
  $$
  0 = r_1tu = r_1r_2s \in r_1r_2\sigma ,
  $$
  a contradiction.  This proves that the whole right hand side of \ref {CondForDomCompos} holds, and hence that $\sigma
$ lies in the domain of $\theta ^\star _{r_1}\circ \theta ^\star _{r_2}$, as needed.

In order to prove \ref {CondForHomo}, let $\sigma $ be in $F^\star _{r_1r_2}$, and notice that
  $$
  r_1r_2\sigma \explain {ProdInForw}\subseteq
  r_1(r_2*\sigma ) \explain {ProdInForw}\subseteq
  r_1*(r_2*\sigma )
  $$
  so the hereditary closure of $r_1r_2\sigma $, namely $(r_1r_2)*\sigma $, is contained in $r_1*(r_2*\sigma )$.  On the
other hand, if $t\in r_1*(r_2*\sigma )$, then $tu=r_1y$, for suitable $u\in \tS $, and $y\in r_2*\sigma $.  We may
moreover write $yv=r_2s$, with $v\in \tS $, and $s\in \sigma $, so
  $$
  tuv = r_1yv = r_1r_2s \in r_1r_2\sigma \explain {ProdInForw}\subseteq (r_1r_2)*\sigma ,
  $$
  whence
  $t\in (r_1r_2)*\sigma $.  This shows that $r_1*(r_2*\sigma )\subseteq (r_1r_2)*\sigma $, and hence verifies \ref
{CondForHomo}, concluding the proof of (ii).
  \endProof

Useful alternative characterizations of $F^\star _r$ and $E^\star _r$ are as follows:

\state Proposition \label FRFstarR
  Given $r$ in $S$, and given any string $\sigma $ in $S^\star $, one has that:
  \iItemize \def \iff {\ \mathrel {\Leftrightarrow }\ }
  \iItem $\sigma \in F^\star _r \Iff \sigma \subseteq F^\theta _r$,
  \iItem $\sigma \in E^\star _r \Iff \sigma \cap E^\theta _r\neq \emptyset ,$
  \iItem $\sigma \in E^\star _r \Imply r\in \sigma $.  In addition the converse holds provided $r\in rS$ (e.g.~if $S$
has right local units).

\Proof (i)\enspace A given string $\sigma $ lies in $F^\star _r$, if and only if $rs\neq 0$, for every $s$ in $\sigma $,
which is to say that $\sigma \subseteq F^\theta _r$.

\itmproof (ii)
  A string $\sigma $ belongs to $E^\star _r$, if and only if $\sigma \cap rS\neq \emptyset $, but since $0$ is not in
$\sigma $, this is equivalent to
  $$
  \emptyset \neq \sigma \cap (rS\setminus \{0\}) = \sigma \cap E^\theta _r.
  $$

\itmproof (iii)
  If $\sigma \in E^\star _r$, there exists some $s$ in $S$ such that $rs\in \sigma $, hence $r\in \sigma $, by \ref
{StringHereditary}.
  Conversely, if $r\in rS$, and $r\in \sigma $, then
  $r\in rS\cap \sigma $, whence $rS\cap \sigma $ is nonempty and we see that $\sigma \in E^\star _r$.
  \endProof

\state Remark \label RmkrSasrightideal If $rS$ is generated by $X\subseteq S$ as a right ideal, that is, $rS=\bigcup
_{s\in X} s\tS $, then $\sigma \in E^\star _r$ if and only if $\sigma \cap X\neq \emptyset $.

Recall from \ref {FormOfHull} that, when $S$ has least common multiples, every $\theta $-constructible subset of $S'$
has the form $\theta _u(F^\theta _\Lambda )$, where $\Lambda \subseteq \tS $ is finite, $\Lambda \cap S\neq \emptyset $,
and $u\in \Lambda $.
  By analogy this suggests that it might also be useful to have a characterization of $\theta ^\star _u(F^\star _\Lambda
)$ along the lines of \ref {FRFstarR}.

\state Proposition \label FRFstarRTwo
  Let $\Lambda $ be a finite subset of $\tS $ having a nonempty intersection with $S$, and let $u\in \Lambda $.  Then
  $\theta ^\star _u(F^\star _\Lambda )$
  consists precisely of the strings $\sigma $ such that
  $$
  \emptyset \neq
  \sigma \cap E^\theta _u \subseteq
  \theta _u(F^\theta _\Lambda ).
  $$

\Proof
  Noticing that $\theta ^\star _u(F^\star _\Lambda )$ is a subset of $E^\star _u$, one has that any given string $\sigma
$ lies in $\theta ^\star _u(F^\star _\Lambda )$ if and only if
  $$
  \big (\sigma \in E^\star _u\big ) \wedge \big ({\theta _u^\star }\inv (\sigma ) \in F^\star _\Lambda \big ) \explain
{FRFstarR.i-ii}{\iff }
  $$$$
  \big (\sigma \cap E^\theta _u\neq \emptyset \big ) \wedge \big (u\inv * \sigma \subseteq F^\theta _\Lambda \big ),
  \equationmark EquivalencesTrFLambda
  $$
  where we observe that our application of \ref {FRFstarR.i-ii} for $u$ in $\tS $, above, is legitimate even though $u$
might not be in $S$, which is required by \ref {FRFstarR}: in the only exceptional case, namely when $u=1$, one has that
\ref {FRFstarR.i-ii} is trivially true.

Note that, since $\sigma $ does not contain $0$, the definition of $u\inv * \sigma $ is equivalent to
  $$
  u\inv * \sigma = \{t\in F^\theta _u : \theta _u(t) \in \sigma \cap E^\theta _u\} = \theta _u\inv (\sigma \cap E^\theta
_u).
  $$
  So \ref {EquivalencesTrFLambda} is further equivalent to
  $$
  (\sigma \cap E^\theta _u\neq \emptyset ) \wedge \big (\theta _u\inv (\sigma \cap E^\theta _u) \subseteq F^\theta
_\Lambda \big ) \iff
  $$$$
  \emptyset \neq \sigma \cap E^\theta _u \subseteq \theta _u(F^\theta _\Lambda ).
  \closeProof
  $$ \endProof

After \ref {IntroStarAction} we now have two natural representations of $S$, namely the regular representation $\theta $
acting on $S'$, and $\theta ^\star $ acting on $S^\star $.  We shall next prove that the correspondence $r\mapsto \delta
_r$ is covariant (i.e., $S$-equivariant) for these representation.  We begin by proving a technical result designed,
among other things, to show that the domains and ranges are matched accordingly.

\state Lemma \label DeltaMatchCosntru
  Let $\Lambda $ be a finite subset of $\tS $, with $\Lambda \cap S\neq \emptyset $, and let $u\in \Lambda $.  Then
  for every $r$ in $S$, one has that
  $$
  r\in \theta _u(F^\theta _\Lambda ) \iff \delta _r\in \theta ^\star _u(F^\star _\Lambda ).
  $$

\Proof
  If $r\in \theta _u(F^\theta _\Lambda )$, then $\delta _r\cap E^\theta _u$ is nonempty because $r$ belongs to this set.
In order to show that
 $\delta _r\in \theta ^\star _u(F^\star _\Lambda )$ it therefore suffices to show that
  $$
  \delta _r\cap E^\theta _u \subseteq \theta _u(F^\theta _\Lambda ),
  $$
  by \ref {FRFstarRTwo}.
  Given $x$ in $\delta _r\cap E^\theta _u$, we may write $x=us$, for some $s$ in $S$ and, observing that $x\|r$, we have
that $xv=r$, for some $v$ in $\tS $.  Then
  $
  r=xv=usv,
  $
  whence $sv=\theta _u\inv (r)\in F^\theta _\Lambda $, and it easily follows that $s\in F^\theta _\Lambda $, as well, as
$F^{\theta }_\Lambda $ is hereditary.  Thus
  $x=us\in \theta _u(F^\theta _\Lambda )$.  Having checked the inclusion displayed above, we have proven that $\delta
_r\in \theta ^\star _u(F^\star _\Lambda )$.  Conversely, assuming that $\delta _r\in \theta ^\star _u(F^\star _\Lambda
)$, we have by \ref {FRFstarRTwo} that
  $$
  \emptyset \neq \delta _r\cap E^\theta _u \subseteq \theta _u(F^\theta _\Lambda ).
  $$
  We may thus pick $x$ in $\delta _r\cap E^\theta _u$, so $x=us$, and $xv=r$, where $s\in S$, and $v\in \tS $, as above.
It follows that
  $$
  r=xv=usv\in \delta _r\cap E^\theta _u \subseteq \theta _u(F^\theta _\Lambda ),
  $$
  concluding the proof.  \endProof

\state Proposition The map
  $$
  \delta : s\in S'\mapsto \delta _s\in S^\star ,
  $$
  where $\delta _s$ is defined in \ref {StringDivisors}, is covariant relative to $\theta $ and $\theta ^\star $.

\Proof
  For every $r$ in $S$, we need to prove that, $\delta (F^\theta _r)\subseteq F^\star _r$, \ $\delta (E^\theta
_r)\subseteq E^\star _r$, and that the diagram
  $$
  \matrix { \offinterlineskip
    F^\theta _r & \buildrel \ds \theta _r \over \longrightarrow & E^\theta _r \cr
    \pilar {15pt}\delta \downarrow \ && \ \downarrow \delta \cr
    F^\star _r & \buildrel \ds \theta ^\star _r \over \longrightarrow & E^\star _r
    }
  $$
  commutes.

The first two facts follow immediately from \ref {DeltaMatchCosntru}, so we need only check the commutativity of the
diagram, which boils down to proving that
  $$
  r*\delta _s=\delta _{rs}, \for s \in F^\theta _r.
  $$ Given $s$ in $F^\theta _r$, for obvious reasons we have that
  $r\delta _s\subseteq \delta _{rs}$, so by \ref {StringContained} we obtain
  $
  r*\delta _s\subseteq \delta _{rs}.
  $
  On the other hand, if $t\in \delta _{rs}$, then $t\|rs$, so there exists $u$ in $\tS $, such that
  $$
  tu=rs\in r\delta _s \explain {ProdInForw}\subseteq r*\delta _s.
  $$
  Therefore
  $t$ is in $r* \delta _s$, showing the reverse inclusion $\delta _{rs}\subseteq r * \delta _s$, and hence proving our
diagram to be commutative.
  \endProof

Let us now consider the question of whether $\theta ^\star $ is an essential representation.  Recall from \ref
{EssSubset} that the essential subset for $\theta ^\star $ is
  $$
  S^\star _\ess =
  \big (\bigcup _{s\in S}F^\star _s \big ) \cup \big (\bigcup _{s\in S} E^\star _s\big ).
  $$

\state Proposition \label ZeroChars
  Let $\sigma $ be a string in $S^\star $. Then:
  \iItemize
  \iItem $\sigma $ does not belong to\/ $\bigcup _{t\in S}E^\star _t$, if and only if $\sigma =\{s\}$, for some
irreducible $s\in S$,
  \iItem $\sigma $ does not belong to $S^\star _\ess $, if and only if $\sigma =\{s\}$, for some degenerate element $s$
in $S$.

\Proof (i)\enspace Observe that, for every $t$ in $S$, one has that
  $$
  \sigma \notin E^\star _t \explain {FRFstarR.ii}\iff \sigma \cap E^\theta _t=\emptyset .
  \equationmark SigmaNotE
  $$
  Assuming that $\sigma \notin E^\star _t$, for every $t$ in $S$, we will now show that $\sigma $ is a singleton, so we
suppose that $s, t\in \sigma $.  Then, by \ref {StringDirected} there are $u,v\in \tS $, such that
  $
  su=tv\in \sigma .
  $
  If $v\neq 1$, then $v\in S$, whence
  $$
  tv\in tS\setminus \{0\} = E^\theta _t,
  $$
  contradicting \ref {SigmaNotE}.  Thus $v=1$, and the same reason shows that $u=1$, whence $s=t$, and we see that
$\sigma $ contains precisely one element.

Letting $s$ denote the single element in $\sigma $, one has that $s$ is necessarily irreducible since otherwise we could
write $s=tr$, with $t, r\in S$, and then $s\in \sigma \cap E^\theta _t$, again contradicting \ref {SigmaNotE}.

Conversely, assuming that $s$ is irreducible then $\sigma :=\{s\}$ is a string by \ref {PrimeStrings}.
  Observing that for every $t$ in $S$, one has $E^\theta _t\subseteq S^2$, it is clear that $\sigma \cap E^\theta
_t=\emptyset $, whence $\sigma \notin E^\star _t$, thanks to \ref {SigmaNotE}.

\itmproof (ii) In order to prove the ``only if" part, let $\sigma \in S^\star \setminus S^\star _\ess $.  Then clearly
$\sigma $ does not belong to\/ $\bigcup _{t\in S}E^\star _t$, so $\sigma =\{s\}$, for some irreducible $s\in S$, by (i),
and all we need to do is show that $Ss=0$.  Arguing by contradiction, suppose that $rs\neq 0$, for some $r$ in $S$.
Then $\sigma \subseteq F^\theta _r$, hence $\sigma \in F^\star _r$, by \ref {FRFstarR.i} contradicting the hypothesis.

Conversely, suppose that $\sigma =\{s\}$, where $s$ is degenerate, hence in particular irreducible.  By (i) we then have
that $\sigma $ is not in any $E^\star _t$, so it suffices to prove that $\sigma $ is not in any $F^\star _t$, either.
But since $Ss=\{0\}$, we see that $s$ is not in any $F^\theta _t$, whence $\sigma \not \subseteq F^\theta _t$, so
$\sigma \not \in F^\star _t$, by \ref {FRFstarR.i}.
  \endProof

We then obtain a result similar to \ref {RegRepEssential}:

\state Proposition \label EssThetaStar
  Denoting the essential subset for $\theta ^\star $ by $S^\star _\ess $, one has that
  $$
  S^\star \setminus S^\star _\ess = \big \{ \{s\}: s \hbox { is degenerate}\big \}.
  $$
  Therefore $\theta ^\star $ is essential if and only if $S$ possesses no degenerate elements.

Observe that when $S$ has right local units, then it has no irreducible elements, much less degenerate ones.  We
therefore obtain the following consequence of the above result and of \ref {RegRepEssential}:

\state Corollary
  If $S$ has right local units, then both $\theta $ and $\theta ^\star $ are essential representations.

Observe that the union of an increasing family of strings is a string, so any string is contained in a maximal one by
Zorn's Lemma.

\definition
  The subset of $S^\star $ formed by all maximal strings will be denoted by $S^\infty $.
  \# $S^\infty $; The set of all maximal strings;

  Our next result says that $S^\infty $ is invariant under $\theta ^\star $.

\state Proposition \label ForwardInvariance
  For every $r$ in $S$, and for every maximal string $\sigma $ in
  $F^\star _r$,
  one has that $\theta ^\star _r(\sigma )$ is maximal.

\Proof
  Suppose that
  $r*\sigma $ is contained in another string $\mu $.  Since $r*\sigma $ is in $E^\star _r$, we have that
  $$
  \emptyset \neq (r*\sigma ) \cap rS \subseteq \mu \cap rS,
  $$
  so $\mu $ is in $E^\star _r$.  It is then easy to see that
  $$
  \sigma = r\inv *(r*\sigma )\subseteq r\inv *\mu ,
  $$
  so $\sigma =r\inv *\mu $, by maximality, whence,
  $$
  r*\sigma = r*(r\inv *\mu )=\mu ,
  $$
  proving that $r*\sigma = \theta ^\star _r(\sigma )$ is maximal.
  \endProof

Observe that the above result says that $S^\infty $ is invariant under each $\theta ^\star _r $, but not necessarily
under ${\theta ^\star _r}\inv $.

An example
  to show that $S^\infty $ may indeed not be invariant under ${\theta ^\star _r}\inv $ is as follows.  Consider the
language $L$ on the alphabet $\Sigma =\{a, b\}$ given by
  $$
  L=\{a, b, aa, ba\}.
  $$
  Then $\sigma =\{b, ba\}$ is a maximal string, while ${\theta ^\star _b}\inv (\sigma ) =\{a\}$ is not maximal.

However, when $S$ is categorical at zero, the situation is much better, as we shall now see.

\state Proposition \label BackInvariance
  Suppose that $S$ is categorical at zero.
  Then, for every $r$ in $S$, and for every maximal string $\sigma $ in
  $E^\star _r$,
  one has that ${\theta ^\star _r}\inv (\sigma )$ is maximal.

\Proof
  Suppose that $\mu $ is a string with
  $$
  {\theta ^\star _r}\inv (\sigma ) \subseteq \mu .
  \equationmark Inclusions
  $$

  Notice that ${\theta ^\star _r}\inv (\sigma )$ lies in $F^\star _r$ for obvious reasons, but it is far from obvious
that the same applies to $\mu $.  However, under the present hypothesis we shall show that indeed
  $
  \mu \in F^\star _r.
  $

Given $t$ in $\mu $, choose any $s$ in ${\theta ^\star _r}\inv (\sigma )$.  Since both $t$ and $s$ lie in $\mu $, we may
pick $u, v\in \tS $, such that $su=tv\in \mu $.  Observing that ${\theta ^\star _r}\inv (\sigma )\in F^\star _r$, we see
that
  $$
  s\in {\theta ^\star _r}\inv (\sigma )\explain {FRFstarR}\subseteq F^\theta _r,
  $$
  so $rs\neq 0$.  Clearly also $su\neq 0$, so we may use the fact that $S$ is categorical at zero (even though $u$ is
not necessarily an element of $S$) to conclude that
  $$
  0\neq rsu = rtv,
  $$
  whence $rt\neq 0$, and we see that $t\in F^\theta _r$, thus proving that $\mu \subseteq F^\theta _r$.  Again by \ref
{FRFstarR} we then get that $\mu \in F^\star _r$, as desired.

  Therefore from \ref {Inclusions} it follows that
  $$
  \sigma =\theta ^\star _r\big ({\theta ^\star _r}\inv (\sigma )\big ) \subseteq
  \theta ^\star _r(\mu ),
  $$
  so $\sigma = \theta ^\star _r(\mu )$, by maximality, and
  $$
  {\theta ^\star _r}\inv (\sigma ) =
  {\theta ^\star _r}\inv \big (\theta ^\star _r(\mu )\big ) = \mu ,
  $$
  thus proving that ${\theta ^\star _r}\inv (\sigma )$ is maximal.
  \endProof

As an immediate consequence we have the following:

\state Corollary \label BackInvarianceCorol
  Let $S$ be a $0$-left cancellative, categorical at zero semigroup, and let $r\in S$.  Then:
  \iItemize
  \iItem $\theta ^\star _r(F^\star _r\cap S^\infty )\subseteq S^\infty $,
  \iItem ${\theta ^\star _r}\inv (E^\star _r\cap S^\infty )\subseteq S^\infty $, and
  \iItem $S^\infty $ is an invariant subset of $S^\star $ under the natural action of \/ $\I (S^\star ,\theta ^\star )$.

\Proof Points (i) and (ii) are restatements of \ref {ForwardInvariance} and \ref {BackInvariance}, respectively, while
(iii) is a consequence of (i--ii), as well as the fact that $\I (S^\star ,\theta ^\star )$ is generated by all of the
$\theta ^\star _r$, together with their inverses.
  \endProof

\section Open strings

\fix Throughout this section $S$ will be a fixed $0$-left cancellative semigroup, as before.

\medskip Given a string $\sigma $ in $S^\star $, and given $s,p\in S$, one has by the very definition of strings that
  $$
  sp\in \sigma \Imply s\in \sigma .
  $$
  Clearly the converse of the above implication is not true and in fact, knowing that $s$ lies in $\sigma $, does not
even guarantee the existence of some $p$ in $S$ such that $sp\in \sigma $.  To clarify this issue we introduce the
following concept:

\definition \label DefOpenString
  Given any string $\sigma $ in $S^\star $, the \emph {interior} of $\sigma $ is the subset of $S$ given by
  $$
  \interior \sigma = \{s\in S: \exists p\in S,\ sp\in \sigma \}.
  $$
  In addition, we shall say that $\sigma $ is an \emph {open} string when $\sigma =\interior \sigma $.

 Observe that when $S$ has right local units, every string is automatically open.

\state Proposition \label PropsOpenString
  Given any string $\sigma $, one has that
  \iItemize
  \iItem $\interior \sigma \subseteq \sigma $,
  \# $\interior \sigma $; The interior of the string $\sigma $;
  \iItem
  $\sigma $ fails to be open
  if and only if $\sigma =\delta _r$ for some $r$ in $S$ such that $r\notin rS$; in this case such an $r$ is unique,
  \iItem if $s_1, s_2\in \interior \sigma $, and if $r$ is a least common multiple for $s_1$ and $s_2$, then $r\in
\interior \sigma $,
  \iItem if $\interior \sigma $ is nonempty, and if $S$ admits least common multiples, then $\interior \sigma $ is a
string.

\Proof (i) \enspace Obvious in view of \ref {StringHereditary}.

\itmproof (ii)
  Assuming that $\interior \sigma \neq \sigma $ we may
  choose $r\in \sigma \setminus \interior \sigma $, so that $r\in \sigma $, but $rp\notin \sigma $, for all $p$ in $S$.
Observing that $\delta _r\subseteq \sigma $, we will show that in fact $\delta _r=\sigma $.  For this, let $s\in \sigma
$, and use \ref {StringDirected} to find $u,v\in \tS $ such that
  $$
  su=rv\in \sigma .
  $$
  Notice that $v=1$, since otherwise $r\in \interior \sigma $.  Therefore $su=r$, so $s$ divides $r$, whence $s\in
\delta _r$.  To see that $r\notin rS$, notice that otherwise there would be some $t$ in $S$ such that $r=rt$, and again
this will conflict with the choice of $r$ outside $\interior \sigma $.

To show that $r$ is unique let us assume that $\delta _r=\delta _{r'}$.  Then $r$ and $r'$ divide each other, so there
exist $u$ and $v$ in $\tS $, such that $ru=r'$, and $r'v=r$.  From this we get that $r = ruv$, and then necessarily
$u=v=1$, or otherwise $r\in rS$.  Thus $r=r'$.

Conversely, if $\sigma =\delta _r$, with $r\notin rS$, we claim that $r\in \delta _r\setminus \interior \delta _r$.  On
the one hand it is evident that $r\in \delta _r$. On the other hand, supposing by contradiction that $r\in \interior
\delta _r$, there exists $p$ in $S$ such that $rp\in \delta _r$, whence $rp\|r$, so we may find $u$ in $\tS $ with
$rpu=r$.  This implies that $r\in rS$, contradicting the assumptions.  Therefore $\delta _r\neq \interior \delta _r$, as
desired.

\itmproof (iii) Given that $s_1$ and $s_2$ are in $\interior \sigma $, choose $p_1$ and $p_2$ in $S$, such that $s_1p_1$
and $s_2p_2$ belong to $\sigma $.  Using \ref {StringDirected}, we furthermore choose $u_1$ and $u_2$ in $\tS $, such
that
  $$
  s:= s_1p_1u_1=s_2p_2u_2\in \sigma .
  $$
  It then follows that
  $$
  s\in s_1S\cap s_2S = rS,
  $$
  so there is some $t$ in $S$ such that $s=rt$, so $r$ divides $s$ and we deduce that $r\in \sigma $.

\itmproof (iv)
  It is easy to see that $\interior \sigma $ satisfies \ref {DefineStrin.i-ii}, while \ref {StringDirected} follows
immediately from (iii) and the existence of least common multiples.  \endProof

\state Proposition \label MaxOpenOrDeadEnd
  Let $\sigma $ be a maximal string.  Then either $\sigma $ is open, or $\sigma =\delta _r$ for some $r$ in $S$, such
that $rS=\{0\}$.

\Proof
  Supposing that the maximal string $\sigma $ is not open, we have that $\sigma =\delta _r$, with $r\notin rS$, by \ref
{PropsOpenString.ii}, so it suffices to prove that $rS=\{0\}$.  Supposing otherwise, let $s\in S$ be such that
  $$
  t:= rs\neq 0.
  $$
  Since $r\|t$, we have that $\delta _r\subseteq \delta _t$, so $\delta _r=\delta _t$ by maximality.  It follows that
$t\|r$, so we may find $u$ in $\tS $ such that $tu=r$.  Therefore
  $$
  r = tu = rsu \in rS,
  $$
  contradicting the fact that $r$ is not in $rS$. This concludes the proof.
  \endProof

Let us now study open strings in relation to the representation $\theta ^\star $.

\state Proposition \label OpenInvarUnderTheta
  Let $\sigma $ be an open string in $S$, and let $r\in S$.
  \iItemize
  \iItem If $\sigma \in F^\star _r$, then $\theta ^\star _r(\sigma )$ is open.
  \iItem If $\sigma \in E^\star _r$, then $\theta ^{\star -1}_r(\sigma )$ is open.

\Proof (i) \enspace
  Let $t\in \theta ^\star _r(\sigma )$, so that $t\|rs$, for some $s$ in $\sigma $, and hence we may write $tx=rs$, for
a suitable $x$ in $\tS $.  Since $\sigma $ is open, we may pick some $p$ in $S$ such that $sp\in \sigma $, whence
  $$
  txp=rsp\in r\sigma \subseteq r*\sigma =\theta ^\star _r(\sigma ).
  $$
  Since $xp\in S$, this proves (i).

  \itmproof (ii) Let $t\in \theta ^{\star -1}_r(\sigma )$, so that $rt\in \sigma $.  Since $\sigma $ is open, we may
pick some $p$ in $S$ such that $rtp\in \sigma $.  Therefore
  $$
  tp\in r\inv *\sigma =\theta ^{\star -1}_r(\sigma ),
  $$
  proving (ii).  \endProof

We will return to the study of open strings in future sections.

\section Representing the inverse hull on strings

In \ref {DefRepre} we introduced the notion of semigroup representations on a set $\Omega $.  The semigroups we had in
mind there were simply associative semigroups with zero, but from now on we will also consider representations of
inverse semigroups, such as $\hull $.

There is no need to amend Definition \ref {DefRepre} when the semigroup considered is an inverse semigroup but it is
worth noticing that if $\rho $ is a representation of a given inverse semigroup $\S $ on a set $\Omega $, then
  $$
  \rho (s\inv )=\rho (s)\inv , \for s\in \S ,
  $$
  a fact that follows easily from the uniqueness of inverses.

Our goal now is to show that the action of a $0$-left cancellative semigroup $S$ on strings extends to the inverse hull
$\hull $.  We do this by first proving a more general result.

Let $X$ be a set equipped with a preorder, that is, a transitive and reflexive relation ``$\leq $".  Thus, $X$ is a
partially ordered set except for the fact that the anti-symmetric property is not required to hold.

The goal is to later allow the example in which $X$ is a semigroup and ``$\leq $" is the order given by \emph
{division}.

By abuse of language we will refer to ``$\leq $" as the \emph {order} on $X$, even though, strictly speaking, this is
not an order relation.

If $x$ and $y$ are in $X$, the \emph {interval from $x$ to $y$} is the set
  $$
  [x,y] = \{z\in X: x\leq z\leq y\}.
  $$
  We will also work with the \emph {unbounded intervals}
  $$
  (-\infty ,y] = \{z\in X: z\leq y\}, \and [x,+\infty ) = \{z\in X: x\leq z\}.
  $$

\definition \label OrderProps
  A subset $A\subseteq X$ is said to be:
  \iItemize
  \iItem \emph {convex}, if whenever $x,y\in A$, one has that $[x, y]\subseteq A$,
  \iItem \emph {hereditary}, if whenever $x\in A$, one has that $(-\infty , x]\subseteq A$,
  \iItem \emph {directed}, if $X$ is nonempty and whenever $x,y\in A$, there exists some $z$ in $A$ such that $x\leq z$
and $y\leq z$,
  \iItem a \emph {string}, if $A$ is hereditary and directed.
  \medskip \noindent The set of all directed subsets of $X$ will be denoted by $\Delta (X)$, and the set of all strings
in $X$ will be written $\Sigma (X)$.  Note that a set $D$ is directed if and only if each of its finite subsets
(including the empty set) has an upper bound.  What we call a string is usually called an ideal in order theory but
since the term ``ideal'' has other meanings in this paper we shall use the alternative terminology.  Convex sets are
precisely the sets which are the intersection of a hereditary set and the complement of a hereditary set.

For the sake of symmetry, observe that the notion of directed sets may also be expressed using intervals: $A\neq
\emptyset $ is directed if and only if, whenever $x,y\in A$, one has that
  $$
  [x,+\infty )\cap [y,+\infty )\cap A\neq \emptyset .
  $$

Observe that the empty set is convex and hereditary by vacuity.  It is however not directed since it has been explicitly
ruled out in \ref {OrderProps.iii}

\definition If $A$ is any subset of $X$, we will denote by $h(A)$ the \emph {hereditary closure} of $A$, namely
  $$
  h(A)=\{x\in X: \exists y\in A, \ x\leq y\} = \bigcup _{y\in A}(-\infty , y].
  $$

It is easy to see that $h(A)$ is the smallest hereditary subset of $X$ containing $A$.  Consequently $h(h(A))=h(A)$, and
$A$ is hereditary if and only if $h(A)=A$.

If $D$ is a directed subset of $X$, it is easy to see that $h(D)$ is also directed, and hence $h(D)$ is a string.  We
may therefore view $h$ as a map
  $$
  h\colon \Delta (X)\to \Sigma (X).
  \equationmark HSendsDelToSigma
  $$

If $A$ and $B$ are nonempty subsets of $X$, with $A\subseteq B$, recall that $A$ is said to be \emph {cofinal}\/ in $B$,
provided for every $x$ in $B$, there exists some $y$ in $A$, such that $x\leq y$.  Notice that if $B$ is nonempty, then
$A$ must be nonempty as well.

\state Proposition \label CofinDirec
  Let $A$ and $D$ be subsets of $X$, with $A\subseteq D$.  If $A$ is cofinal in $D$, and $D$ is directed, then $A$ is
also directed.

\Proof As observed above, $A\neq \emptyset $.  Given $x$ and $y$ in $A$, use that $D$ is directed to produce some $z$ in
$D$, such that $x,y\leq z$.  By cofinality, there is $w$ in $A$, with $z\leq w$, whence $x,y\leq w$.  \endProof

The following is a useful condition, similar to cofinality, that applies even when $A$ is not a subset of $B$.

\definition Given subsets $D$ and $E$ of $X$, with $D$ directed, we will say that $D$ is \emph {asymptotically
contained} in $E$, in symbols
  $$
  D\sqsubseteq E,
  $$
  provided $E\cap D$ is cofinal in $D$.

The above notion could easily be applied to a not necessarily directed set $D$, but we will have no use for it ouside
the situation outlined above.  Moreover, in most applications of that notion, $E$ will be a convex set.

The following will be useful later.

\state Lemma
  Let $D$ be a directed subset of $X$, and let $E\subseteq X$.
  \iItemize
  \iItem If $D\sqsubseteq E$, then $h(D)\sqsubseteq E$. \iItemmark ClosureSwept
  \iItem If $D\sqsubseteq E$, then $E\cap D$ is directed. \iItemmark InterDirec

\Proof
  (i) \ We need to show that $E\cap h(D)$ is cofinal in $h(D)$, so let $x\in h(D)$.  Then there exists $y$ in $D$ such
that $x\leq y$.  By hypothesis, there exists $z$ in $E\cap D$, such that $y\leq z$.  Therefore
  $$
  x\leq y\leq z\in E\cap D\subseteq E\cap h(D),
  $$
  as desired.

\itemproof {ii} follows immediately from \ref {CofinDirec}.
  \endProof

The relevance of asymptotical containment for strings is emphasized next.

\state Lemma \label ClosureIntersec
  Let $\sigma $ be a string in $X$, and let $E$ be any subset of $X$.  Then the following are equivalent:
  \iItemize
  \iItem $\sigma \sqsubseteq E$,
  \iItem $\sigma =h(E\cap \sigma )$,
  \iItem $\sigma =h(A)$, for some subset $A\subseteq E$.

\Proof (i)$\imply $(ii) \
  Since $\sigma $ is hereditary and $E\cap \sigma \subseteq \sigma $, it is clear that $h(E\cap \sigma )\subseteq \sigma
$.  To prove the reverse inclusion, pick any $x$ in $\sigma $.  Then by (i) there exists $y$ such that $x\leq y\in E\cap
\sigma $, so we see that $x\in h(E\cap \sigma )$.

\itemproof {ii)$\imply $(iii} Obvious.

\itemproof {iii)$\imply $(i} Since $A$ is evidently cofinal in $h(A)$, the hypothesis implies that $A$ is cofinal in
$\sigma $.  Since $A\subseteq E\cap \sigma $, we have that $E\cap \sigma $ is cofinal in $\sigma $, whence $\sigma
\sqsubseteq E$.
  \endProof

Recall that $\I (X)$ denotes the collection of all partial bijections on $X$.  It is well known that $\I (X)$ is an
inverse semigroup, known as the \emph {symmetric inverse semigroup} of $X$.

\definition Let $\varphi \in \I (X)$.
  \iItemize
  \iItem The domain of $\varphi $ will be denoted by $F_\varphi $ and its range will be written $E_\varphi $.
  \iItem We will say that $\varphi $ is \emph {order preserving} if, for every $x$ and $y$ in $F_\varphi $, one has that
  $$
  x\leq y \iff \varphi (x)\leq \varphi (y).
  $$
  \iItem We shall denote by $\I _+(X)$, the collection of all order preserving $\varphi $ in $\I (X)$, such that
$F_\varphi $ and $E_\varphi $ are convex.

\state Proposition
  $\I _+(X)$ is an inverse subsemigroup of $\I (X)$.

\Proof It is easy to see that if $\varphi $ is order preserving then so is $\varphi \inv $, so $\I _+(X)$ is seen to be
invariant under taking inverses.  Given $\varphi $ and $\psi $ in $\I _+(X)$, recall that the domain and range of the
composition $\varphi \psi $ are given respectively by
  $$
  F_{\varphi \psi } = \psi \inv (F_\varphi \cap E_\psi ), \and E_{\varphi \psi } = \varphi (F_\varphi \cap E_\psi ).
  $$
  We will next prove that $F_{\varphi \psi }$ is convex.  For this, suppose that $x, y\in F_{\varphi \psi }$, and $z\in
[x, y]$.  Since $F_{\varphi \psi }\subseteq F_\psi $, and the latter is convex, we have that $z\in F_\psi $.  So
  $$
  \psi (x)\leq \psi (z)\leq \psi (y),
  $$
  and since $\psi (x)$ and $\psi (y)$ lie in the convex set $F_\varphi $, it follows that $\psi (z)\in F_\varphi \cap
E_\psi $, whence $z\in F_{\varphi \psi }$.

This proves that $F_{\varphi \psi }$ is convex and a similar argument shows that $E_{\varphi \psi }$ is also convex.
Trivially, $\varphi \psi $ is order preserving.  \endProof

Our main goal will be to describe a canonical action of $\I _+(X)$ on $\Sigma (X)$.  As an intermediate step, for each
$\varphi $ in $\I _+(X)$, we will build a partial mapping $\varphi ^\Delta $ on $\Delta (X)$.

\definition Given any subset $E\subseteq X$, we will denote by $E^\Delta $ the subset of $\Delta (X)$ consisting of all
directed subsets $D$ of $X$, such that $D\sqsubseteq E$.

  Since every nonempty directed subset of $X$ is clearly asymptotically contained in $X$, we have that
  $$
  X^\Delta =\Delta (X).
  $$

Given $\varphi $ in $\I _+(X)$, and given $D$ in $F^\Delta _\varphi $, we have that $F_\varphi \cap D$ is directed by
\ref {InterDirec}.  Since $\varphi $ is an order-isomorphism from $F_\varphi $ to $E_\varphi $, it follows that $\varphi
(F_\varphi \cap D)$ is also directed.
  In addition, since $\varphi (F_\varphi \cap D)$ is contained in $E_\varphi $, it is obvious that $\varphi (F_\varphi
\cap D)\sqsubseteq E_\varphi $, meaning that $\varphi (F_\varphi \cap D)\in E^\Delta _\varphi $.  We therefore have a
well defined mapping
  $$
  D\in F^\Delta _\varphi \ \longmapsto \ \varphi (F_\varphi \cap D) \in E^\Delta _\varphi .
  $$

\definition Given any $\varphi $ in $\I _+(X)$, we will denote by $\varphi ^\Delta $ the above map from $F^\Delta
_\varphi $ to $E^\Delta _\varphi $, namely
  $$
  \varphi ^\Delta (D) = \varphi (F_\varphi \cap D), \for D\in F^\Delta _\varphi .
  $$ We view $\varphi ^\Delta $ as a partial mapping on $\Delta (X)$.

If $\psi =\varphi \inv $, we may ask ourselves what is the relationship between $\psi ^\Delta $ and $\varphi ^\Delta $,
but one should not expect these maps to be the inverse of each other.  The reason is that
  $$
  \psi ^\Delta \big (\varphi ^\Delta (D)\big ) = F_\varphi \cap D, \for D \in F^\Delta _\varphi ,
  $$
  and so when $D$ is asymptotically contained in $F_\varphi $, without being a subset of $F_\varphi $, one would have
that $\psi ^\Delta (\varphi ^\Delta (D)) \neq D$.
  By turning our attention to strings, rather than directed sets, we will soon fix this anomaly.
  Nonetheless, the assignment $\varphi \mapsto \varphi ^\Delta $ is a homomorphism.

\state Proposition \label Ishom For $\varphi ,\psi \in \I _+(X)$, one has $\varphi ^\Delta \psi ^\Delta =(\varphi \psi
)^\Delta $.

\Proof Suppose that $D\in F_{\varphi \psi }^\Delta $.  Then $D\cap F_{\varphi \psi }$ is cofinal in $D$.  So if $x\in
D$, then there is $y\in D\cap F_{\varphi \psi }\subseteq F_{\psi }$ with $x\leq y$.  Thus $D\in F_{\psi }^\Delta $ and
$\psi ^\Delta (D) = \psi (D\cap F_\psi )$. If $z\in \psi (D\cap F_{\psi })$ and $z=\psi (d)$ with $d\in D$, then by
assumption, there is $d'\in D\cap F_{\varphi \psi }$ with $d\leq d'$.  Then $\psi (d')\geq \psi (d)=z$ and $\psi (d')\in
\psi (D\cap F_{\psi })\cap F_{\varphi }$.  Thus $\psi ^\Delta (D)\in F^\Delta _\varphi $ and $\varphi ^\Delta \psi
^\Delta (D) = \varphi (\psi (D\cap F_\psi )\cap F_\varphi ) = \varphi \psi (D\cap F_{\varphi \psi })$.

It remains to show that the domain of $\varphi ^\Delta \psi ^\Delta $ is contained $F^\Delta _{\varphi \psi }$.  If
$D\in F^\Delta _\psi $ and $\psi ^\Delta (D)\in F^\Delta _{\varphi }$, then $D\cap F_\psi $ is cofinal in $D$ and $\psi
(D\cap F_\psi )\cap F_\varphi $ is cofinal in $\psi (D\cap F_\psi )$.  So if $x\in D$, then there exists $y\in D\cap
F_\psi $ with $x\leq y$. Then $\psi (y)\leq z$ with $z\in \psi (D\cap F_\psi )\cap F_\varphi $.  So $z=\psi (d)$ with
$d\in D\cap F_\psi $.  Then $d\in D\cap F_{\varphi \psi }$ and $x\leq y\leq d$ (as $\psi (y)\leq \psi (d)$ implies
$y\leq d$).  This completes the proof.  \endProof

\definition Given any subset $E$ of $X$, we will denote by $E^\Sigma $ the collection of all strings $\sigma $ in $X$
such that $\sigma \sqsubseteq E$.  Equivalently
  $$
  E^\Sigma =\Sigma (X)\cap E^\Delta .
  $$

   Since every string in $X$ is clearly asymptotically contained in $X$, we have that $X^\Sigma =\Sigma (X)$.

\state Proposition
  For every subset $E$ of $X$, one has that $h$ maps $E^\Delta $ into $E^\Sigma $.

\Proof
  Given $D$ in $E^\Delta $, we have that $h(D)\sqsubseteq E$, by \ref {ClosureSwept}.  Since $h(D)$ is a string by \ref
{HSendsDelToSigma}, we conclude that $h(D)\in E^\Sigma $.  \endProof

\definition Given any $\varphi $ in $\I _+(X)$, consider the map $\varphi ^\Sigma $ from $F^\Sigma _\varphi $ to
$E^\Sigma _\varphi $ given by $\varphi ^\Sigma =h\circ \varphi ^\Delta |_{F^\Sigma _\varphi }$.  Diagramatically:
  $$
  \matrix {
    & \varphi ^\Sigma \cr
    F^\Sigma _\varphi & \longrightarrow & E^\Sigma _\varphi \cr \vrule height 14pt depth 10pt width 0pt
    \raise 7pt
    \hbox {\rotatebox {-90}{$\subseteq $}}
    && \ \uparrow h\cr
    F^\Delta _\varphi & \longrightarrow & E^\Delta _\varphi \cr
    & \varphi ^\Delta
  }
  $$

Our main goal is to prove that the correspondence $\varphi \to \varphi ^\Sigma $ is an action of $\I _+(X)$ on $\Sigma
(X)$ by partial bijections.

The slightly complex definition of $\varphi ^\Sigma $ tends to cause formulas to quickly grow in size.  The following
technical fact is designed to contain the buildup of formulas by knocking off an ``$h$''.

\state Lemma \label KnockOffh
  If $\varphi \in \I _+(X)$, then $h\varphi ^\Delta =h\varphi ^\Delta h$ as partial mappings.

\Proof Suppose that $D\sqsubseteq F_{\varphi }$.  Then $h(D)\sqsubseteq F_{\varphi }^\Delta $ by \ref {ClosureSwept} and
so the domain of $h\varphi ^\Delta $ is contained in the domain of $h\varphi ^\Delta h$.  Suppose that $D$ is in the
domain of $h\varphi ^\Delta h$.  Then $h(D)\in F_\varphi ^\Delta $ and so $h(D)\cap F_{\varphi }$ is cofinal in $h(D)$.
Let $x\in D$.  Then since $x\in h(D)$, there is $y\in h(D)\cap F_{\varphi }$ with $x\leq y$.  Then $y\leq z$ with $z\in
D$ by definition of $h(D)$.  Then $z\leq w$ with $w\in h(D)\cap F_{\varphi }$, again by cofinality of $h(D)\cap
F_\varphi $ in $h(D)$.  As $y\leq z\leq w$ and $y,w\in F_{\varphi }$, we have $z\in F_{\varphi }$ by convexity of
$F_{\varphi }$.  We conclude that $x\leq z\in D\cap F_{\varphi }$ and so $D\in F^\Delta _{\varphi }$.  Thus the partial
mappings $h\varphi ^\Delta $ and $h\varphi ^\Delta h$ have the same domain.

Now we prove if $D\in F_{\varphi }^\Delta $, then $$
  h\big (\varphi \big (F_\varphi \cap h(D)\big )\big ) = h\big (\varphi (F_\varphi \cap D)\big ).  $$ Since $h(D)$
contains $D$, the inclusion ``$\supseteq $" is evident.  In order to prove the reverse inclusion, let $x$ be any element
of $h\big (\varphi \big (F_\varphi \cap h(D)\big )\big )$.  Then, there exists some $y$ in $F_\varphi \cap h(D)$, such
that $x\leq \varphi (y)$.  We may then pick $z$ in $D$, with $y\leq z$.

Since $D\sqsubseteq F_\varphi $, we have that $F_\varphi \cap D$ is cofinal in $D$, so there exists $w$ such that $z\leq
w\in F_\varphi \cap D$.  It follows that $y\leq w$, and since both $y$ and $w$ lie in $F_\varphi $, we have
  $$
  x\leq \varphi (y)\leq \varphi (w) \in \varphi (F_\varphi \cap D),
  $$
  whence $x\in h\big (\varphi (F_\varphi \cap D)\big )$.  \endProof

\state Proposition \label Inverse
  For every $\varphi $ in $\I _+(X)$, one has that $\varphi ^\Sigma $ is a bijective mapping from $F^\Sigma _\varphi $
to $E^\Sigma _\varphi $, and its inverse is given by $(\varphi \inv )^\Sigma $.

\Proof Initially notice that $F_{\varphi \inv }=E_\varphi $, and $E_{\varphi \inv }=F_\varphi $, whence $(\varphi \inv
)^\Sigma $ is a map from $E^\Sigma _\varphi $ to $F^\Sigma _\varphi $, as expected.  Next observe that $h(\varphi \inv
)^\Delta h\varphi ^\Delta =h(\varphi \inv )^\Delta \varphi ^\Delta =h(\varphi \inv \varphi )^\Delta =h(1_{F_\varphi
})^\Delta $ by \ref {KnockOffh} and \ref {Ishom}.  But if $\sigma $ is a string in $F^\Sigma _\varphi $, then
$h(1_{F_\varphi })^\Delta (\sigma ) = h(\sigma \cap F_{\varphi })=\sigma $ where the last equality is from \ref
{ClosureIntersec}.  It follows that $(\varphi \inv )^\Sigma \varphi ^\Sigma = 1_{F^\Sigma _{\varphi }}$.  A dual
argument completes the proof.  \endProof

The result above shows that $\varphi ^\Sigma $ is indeed an element of the inverse semigroup $\I (\Sigma (X))$, and we
next plan to prove that the correspondence
  $$
  \varphi \in \I _+(X)\ \mapsto \ \varphi ^\Sigma \in \I (\Sigma (X)),
  $$
  is an inverse semigroup homomorphism.

\state Proposition \label Compose
  For every $\varphi $ and $\psi $ in $\I _+(X)$, one has that $(\varphi \psi )^\Sigma =\varphi ^\Sigma \psi ^\Sigma $.
That is, the assignment $\I _+(X)\to \I (\Sigma (X))$ given by $\varphi \mapsto \varphi ^\Sigma $ is a homomorphism of
inverse semigroups.

\Proof We compute
  $$
  \varphi ^\Sigma \psi ^\Sigma =
  (h\varphi ^\Delta )|_{\Sigma (X)} (h\psi ^\Delta )|_{\Sigma (X)} =
  (h\varphi ^\Delta h\psi ^\Delta )|_{\Sigma (X)} \quebra =
  (h\varphi ^\Delta \psi ^\Delta )|_{\Sigma (X)} =
  (h(\varphi \psi )^\Delta )|_{\Sigma (X)} =
  (\varphi \psi )^\Sigma
  $$
  by \ref {KnockOffh} and \ref {Ishom}.  This completes the proof.  \endProof

Our main result is now an easy consequence of \ref {Inverse} and \ref {Compose}:

\state Theorem \label TheRep
  Let $X$ be a set equipped with a transitive and reflexive relation ``$\leq $".  Then there exists a semigroup
homomorphism
  $$
  \varphi \in \I _+(X)\mapsto \varphi ^\Sigma \in \I \big (\Sigma (X)\big )
  $$
  such that, for all $\varphi $ in $\I _+(X)$, and for every string $\sigma $ in $F^\Sigma _\varphi $, one has that
  $$
  \varphi ^\Sigma (\sigma )=h\big (\varphi (F_\varphi \cap \sigma )\big ).
  $$

As an application, let $S$ be a $0$-left cancellative semigroup, and let $S'=S\setminus \{0\}$.  Given $s$ and $t$ in
$S'$, recall that $s$ is said to divide $t$, in symbols $s\,|\,t$, if there exists $u$ in $\tilde S = S\cup \{1\}$ such
that $su=t$.  Setting
  $$
  s\leq t \iff s\,|\,t,
  $$
  we have that $S'$ becomes a (possibly not anti-symmetric) ordered set.  For each $s$ in $S$, recall that
  $$
  F_s=\{x\in S':sx\neq 0\}, \and E_s=sS\setminus \{0\},
  $$
  and that
  $$
  \theta _s\colon x\in F_s\mapsto sx \in E_s.
  $$

\state Proposition For every $s$ in $S$, one has that $F_s$ and $E_s$ are convex in $S'$, and $\theta _s$ is order
preserving.

\Proof
  Given $x,y,z\in S'$, with $x\leq y\leq z$, and $x,z\in F_s$, pick $u$ in $\tilde S$ such that $yu=z$.  Then
  $$
  syu = sz \neq 0,
  $$
  so $sy\neq 0$, whence $y\in F_s$.  This shows that $F_s$ is convex.  Notice that $x$ did not play any role above,
which means that $F_s$ is in fact a hereditary set.

Now let $x,y,z\in S'$, with $x\leq y\leq z$, and $x,z\in E_s$.  We may then pick $u$ in $\tilde S$ such that $xu=y$.  We
may also write $x=st$, for some $t$ in $S$, so
  $$
  y = xu = stu,
  $$
  whence $y\in E_s$.  This shows that $E_s$ is convex.  Notice that $z$ did not play any role above, which means that
$E_s$ is in fact a hereditary set for the reverse order.

Given $x$ and $y$ in $F_s$, with $x\leq y$, pick $u$ in $\tilde S$ such that $xu=y$.  Then
  $$
  \theta _s(x)u=sxu = sy = \theta _s(y),
  $$
  whence $\theta _s(x)\leq \theta _s(y)$.

Instead of assuming that $x\leq y$, suppose that $\theta _s(x)\leq \theta _s(y)$, so we may find $u$ in $\tilde S$, such
that $\theta _s(x)u=\theta _s(y)$, which is to say that $sxu=sy$, whence $xu=y$ by virtue of $0$-left cancellativity,
proving that $x\leq y$.
  This concludes the proof that $\theta _s$ is order preserving.  \endProof

The result above implies that each $\theta _s$ lies in $\I _+(S')$, and hence the inverse semigroup generated by the
$\theta _s$, namely $\hull $, is a subset of $\I _+(S')$.

We may then view the composition
  $$
  \rho \ \colon \ \hull \hookrightarrow \I _+(S')\to \I \big (\Sigma (S')\big ),
  $$
  where the rightmost arrow is the representation defined by \ref {TheRep},
  as a representation of $\hull $ on $\Sigma (S')$.

In the context of $0$-left cancellative semigroups notice that the notion of strings, as introduced in \ref
{DefineStrin}, coincides with the concept defined in \ref {OrderProps.iv} for the above order relation on $S'$.  In
other words,
  $$
  S^\star =\Sigma (S'),
  $$
  so $\rho $ is seen to be a representation of $\hull $ on $S^\star $.

\state Corollary Regarding the representation $\rho $ above, one has that
  $$
  \rho (\theta _s) = \theta ^*_s,
  $$
  for every $s$ in $S$.

\Proof
  We must first prove that $\theta ^*_s$ and $\rho (\theta _s)$ share domains an ranges.  Recall that the domain of
$\theta _s^*$ is given by
  $$
  F^\star _s=\{\sigma \in S^\star : s\sigma \not \ni 0\},
  $$
  and its range is the set
  $$
  E^\star _s=\{\sigma \in S^\star : \sigma \cap sS \neq \emptyset \}.
  $$
  On the other hand, the domain and range of $\rho (\theta _s)=\theta _s^\Sigma $ are respectively given by $F_s^\Sigma
$ and $E_s^\Sigma $.

  In order to prove that $F^\star _s=F_s^\Sigma $, let $\sigma $ be a string in $F^\star _s$.  Then clearly $\sigma
\subseteq F_s$, whence $\sigma \sqsubseteq F_s$, and we see that $\sigma \in F_s^\Sigma $.  Conversely, if $\sigma \in
F_s^\Sigma $, suppose by contradiction that there exists $x$ in $\sigma $ such that $sx=0$.  By assumption there exists
some $y$ in $F_s\cap \sigma $, such that $x\leq y$, so we may find $u$ in $\tilde S$, such that $xu=y$.  Therefore
  $$
  sy = sxu = 0,
  $$
  contradicting the fact that $y\in F_s$.  This proves that $\sigma $ lies in $F^\star _s$, and hence that $F^\star
_s=F_s^\Sigma $.

We will next prove that $E^\star _s=E_s^\Sigma $, so pick any $\sigma $ in $E^\star _s$.  Then there exists some $x$ in
$\sigma $ of the form $x=st$, with $t\in S$.  Given any $y$ in $\sigma $, we may use the fact that $\sigma $ is directed
to find some $z$ in $\sigma $, such that $x,y\leq z$.  Therefore there exists $u$ in $\tilde S$, such that $xu=z$, so
  $$
  z=xu=stu\in E_s\cap \sigma .
  $$
  Since $y\leq z$, this proves that $E_s\cap \sigma $ is cofinal in $\sigma $, whence $\sigma \in E_s^\Sigma $.  This
proves that $E^\star _s\subseteq E_s^\Sigma $.  Conversely, given $\sigma $ in $E_s^\Sigma $, we have that $E_s\cap
\sigma $ is cofinal in $\sigma $, so in particular $E_s\cap \sigma $ is nonempty, and this clearly implies that $\sigma
\in E^\star _s$.  This concludes the proof that $E^\star _s=E_s^\Sigma $.

Given any $\sigma \in F_s^\Sigma =F^\star _s$, we have seen that $\sigma \subseteq F_s$, whence
  $$
  \rho (\theta _s)(\sigma ) =
  \theta _s^\Sigma (\sigma ) =
  h\big (\theta _s(F_s\cap \sigma )\big ) =
  h\big (\theta _s(\sigma )\big ) =
  \{t\in S': t\leq x, \hbox { for some } x\in s\sigma \} \quebra =
  \{t\in S': t\,|\, sr, \hbox { for some } r\in \sigma \} =
  \theta ^*_s(\sigma ).
  $$
  This completes the proof.
  \endProof

The big conclusion of all this is as follows:

\state Proposition \label MapForISG
  Let $S$ be a $0$-left cancellative semigroup.
  Then there exists a unique representation $\rho $ of $\hull $ on $S^\star $,
  such that the following diagram commutes. \smallskip
  \hfill \beginpicture \setcoordinatesystem units <0.003truecm, -0.00250truecm>
  \put {$S$} at 0000 0000
  \put {$\I (S^\star )$} at 1100 0000
  \put {$\hull $} at 500 580
  \arrow <0.11cm> [0.3,1.2] from 150 000 to 850 000
  \put {$\theta ^\star $} at 500 -130 \arrow <0.11cm> [0.3,1.2] from 100 100 to 400 450
  \put {$\theta $} at 150 350
  \arrow <0.11cm> [0.3,1.2] from 600 450 to 900 100
  \put {$\rho $} at 830 350
  \endpicture \hfill \null

Observing that a homomorphism of inverse semigroups must restrict to the corresponding idempotent {\sla }s, we obtain
the following:

\state Corollary \label IntroEpsilon
  Let $S$ be a $0$-left cancellative semigroup.
  Then there exists a {\sla } representation
  $$
  \varepsilon \colon \ehull \to \P (S^\star ),
  $$
  \# $\varepsilon $; Representation of $\ehull $ on $S^\star $;
  such that
  $$
  \varepsilon \big (\theta _u(F^\theta _\Lambda )\big )=\theta ^\star _u(F^\star _\Lambda ),
  $$
  whenever $\Lambda $ is a finite subset of $\tS $ intersecting $S$ nontrivially, and $u\in \Lambda $.

\Proof
  Identifying the idempotent {\sla }s of $\hull $ and $\I (S^\star )$ with $\ehull $ and $\P (S^\star )$, respectively,
it is enough to take $\varepsilon $ to be the restriction of the representation $\rho $ of \ref {MapForISG} to $\ehull
$.
  \endProof

Observing that $E^\theta _r=\theta _r(F^\theta _r)$, notice that
  $$
  \varepsilon (E^\theta _r) = \theta ^\star _r(F^\star _r)=E^\star _r.
  $$

The subset of $S^\star $ consisting of all open strings was shown in \ref {OpenInvarUnderTheta} to be invariant under
both $\theta ^\star _r$ and $\theta ^{\star -1}_r$, for every $r$ in $S$.  As an immediate consequence we thus obtain
the following:

\state Proposition \label OpenStringsInvariant
  The set of all open strings in $S$ is invariant under the representation of $\hull $ on $S^\star $ described in \ref
{MapForISG}.

We now prove that the mapping $r\mapsto \delta _r$ from $S'=S\setminus \{0\}$ to $\Sigma (S)$ is covariant with respect
to $\rho $.  More precisely, we prove the following.

\state Theorem \label BiggerCovariance Let $\varphi \in \hull $ and $r\in S'=S\setminus \{0\}$.  Then $r\in F_{\varphi
}$ if and only if $\delta _r\in F^\Sigma _{\varphi }$.  Moreover, if $r\in F_{\varphi }$, then $\delta _{\varphi (r)} =
\varphi ^\Sigma (\delta _r)$.

\Proof Suppose first that $r\in F_{\varphi }$.  Thence since $r$ is the maximum element of $\delta _r$, clearly $\delta
_r\sqsubseteq F_{\varphi }$ and so $\delta _r\in F^\Sigma _{\varphi }$.  Conversely, if $\delta _r\in F^\Sigma _{\varphi
}$, then $\delta _r\cap F_{\varphi }$ is cofinal in $\delta _r$.  But then there exists $s\in \delta _r\cap F_{\varphi
}$ such that $r\leq s$.  But also, $s\leq r$.  It follows since $F_{\varphi }$ is convex that $r\in F_{\varphi }$.

Assume that $r\in F_{\varphi }$.  By definition $\varphi ^\Sigma (\delta _r) = h(\varphi (F_{\varphi }\cap \delta _r))$.
So $\varphi (r)\in \varphi ^\Sigma (\delta _r)$, whence $\delta _{\varphi (r)}\subseteq \varphi ^\Sigma (\delta _r)$.
If $s\in \varphi ^\Sigma (\delta _r)$, then $s\leq \varphi (t)$ with $t\in \delta _r\cap F_{\varphi }$.  Since $\varphi
$ is order preserving, $s\leq \varphi (t)\leq \varphi (r)$ and so $s\in \delta _{\varphi (r)}$.  This completes the
proof.  \endProof

\section Unbounded strings and backward invariance

Recall from \ref {BackInvarianceCorol} that, for semigroups which are categorical at zero, the space $S^\infty $ of
maximal strings is invariant under every ${\theta ^\star _r}\inv $.  However, even though semigroups obtained from
subshifts (see \ref {DefineSubshSgrp}) are not necessarily categorical at zero, the above invariance may be shown to
hold.

Here we would like to present a general result about invariance of $S^\infty $ under ${\theta ^\star _r}\inv $, which
does not rely on the property of being categorical at zero and hence applies to subshift semigroups.

\fix Let $S$ be a $0$-left cancellative semigroup, fixed throughout this section.

\definition \label LengtFun
  Let $N$ be a totally ordered set.  An \emph {$N$-valued length function} for $S$ is a function
  $$
  \ell \colon S'=S\setminus \{0\}\to N,
  $$
  \# $\ell $; Length function;
  such that for every $r,s,t$ in $S'$, one has that
  \iItemize
  \iItem if $s\|t$, then $\ell (s)\leq \ell (t)$,
  \iItem if $r\|st\neq 0$, and $\ell (r)\leq \ell (s)$, then $r \|s$.

For an example consider the semigroup $S=L\cup \{0\}$, given by a language $L$ invariant under prefixes and suffixes, as
in \ref {LanguageSemigroup}.  Letting
  $$
  \ell \colon S'=L\to {\bf N}
  $$
  be defined by setting $\ell (s)$ equal to the usual word-length of $s$, it is easy to see that $\ell $ is a length
function in the sense of \ref {LengtFun}.

\definition \label BoundedSets
  Let $\ell $ be an $N$-valued length function on $S$.
  \iItemize
  \iItem We will say that a subset $X\subseteq S'$ is bounded
  (relative to $\ell $), provided there exists $n_0\in N$, such that $\ell (s)\leq n_0$, for all $s$ in $X$.
  \iItem We will say that $\ell $ is \emph {homogeneous}
  if, for every $r$ in $S$, and for every bounded subset $X\subseteq F^\theta _r$, one has that $rX$ is bounded.

If $\ell $ is the word-length mentioned above, one has that
  $$
  \ell (rs)=\ell (r)+\ell (s),
  $$
  whenever $rs\neq 0$, and from this it easily follows that $\ell $ is homogeneous.

\state Lemma \label MaxUnbdd
  Let $S$ be a $0$-left cancellative semigroup equipped with a homogeneous length function
  $
  \ell \colon S' \to N.
  $
  Given $r$ in $S$, let $\sigma $ be an unbounded string in $E^\star _r$.  Then
  \iItemize
  \iItem ${\theta ^\star _r}\inv (\sigma )$ is unbounded,
  \iItem if $\mu $ is a string with ${\theta ^\star _r}\inv (\sigma ) \subseteq \mu $, then $\mu \in F^\star _r$,
  \iItem if $\sigma $ is maximal, then ${\theta ^\star _r}\inv (\sigma )$ is also maximal.

\Proof
  (i) Supposing by contradiction that
  $$
  \tau := {\theta ^\star _r}\inv (\sigma ) = r\inv *\sigma
  $$
  is bounded, then $r\tau $ is also bounded by homogeneity.  Noticing that $\sigma = r*\tau $, by \ref
{IntroStarAction}, and hence that $\sigma $ is is the hereditary closure of $r\tau $, it would follow that $\sigma $ is
also bounded, a contradiction.  This shows that ${\theta ^\star _r}\inv (\sigma )$ is unbounded.

\medskip \noindent (ii)
  Arguing again by contradiction, suppose that
  $$
  {\theta ^\star _r}\inv (\sigma ) = r\inv *\sigma \subseteq \mu ,
  $$
  and that $\mu $ is not in $F^\star _r$, so there exists some $x\in \mu $, with $rx=0$.

By (i) we have that $r\inv *\sigma $ is unbounded, so in particular $\ell (x)$ cannot be a bound for $r\inv *\sigma $.
Therefore there exists some $y\in r\inv *\sigma $, such that $\ell (y)\not \leq \ell (x)$, and because $N$ is totally
ordered this means that $\ell (y)>\ell (x)$.

Since $x,y\in \mu $, there are $u,v\in \tS $, such that
  $
  xu=yv\in \mu .
  $
  It follows that $x\| yv$, and hence by \ref {LengtFun.ii} we deduce that $x\|y$.  So $y=xw$, for some $w$ in $\tS $,
whence
  $$
  0 = rxw = ry \in \sigma ,
  $$
  a contradiction.  This proves that $\mu \in F^\star _r$.

\medskip \noindent (iii) Suppose that $\mu $ is a string with
  $$
  {\theta ^\star _r}\inv (\sigma )\subseteq \mu .
  $$
  We then have by (ii) that $\mu \in F^\star _r$, so
  $$
  \sigma =\theta ^\star _r\big ({\theta ^\star _r}\inv (\sigma )\big ) \subseteq
  \theta ^\star _r(\mu ),
  $$
  and then $\sigma = \theta ^\star _r(\mu )$, by maximality, whence
  $$
  {\theta ^\star _r}\inv (\sigma ) =
  {\theta ^\star _r}\inv \big (\theta ^\star _r(\mu )\big ) = \mu ,
  $$
  thus proving that ${\theta ^\star _r}\inv (\sigma )$ is maximal.
  \endProof

\state Corollary \label BackInvarCorol
  Let $S$ be a $0$-left cancellative semigroup, admitting a homogeneous length function relative to which every maximal
string is unbounded.  Then $S^\infty $ is an invariant subset of $S^\star $ under the natural action
  of \/ $\I (S^\star ,\theta ^\star )$.

\Proof
  Follows immediately from the above result as in \ref {BackInvarianceCorol}.
  \endProof

If $\X $ is a subshift, it is easy to see that the maximal strings in $S_\X ^\star $ are unbounded, so the result above
applies even though such semigroups are not always categorical at zero.

\part {TWO}{Semilattices}

\noindent Given a $0$-left cancellative semigroup $S$, we will now concentrate on studying the semillatice $\ehull $ of
constructible sets, focusing in particular on the question of determining its spectrum.

We begin by discussing some general aspects of {\sla s}.

\section Preliminaries on {\sla s}

Let $\E $ be a {\sla } with zero.  By a \emph {character} on $\E $ \cite [12.4]{actions} we mean any nonzero map
  $$
  \varphi \colon \E \to \{0, 1\},
  $$
  such that $\varphi (0)=0$, and
  $\varphi (xy)=\varphi (x)\varphi (y)$, for all $x$ and $y$ in $\E $.  The set of all characters on $\E $, usually
  denoted\fn {The spectrum of $\E $ is however denoted by $\hat \E _0$ in \cite {actions}, which incidentally is our
basic reference for the theory of {\sla s}.}
  by $\hat \E \pilar {12pt}$, is called the \emph {spectrum} of $\E $.  It is well known that $\hat \E $ is a locally
compact, Hausdorff topological space when equipped with the product topology induced from $\{0, 1\}^\E $.

The topology on $\hat \E $ therefore admits a basis of open sets of the form
  $$
  D(f_1)\cap \cdots \cap D(f_m) \cap D(e_1)^c\cap \cdots \cap D(e_n)^c,
  \equationmark BasicOpen
  $$
  with $f_1, \ldots f_m, e_1,\ldots , e_n\in \E $, where
  $$
  D(e):=\{\varphi \in \E : \varphi (e)=1\},
  \equationmark SubBasicOpen
  $$
  and $X^c$ is the complement of $X$ in $\hat \E $.  Setting $e =f_1 \ldots f_m$, and observing that
  $$
  \varphi (f_1)=\cdots =\varphi (f_m)=1 \IFF \varphi (e)=1,
  $$
  we see that the set in \ref {BasicOpen} coincides with
  $$
  D(e) \cap D(e_1)^c\cap \cdots \cap D(e_n)^c,
  \equationmark BasicOpenBis
  $$
  so we may take the above to be the general form of a basic open set.  Note that one may assume without loss of
generality that $e_i\leq e$, for $i=1,\ldots ,n$.

A \emph {filter} on $\E $ is by definition \cite [12.1 \& 12.2]{actions} a nonempty subset $\xi \subseteq \E $, such
that
  \medskip \item {$\bullet $} $0\notin \xi $,
  \medskip \item {$\bullet $} $x,y\in \xi \IMPLY xy\in \xi $
  \medskip \item {$\bullet $} $x\geq y\in \xi \IMPLY x\in \xi $,
  \medskip \noindent
  for every $x$ and $y$ in $\E $.

Recall from \cite {actions} that there is a one-to-one correspondence between
  characters and filters, as follows: given a character $\varphi $, the corresponding filter is given by
  $$
  \xi =\{x\in \E : \varphi (x)=1\}.
  $$
  On the other hand, given a filter $\xi $, the corresponding character is given by
  $$
  \varphi (x) = \bool {x\in \xi }, \for x\in \E ,
  $$
  where brackets stand for Boolean value.

By definition an \emph {\uf } is a filter which is not properly contained in any other filter.  If a character $\varphi
$ corresponds to an {\uf }, then $\varphi $ is called an \emph {ultracharacter}.  The set of all ultracharacters is
denoted $\hat \E _\infty $.

It is well known \cite [12.3]{actions} that a filter $\xi $ is an ultrafilter if and only if $e\notin \xi $ implies that
$ef=0$ for some $f\in \xi $.  Also, every filter is contained in an ultrafilter.

The topology on the set of filters, inherited from the product topology via the above correspondence of filters and
characters, still has basic sets of the form \ref {BasicOpenBis}, except that now we should interpret \ref
{SubBasicOpen} as
  $$
  D(e):=\{\xi : e\in \xi \}.
  $$

If $\xi $ is an ultrafilter, then its neighborhoods of the form $D(e)$, with $e\in \xi $, form a neighborhood base for
$\xi $.  Indeed, if
  $$
  \xi \in D(e)\cap D(e_1)^c\cap \dots \cap D(e_n)^c,
  $$
  then $e_1,\ldots ,e_n\notin \xi $, so we can find $z_1,\ldots , z_n\in \xi $, with $z_ie_i=0$.  Then $e'=ez_1\cdots
z_n\in \xi $ and
  $$
  \xi \in D(e')\subseteq D(e)\cap D(e_1)^c\cap \dots \cap D(e_n)^c.
  $$

  Given any subset $F\subseteq \E $, we shall say that a subset $Z\subseteq F$ is a \emph {cover for} $F$ \cite
[11.5]{actions} if, for every nonzero $x\in F$, there exists $z\in Z$ such that $z x\neq 0$.  If $y\in \E $ and $Z$ is a
cover for $F:=\{x\in \E : x\leq y\}$, we will say that $Z$ is a \emph {cover for} $y$.

  Given finite subsets $X,Y\subseteq \E $, we shall denote by $\E ^{X,Y}$ the subset of\/ $\E $ given by
  $$
  \E ^{X,Y} = \{z\in \E : z\leq x,\ \forall x\in X,\hbox { and } zy=0,\ \forall y\in Y\}.
  $$

A character $\varphi $ of $\E $ is said to be \emph {tight} \cite [11.6]{actions} if, for all finite subsets
$X,Y\subseteq \E $, and for every finite cover $Z$ for
  $\E ^{X,Y}$, one has that
  $$
  \bigvee _{z\in Z}\varphi (z) =
  \bigwedge _{x\in X} \varphi (x) \wedge \bigwedge _{y\in Y} \neg {\varphi (y)}.
  $$

In view of the fact that characters are nonzero by definition, and hence satisfy \cite [11.7.(i)]{actions}, one has by
\cite [11.8]{actions} that a character $\varphi $ is tight if and only if, for every $x\in \E $ and for every finite
  cover $Z$ for $x$, one has that
  $$
  \bigvee _{z\in Z}\varphi (z) = \varphi (x).
  $$

Every ultracharacter is necessarily tight \cite [12.7]{actions}, and in fact the set $\hat \E _\tight $ formed by the
tight characters coincides with the closure of $\hat \E _\infty $ in $\hat \E $ \cite [12.9]{actions}.

If $\xi _1$ and $\xi _2$ are two filters in $\E $, and if $\varphi _1$ and $\varphi _2$ are the corresponding
characters, observe that
  $$
  \xi _1\subseteq \xi _2 \iff \varphi _1\leq \varphi _2.
  \equationmark InclusionVsLess
  $$
  We note that characters are functions taking values in the ordered set $\{0,1\}$, and that the order among characters
mentioned above is supposed to mean pointwise order.

Based on this one may give a characterization of ultracharacters which does not explicitly mention the associated
filters:

\state Proposition \label UltraCharCharac
  Let $\E $ be a {\sla } and let $\varphi \in \hat \E $.  Then $\varphi $ is an ultracharacter if and only if,
  $$
  \forall \psi \in \hat \E ,\ \varphi \leq \psi \Imply \varphi =\psi .
  $$

\Proof
  Follows immediately from \ref {InclusionVsLess}.
  \endProof

\def \fgen #1{\langle #1\rangle \kern -1pt_{\scriptscriptstyle E}} \def \ne {\fgen \eta }

  Let $\E $ be a {\sla } with zero and let $J$ be an ideal in $\E $ in the sense that $J\E \subseteq \E $.
  For every filter $\eta $ in $J$, consider the filter in $\E $ generated by $\eta $, namely
  $$
  \ne =\{x\in \E :\ \exists y\in \eta , \ y\leq x\}.
  $$
  On the other hand, for every filter $\xi $ in $\E $, such that $\xi \cap J\neq \emptyset $, it is evident that $\xi
\cap J$ is a filter in $J$.

\state Proposition \label PandI
  The correspondence
  $$
  i: \eta \in \hat J\mapsto \ne \in \hat \E ,
  $$
  is a homeomorphism
  from $\hat J$ onto the open subset of $\hat \E $ given by
  $$
  U = \{\xi \in \hat \E : \xi \cap J\neq \emptyset \}.
  $$
  Moreover, the inverse of the above correspondence is given by
  $$
  p:\xi \in U\mapsto \xi \cap J\in \hat J.
  $$ Furthermore, $i$ and $p$ are order isomorphisms where we order filters by inclusion.

\Proof
  Given any filter $\eta $ in $\hat J$, it is evident that
  $$
  \eta \subseteq \ne \cap J,
  \equationmark EtaInGen
  $$
  so $\ne \cap J$ is nonempty and hence we see that $\ne \in U$.

Observe that \ref {EtaInGen} is in fact an equality, since for every $x\in \ne \cap J$, there exists $y\in \eta $, with
$y\leq x$, so obviously $x\in \eta $.  The fact that $\eta =\ne \cap J$ may then be expressed as $\eta =p(i(\eta ))$, so
the composition $p\circ i$ is the identity on $\hat J$.

Given $\xi $ in $U$, it is easy to see that
  $$
  \fgen {\xi \cap J}\subseteq \xi ,
  \equationmark GenInXi
  $$
  and we again claim that this inclusion is an equality.  In fact, given any $x$ in $\xi $, choose any $y$ in the
nonempty set $\xi \cap J$, and notice that that
  $$
  x\geq xy \in \xi \cap J,
  $$
  so $x\in \fgen {\xi \cap J}$.  This proves that $\fgen {\xi \cap J}=\xi $ or, equivalently, that $i\circ p$ is the
identity on $U$.  Obviously, $i$ and $p$ are order preserving.

To see that $U$ is open, notice that we may write it as
  $$
  U = \bigcup _{x\in J}\{\xi \in \hat \E : x\in \xi \},
  $$
  which is a union of basic open subsets of $\hat \E $.

We check now that $i$ and $p$ are continuous.  We will subscript $D$ by the space we are working in.  If $\eta \in \hat
J$, then a basic neighborhood of $i(\eta )$ is of the form
  $$
  D_{\hat \E }(e)\cap D_{\hat \E }(e_1)^c\cdots \cap D_{\hat \E }(e_n)^c
  $$
  with $e,e_1,\ldots , e_n\in \E $.  Since $e\in i(\eta )$, there exists $f\in \eta $ with $f\leq e$.  Let $f_i=fe_i$;
note that $f_i\in J$ and $f_i\notin \eta $.  Then $\eta \in D_{\hat J}(f)\cap D_{\hat J}(f_1)^c\cap \cdots \cap D_{\hat
J}(f_n)^c$.  Moreover, if $\xi \in D_{\hat J}(f)\cap D_{\hat J}(f_1)^c\cap \cdots \cap D_{\hat J}(f_n)^c$, then $e\in
i(\xi )$, but $e_1,\ldots , e_n\notin i(\xi )$.  This shows that $i$ is continuous.

Next, if $\xi \in U$ and
  $$
  V=D_{\hat J}(e)\cap D_{\hat J}(e_1)^c\cap \cdots \cap D_{\hat J}(e_n)^c
  $$
  with $e,e_1,\ldots , e_n\in J$ is a basic neighborhood of $p(\xi )$ in $\hat J$, then $\xi \in D_{\hat \E }(e)\cap
D_{\hat \E }(e_1)^c\cap \cdots \cap D_{\hat \E }(e_n)^c$ and this neighborhood clearly maps into $V$ under $p$.  This
completes the proof that $i$ and $p$ are homeomorphisms.  \endProof

As mentioned earlier, filters are in one-to-one correspondence with characters.  Seeing things from the latter point of
view, one has:

\state Proposition \label SpecIdealNoVanish
  Identifying $\hat J$ with its image in $\hat \E $ under the map $i$ of \ref {PandI}, one has that
  $$
  \hat J = \{\varphi \in \hat \E : \varphi |_J \neq 0\}.
  $$

\Proof
  Given $\varphi $ in $\hat \E $, let $\xi $ be the associated filter, namely
  $$
  \xi =\{x\in \E : \varphi (x)=1\}.
  $$
  Thus, a necessary and sufficient for $\varphi |_J$ to be nonzero is that $\xi \cap J\neq \emptyset $, which is to say
that $\xi $ lies in the set $U$ of \ref {PandI}, namely the range of $i$.  This completes the proof.
  \endProof

\state Proposition Let $\eta $ be a filter in $J$.  Then
  $\eta $ is an ultrafilter if and only if $\fgen \eta $ is an ultrafilter in $\hat \E $.

\Proof Observe that the open set $U$ in \ref {PandI} is an upper set.  Hence a filter in $U$ is maximal in $U$ if and
only if it is an ultrafilter.  The proposition now follows because $i$ and $p$ in \ref {PandI} are order isomorphisms.
  \endProof

\section Representations of {\sla }s

Given a {\sla } $\E $, and and given any representation
  $$
  \pi \colon \E \to \I (\Omega ),
  $$
  on some set $\Omega $, it is easy to see that the range of $\rho $ must in fact be contained in the {\sla } $\P
(\Omega )$.  This motivates the following:

\definition \label DefRepOfSLA
  Let $\E $ be a {\sla } and let $\B $ be a Boolean
  algebra\fn {We do not require Boolean algebras to have a top element: for us they have meets, joins, a bottom element
and relative complements.}.
  A map $\pi \colon \E \to \B $ will be called a \emph {representation of $\E $ in $\B $}, provided $\pi (0)=0$, and
$\pi (xy)=\pi (x)\wedge \pi (y)$, for every $x$ and $y$ in $\E $.

In case $\B $ coincides with the Boolean algebra $\P (\Omega )$, for a given set $\Omega $, the above concept therefore
reduces to the notion of a representation of $\E $ on $\Omega $, as defined in \ref {DefRepre}.

If $\E $ is a sub{\sla } of $\P (\Omega )$ (always supposed to include the zero element of $\P (\Omega )$, namely the
empty set), then the inclusion map
  $$
  \E \hookrightarrow \P (\Omega )
  $$
  is evidently a representation of $\E $ on $\Omega $.

Suppose we are given a {\sla } $\E $, and a representation $\pi $ of $\E $ on a set $\Omega $.  For each $\omega $ in
$\Omega $, consider the mapping
  $
  \varphi ^\pi _\omega \colon \E \to \{0, 1\},
  $
  defined by
  $$
  \varphi ^\pi _\omega (x) = \bool {\omega \in \pi (x)}, \for x\in \E ,
  \equationmark CharacterFromRep
  $$
  where brackets stand for Boolean value.

It is clear that $\varphi ^\pi _\omega $ is a multiplicative map, so it is a character
  as long as it is nonzero.  It is moreover clear that $\varphi ^\pi _\omega $ is nonzero if and only if $\omega $ is in
the essential subset $\Omega _\ess $ for $\pi $.  Therefore
  $$
  \hat \E _\Omega := \{\varphi ^\pi _\omega :\omega \in \Omega \}\setminus \{0\} = \{\varphi ^\pi _\omega :\omega \in
\Omega _\ess \}
  \equationmark BigSubset
  $$
  is a subset of the spectrum $\hat \E $ of $\E $.

We would now like to discuss the question of how big is $\hat \E _\Omega $ within $\hat \E $.  The answer will of course
depend on $\pi $ since, when $\pi $ is the identically zero map, for instance, one should not expect $\hat \E _\Omega $
to be very big at all.
  Based on similar results pertaining to other representation theories, such as of C*-algebras,
  one might expect that $\hat \E _\Omega $ is dense in $\hat \E $ when $\pi $ is injective,
  but this is unfortunately not true.

For example, if $\Omega =\{1, 2, 3\}$, $\E =\P (\Omega )$, and $\pi $ is the identity map from $\E $ to $\P (\Omega )$,
then $\hat \E $ is finite and the character $\varphi $ defined by
  $$
  \varphi (X)=\bool {\{1, 2\}\subseteq X}, \for X\in \E ,
  $$
  is neither in $\hat \E _\Omega $, nor in its closure.

The next result will give a measure of the size of $\hat \E _\Omega $,
  provided $\pi $ does not send a nonzero element to the empty set.  Its proof is deceptively simple, yet the result is
significant.

\state Proposition \label ApproxTightChar
  Let $\pi $ be a representation of the {\sla } $\E $ on a set $\Omega $.  If $\pi (x)$ is nonempty for every nonzero
$x$ in $\E $, then
  the closure of\/ $\hat \E _\Omega $ in $\hat \E $ contains all tight characters of $\E $.

\Proof Since the ultracharacters are dense in the space of tight characters, it suffices to show that each neighborhood
of an ultracharacter $\psi $ intersects $\hat \E _\Omega $. As already mentioned, a basic neighborhood of $\psi $ is
given by $D(e)$, with $0\neq e\in \E $.  Let $\omega \in \pi (e)$.  Then $\varphi _{\omega }^{\pi }(e) =1$ and so
$\varphi _{\omega }^{\pi }\in D(e)$.  This completes the proof.  \endProof

As an immediate consequence we:

\state Corollary \label TautoCharsDense
  Let $\E $ be a sub{\sla } of\/ $\P (\Omega )$.  Then every tight character of $\E $ lies in the closure of the set
$\{\varphi _\omega :\omega \in \Omega \}\setminus \{0\}$, where each $\varphi _\omega $ is defined by
  $$
  \varphi _\omega (X)=\bool {\omega \in X}, \for X\in \E .
  $$

\Proof
  Follows by applying \ref {ApproxTightChar} to the identity representation $\id \colon \E \hookrightarrow \P (\Omega
)$.  \endProof

Proposition \ref {ApproxTightChar} speaks about the abundance of characters of $\E $ obtained by composing the
representation $\pi $ with \emph {principal} characters of $\P (\Omega )$, namely characters of the form
  $$
  X\mapsto [\omega \in X].
  $$

Given a representation $\pi $ of $\E $ in a Boolean algebra $\B $, it is therefore interesting to determine which
characters $\varphi $ of $\E $ factor as
  $$
  \varphi =\chi \circ \pi ,
  \equationmark CharFactor
  $$
  for some character $\chi $ of $\B $, preserving meets and
  joins\fn {In the theory of Boolean algebras, in fact also in the theory of lattices, characters are usually assumed to
preserve meets and joins, meaning that $\varphi (x\wedge y)=\varphi (x)\wedge \varphi (y)$, and $\varphi (x\vee
y)=\varphi (x)\vee \varphi (y)$, for all $x$ and $y$.  Virtually all characters of Boolean algebras in this work will be
supposed to preserve meets and joins but, since we are simultaneously dealing with {\sla }s and Boolean algebras, we
will try to be explicit every time these properties are required of a character.  Note that a character of a boolean
algebra preserves joins and meets if and only if it is an ultracharacter.
  }.
  An obvious necessary condition is that if $x, y_1,\ldots ,y_n$ are elements of $\E $ with
  $\pi (x)=\bigvee _{i=1}^n\pi (y_i)$,
  one must have that
  $\varphi (x)=\bigvee _{i=1}^n\varphi (y_i)$.

There are some useful equivalent forms of the above condition which will be important in the sequel.  First recall that
the category of Boolean algebras is equivalent to the category of Boolean rings.  (Recall that a Boolean ring is a ring,
necessarily commutative, in which every element is idempotent.)  If $R$ is a Boolean ring, then the Boolean algebra
structure on $R$ is given by $x\wedge y=xy$, $x\vee y=x+y-xy$ and $x\setminus y=x-y$ where $x\setminus y$ is the
relative complement for $y\leq x$.  Conversely, if $\B $ is a Boolean algebra, the ring structure takes $\wedge $ as the
multiplication and defines addition by $x+y=(x\setminus y)\vee (y\setminus x)$.  Boolean algebra characters on $\B $
correspond to surjective ring homomorphisms to the two-element field $\Bbb F_2$.

\state Proposition \label PropRelTight
  Let $\pi $ be a representation of the {\sla } $\E $ in a Boolean algebra $\B $, and let $\varphi $ be a character on
$\E $.  Then the following are equivalent:
  \iItemize
  \iItem for every $x, y_1,\ldots ,y_n\in \E $, one has that
  $$
  \pi (x)=\bigvee _{i=1}^n\pi (y_i) \Imply \varphi (x)=\bigvee _{i=1}^n\varphi (y_i),
  $$
  \iItem for every $x, y_1,\ldots ,y_n\in \E $, one has that
  $$
  \pi (x)\leq \bigvee _{i=1}^n\pi (y_i) \Imply \varphi (x)\leq \bigvee _{i=1}^n\varphi (y_i).
  $$
  \iItem for every $x, y_1,\ldots ,y_n\in \E $ with $y_i\leq x$, for $i=1,\ldots , n$, one has that
 $$
 \prod _{i=1}^n(\pi (x)-\pi (y_i))=0\Imply \prod _{i=1}^n(\varphi (x)-\varphi (y_i))=0
 $$
 for the Boolean ring structures on $\B $ and $\{0,1\}$.

\Proof
  (i)$\Imply $(ii): Assuming that $\pi (x)\leq \bigvee _{i=1}^n\pi (y_i)$, we have that
  $$
  \pi (x)=
  \pi (x)\wedge \big (\bigvee _{i=1}^n\pi (y_i)\big ) =
  \bigvee _{i=1}^n\pi (x)\wedge \pi (y_i) =
  \bigvee _{i=1}^n\pi (x y_i),
  $$
  so we may use (i) to deduce that
  $$
  \varphi (x) =
  \bigvee _{i=1}^n\varphi (x y_i) =
  \bigvee _{i=1}^n\varphi (x)\wedge \varphi (y_i) \leq
  \bigvee _{i=1}^n\varphi (y_i),
  $$
  proving (ii).

\medskip \noindent (ii)$\Imply $(i): \enspace Notice that if $\pi (x)=\bigvee _{i=1}^n\pi (y_i)$, then for each $i$, one
has that $\pi (y_i)\leq \pi (x)$, so (ii) implies that $\varphi (y_i)\leq \varphi (x)$, whence
  $$
  \bigvee _{i=1}^n\varphi (y_i)\leq \varphi (x).
  $$
  Since the reverse inequality also follows from (ii), the proof is concluded.

\medskip \noindent (ii)$\Imply $(iii): \enspace Note that $$ 0=\prod _{i=1}^n(\pi (x)-\pi (y_i))= \pi (x)\setminus
\bigvee _{i=1}^n\pi (y_i) $$ and so $\pi (x)\leq \bigvee _{i=1}^n\pi (y_i)$.  Therefore, by (ii), we have that $\varphi
(x)\leq \bigvee _{i=1}^n\varphi (y_i)$ and so $$ \prod _{i=1}^n(\varphi (x)-\varphi (y_i))= \varphi (x)\setminus \bigvee
_{i=1}^n\varphi (y_i)=0 $$ establishing (iii).

\medskip \noindent (iii)$\Imply $(ii): \enspace Notice that $\pi (x)\leq \bigvee _{i=1}^n\pi (y_i)$ implies that $\pi
(x)\leq \bigvee _{i=1}^n\pi (y_ix)$.
 Therefore,
 $$ \prod _{i=1}^n(\pi (x)-\pi (y_ix))= \pi (x)\setminus \bigvee _{i=1}^n\pi (y_ix)=0.  $$ So (iii) implies $$ 0=\prod
_{i=1}^n(\varphi (x)-\varphi (y_ix))= \varphi (x)\setminus \bigvee _{i=1}^n\varphi (y_ix) $$ and hence $$ \varphi
(x)\leq \bigvee _{i=1}^n\varphi (y_ix)\leq \bigvee _{i=1}^n\varphi (y_i) $$ as required.
  \endProof

The frequency with which we will use this condition largely justifies giving it a name:

\definition \label IntroRelTight
  Let $\pi $ be a representation of the {\sla } $\E $ in a Boolean algebra $\B $, and let $\varphi $ be a character on
$\E $.  We will say that $\varphi $ is \emph {tight relative to $\pi $}, or simply \emph {$\pi $-tight}, provided it
satisfies the equivalent conditions of \ref {PropRelTight}.  The set of all $\pi $-tight characters on $\E $ will be
written as $\hat \E _\pi $.
  \# $\hat \E _\pi $; Set of all $\pi $-tight characters on the {\sla } $\E $;

Observing that a character is $\pi $-tight if and only if it satisfies a set of equations, namely those in \ref
{PropRelTight.i}, it is clear that the set of all $\pi $-tight characters is closed in $\hat \E $.  For future reference
we record this fact below.

\state Proposition \label RelTightClosed
  Let $\pi $ be a representation of the {\sla } $\E $ in a Boolean algebra $\B $.  Then $\hat \E _\pi $ is a closed
subset of $\hat \E $.

As the reader might have already suspected, the necessary condition for the question raised in \ref {CharFactor} to have
a positive answer is also sufficient, as we shall prove next.

It will be convenient to use some rudiments of Stone duality between Boolean rings (algebras) and locally compact
Hausdorff spaces with a basis of compact open sets (generalized Stone spaces).  Most of this is well known in the unital
case and is folklore in the non-unital case but one must take care with morphisms.  If $X$ is a generalized Stone space,
then the ring $C_c(X,\Bbb F_2)$ of continuous functions with compact support from $X$ to $\Bbb F_2$ (with the discrete
topology) is a Boolean ring (note that every Boolean ring has characteristic $2$ and so is an $\Bbb F_2$-algebra).
Conversely, if $R$ is a Boolean ring, then ${\rm Spec}(R)=\hom (R,\Bbb F_2)$ is a generalized Stone space with the
topology of pointwise convergence.  These two constructions are inverse to each other up to isomorphism.  If $\pi \colon
R\to R'$ is a surjective homomorphism of Boolean rings, then it is obvious that ${\rm Spec}(R')$ embeds as a closed
subspace of ${\rm Spec}(R)$.  Suppose that $R'\leq R$ is a subring.  To show that restriction induces a surjective
continuous map ${\rm Spec}(R)\to {\rm Spec}(R')$, we need that every character of $R'$ extends to $R$.  This is well
known in the unital case and here is the proof in the non-unital case.

\state Proposition \label ExtendBooleanUltra Let $\B $ be a Boolean algebra and $\B '$ a non-zero subBoolean algebra.
Then every ultracharacter of $\B '$ extends to $\B $.

\Proof We prove this in the language of ultrafilters.  Let $\xi $ be an ultrafilter on $\B '$.  Then $\eta =\{x\in \B :
\exists y\in \xi , x\geq y\}$ is a filter on $\B $ with $\eta \cap \B '=\xi $.  It follows by Zorn's lemma that there is
an ultrafilter $\zeta $ on $\B $ containing $\eta $. Then $\zeta \cap \B '$ is a filter on $\B '$ containing $\xi $.
Since $\xi $ is an ultrafilter, we must have $\zeta \cap B'=\xi $.  This completes the proof.  \endProof

We remark that if $\E $ is a {\sla }, the semigroup algebra $\Bbb F_2\E $ is a Boolean ring and each character of $\E $
extends uniquely to a ring homomorphism $\Bbb F_2\E \to \Bbb F_2$.  Thus ${\rm Spec}(\Bbb F_2E)$ can be identified with
$\hat \E $.  We now show that if $\pi \colon \E \to \B $ is a representation into a Boolean algebra, then the $\pi
$-tight characters are precisely those that factor through $\pi $.

\state Theorem \label FactorlemmaTwo
 Let $\pi $ be a representation of the {\sla } $\E $ into the Boolean algebra $\B $ and let $\varphi $ be a character on
$\E $.  Then there exists a character $\chi $ of\/ $\B $, preserving meets and joins, such that $\varphi =\chi \circ \pi
$ if and only if $\varphi $ is $\pi $-tight.

\Proof We work with Boolean rings.  Consider the extension $\pi \colon \Bbb F_2\E \to \B $.  By \ref
{ExtendBooleanUltra}, we may replace $\B $ by $\pi (\Bbb F_2\E )$ and so we assume without loss of generality that $\pi
$ onto.  Thus we want to show that $\varphi \colon \Bbb F_2\E \to \Bbb F_2$ factors through $\pi $ if and only if the
corresponding character is $\pi $-tight.  The surjective homomorphism $\pi $ embeds ${\rm Spec}(\B )$ into $\hat \E $
and so if $X$ is the set of characters $\varphi $ with a factorization as $\chi \circ \pi $, then we can identify $\Bbb
F_2\E $ with $C_c(\hat \E ,\Bbb F_2)$ and $\B $ with $C_c(X,\Bbb F_2)$ and, moreover, we can identify $\pi $ with the
restriction map $f\mapsto f|_X$.  The kernel of the restriction map was determined in \cite [Proposition
5.2]{SteinbergPrimitive} in a more general setting (one must take $\Bbbk =\Bbb F_2$ and $S=\E $ in that theorem).
Namely, $\ker \pi $ is the ideal generated by all products $\prod _{i=1}^n(e-e_i)$ such that $e,e_1,\ldots , e_n\in E$,
$e_i\leq e$ and $$ D(e)\cap D(e_1)^c\cap \cdots \cap D(e_n)^c\cap X=\emptyset .  \equationmark IntersectEq $$ We claim
that \ref {IntersectEq} holds if and only if $\prod _{i=1}^n(\pi (e)-\pi (e_i))=0$.  Indeed, if we have $\prod
_{i=1}^n(\pi (e)-\pi (e_i))=0$, then trivially, for any character $\chi $ of $\B $, $\prod _{i=1}^n(\chi (\pi (e))-\chi
(\pi (e_i)))=0$ and hence either $\chi (\pi (e))=0$ or $\chi (\pi (e))=1=\chi (\pi (e_i))$ for some $i$. Thus $\chi
\circ \pi \notin D(e)\cap D(e_1)^c\cap \cdots \cap D(e_n)^c$ and so \ref {IntersectEq} holds.  Conversely, if $x=\prod
_{i=1}^n(\pi (e)-\pi (e_i))\neq 0$, then there is a character $\chi $ of $\B $ with $\chi (x)=1$ (choose an ultrafilter
containing the principal filter generated by $x$).  Then $1=\chi (x) = \prod _{i=1}^n(\chi (\pi (e))-\chi (\pi (e_i)))$
and so $\chi \circ \pi \in D(e)\cap D(e_1)^c\cap \cdots \cap D(e_n)^c\cap X$ and hence \ref {IntersectEq} fails.

It now follows from \ref {PropRelTight} and \cite [Proposition 5.2]{SteinbergPrimitive} that $\varphi $ factors through
$\pi $ if and only if it is $\pi $-tight.  \endProof

Notice that the above proof shows that if $\pi \colon \E \to \B $ is a homomorphism such that $\pi (\E )$ generates $\B
$ as a Boolean algebra, then ${\rm Spec}(\B )$ is homeomorphic to $\hat \E _\pi $.  The following special case will be
used repeatedly.

\state Theorem \label Factorlemma
  Let $\pi $ be a representation of the {\sla } $\E $ on a set $\Omega $, and let $\varphi $ be a $\pi $-tight character
on $\E $.  Then there exists a character $\chi $ of\/ $\P (\Omega )$, preserving meets and joins, such that $\varphi
=\chi \circ \pi $.

We shall now be interested in representing {\sla }s in associative algebras.

\definition
  Let $\E $ be a {\sla } and let $A$ be an associative algebra over a field $\K $.  A mapping $\pi \colon \E \to A$ is
said to be a representation of $\E $ in $A$, if $\pi (0)=0$, and $\pi (xy)=\pi (x)\pi (y)$, for all $x$ and $y$ in $A$.

Given a {\sla } $\E $, and a representation $\pi $ of $\E $ in an algebra $A$, one may often assume that $A$ is abelian
by replacing $A$ with the subalgebra of $A$ generated by the range of $\pi $.
  If $A$ is indeed abelian, we may view $\pi $ as a representation in a canonically defined Boolean algebra as follows:
  $$
  \BA :=\{e\in A: e^2=e\}.
  $$
  Under the operations
  $$
  e\wedge f= ef, \and e\vee f = e+f-ef, \for e,f\in \BA
  $$
  it is easy to see that $\BA $ is a Boolean algebra and clearly $\pi $ takes values in $\BA $.  All concepts relating
to Boolean algebra representations, such as $\pi $-tightness, therefore immediately apply to representations in
associative algebras by considering the associated representation in $\BA $.  If $A$ is a commutative $\K $-algebra,
then its spectrum $\hat A$ is the space of non-zero $\K $-algebra homomorphisms $f\colon A\to \K $ equipped with the
topology of pointwise convergence (where $\K $ is viewed as a discrete space).

\state Proposition \label SpecAlg
  Let $\pi $ be a representation of a {\sla } $\E $ in an associative $\K $-algebra $A$, such that $A$ is generated by
the range of\/ $\pi $.  Then:
  \iItemize
  \iItem The spectrum of $A$, which we will denote by $\hat A$, is homeomorphic to $\hat \E _\pi $.
  \iItem $A$ is isomorphic to the algebra $C_c(\hat \E _\pi ,\K )$, consisting of all locally constant, compactly
supported, $\K $-valued functions on $\hat \E _\pi $.

\Proof
  Since $A$ is an abelian algebra generated by idempotents, it follows from \cite [Corollaire 1]{Keimel} that $A$ is
naturally isomorphic to $C_c({\rm Spec}(\BA ),\K )$.  Since $\pi \colon \E \to \BA $ extends to a surjective
homomorphism $\pi \colon \Bbb F_2\E \to \BA $, \ref {FactorlemmaTwo} and its proof shows that we can identify ${\rm
Spec}(\BA )$ with $\hat \E _\pi $.  It remains to show that $\hat A$ is homeomorphic to ${\rm Spec}(\BA )$.  The
restriction map $\hat A\to {\rm Spec}(\BA )$ is clearly continuous and injective (the latter since $\pi (\E )$ generates
$A$).  It is onto because we can identify $\BA $ with the characteristic functions of compact open subsets of ${\rm
Spec}(\BA )$ and then if $\xi \in {\rm Spec}(\BA )$, the corresponding character is evaluation at $\xi $.  But this
extends to $C_c({\rm Spec}(\BA ),\K )$ as evaluation at $\xi $.  We show that the restriction map is open.  A basic
neighborhood $U$ in $\hat A$ specifies the value of a character at finitely many elements of $\BA $ (since each element
can be expressed as a finite linear combination of idempotents and hence the value of a character at any element is
determined by its value at finitely many idempotents).  Since a character can only take on values $0$ and $1$ on an
idempotent (because $\K $ is a field) it follows that the image of $U$ is the open subset of ${\rm Spec}(\BA )$ of
characters taking the same values on those specified idempotents.
  \endProof

We now begin to apply some of this machinery in the strongly finitely aligned case.  This will be used in our sequel
paper to show that the tight C*-algebra of the inverse hull of the semigroup associated to a finitely aligned higher
rank graph is the corresponding higher rank graph C*-algebra.  Again, we shall assume our strongly finitely aligned
semigroup has right local units.

\state Lemma \label alignedcover
  Let $S$ be strongly finitely aligned and let $s,t\in S$.  Suppose that $sS\cap tS$ has basis $B=\{w_1,\ldots , w_n\}$.
Write $w_i=sx_i=ty_i$ with $x_i,y_i\in S$.  \iItemize
  \iItem $\theta _s\inv \theta _t\theta _{y_i}\theta _{y_i}\inv = \theta _{x_i}f_{w_i}\theta _{y_i\inv }$.
  \iItem The set
  $\{\theta _{y_i}\theta _{y_i}\inv \theta _t\inv \theta _t\mid 1\leq i\leq n\}$ is a cover of $\theta _t\inv \theta
_s\theta _s\inv \theta _t$.

\Proof (i) If $u$ is in the domain of the left hand side, then $u=y_iz$ with $z\in S$, $tu\neq 0$ and $\theta _s\inv
\theta _t\theta _{y_i}\theta _{y_i}\inv (u) = \theta _s\inv (ty_iz) = \theta _s\inv (sx_iz) = x_iz$.  On the other hand,
$w_iz=ty_iz=tu\neq 0$ and so $\theta _{x_i}f_{w_i}\theta _{y_i}\inv (u)=\theta _{x_i}(z)=x_iz$.  Similarly, if $u$ is in
the domain of $\theta _{x_i}f_{w_i}\theta _{y_i\inv }$, then $u=y_iz$ with $w_iz\neq 0$ and $\theta _{x_i}f_{w_i}\theta
_{y_i\inv }(u) = x_iz$.  But $tu=ty_iz=w_iz\neq 0$ implies that $\theta _s\inv \theta _t\theta _{y_i}\theta _{y_i}\inv
(u)=\theta _s\inv (tu) = \theta _s\inv (sx_iz)=x_iz$ as $tu=ty_iz=sx_iz$.

(ii) First note that if $\theta _{y_i}\theta _{y_i}\inv \theta _t\inv \theta _t(x)=x$, then $x=y_iz$ with
$tx=ty_iz=sx_iz\neq 0$.  Then $\theta _t\inv \theta _s\theta _s\inv \theta _t(x) =\theta _t\inv \theta _s\theta _s\inv
(tx)=\theta _t\inv \theta _s\theta _s\inv (sx_iz)=\theta _t\inv (tx) =x$ and so we have $\theta _{y_i}\theta _{y_i}\inv
\theta _t\inv \theta _t\leq \theta _t\inv \theta _s\theta _s\inv \theta _t$.

  Suppose now that $0\neq f\leq \theta _t\inv \theta _s\theta _s\inv \theta _t$ and that $f(x)\neq 0$.  Then $\theta
_t\inv \theta _s\theta _s\inv \theta _t(x)\neq 0$ and so $tx\neq 0$ and $tx=sy$ for some $y\in S$.  Therefore, $tx\in
sS\cap tS$ and so $tx=w_iz=ty_iz$ for some $z\in S$ and $1\leq i\leq n$.  Then $x=y_iz$ by $0$-left cancellativity and
so $\theta _{y_i}\theta _{y_i}\inv \theta _t\inv \theta _t(x)=x$. Thus $f\theta _{y_i}\theta _{y_i}\inv \theta _t\inv
\theta _t\neq 0$.
  \endProof

Let us say that a representation $\pi \colon T\to A$ of $T$ in an associative algebra $A$ is \emph {cover-to-join} if
whenever $\{f_1,\ldots , f_n\}$ is a cover of an idempotent $e$ of $E(T)$ with $f_i\leq e$ for all $i=1,\ldots , n$,
then $\pi (e)=\bigvee _{i=1}^n\pi (f_i)$ (where the join is taken in the commutative algebra generated by $\pi (E(T))$),
i.e, $\prod _{i=1}^n(\pi (e)-\pi (f_i))=0$.

  \state Corollary \label formulaalignedcase
  Let $S$ be a strongly finitely aligned $0$-left cancellative semigroup and let $s,t\in S$.  Suppose that $sS\cap tS$
has basis $B=\{w_1,\ldots , w_n\}S$.  Write $w_i=sx_i=ty_i$ with $x_i,y_i\in S$.
 If $\pi \colon \hull \to A$ is a cover-to-join $\ast $-representation into an associative $\ast $-algebra $A$ and if we
put $\pi _s=\pi (\theta _s)$, then
 $$\pi _s^*\pi _t = \sum _{i=1}^m \pi _{x_i}\pi (f_{w_i})\pi _{y_i}^*.$$
 holds.  In particular, if $S$ is categorical at zero and right reductive, then the formula $$\pi _s^*\pi _t = \sum
_{i=1}^m \pi _{x_i}\pi _{y_i}^*$$ holds.

 \Proof Since $w_iS\cap w_jS=\{0\}$ for $i\neq j$, we have by \ref {alignedcover} that
  $$
  \pi _t^*\pi _s\pi _s^*\pi _t = \sum _{i=1}^n \pi _t^*\pi _t\pi _{y_i}\pi _{y_i}^*
  $$
  and so
  $$
  \pi _s^*\pi _t = \sum _{i=1}^n\pi _s^*\pi _t\pi _{y_i}\pi _{y_i}^* = \sum _{i=1}^n \pi _{x_i}\pi (f_{w_i})\pi
_{y_i}^*,
  $$
  where the last equality follows from \ref {alignedcover}.  The final statement follows because $w_i=sx_i$ implies
$w_i^+=x_i^+$ and hence $\theta _{x_i}f_{w_i} = \theta _{x_i}\theta _{w_i^+} = \theta _{x_i}$.
    \endProof

As a corollary, we deduce that the image of $\hull $ of a strongly finitely aligned semigroup under a cover-to-join
$\ast $-representation $\pi $ is spanned by elements of the form $\pi _s\pi (f_{\Lambda })\pi _t\inv $ where $\Lambda $
is a finite subset of $S$, which is similar to what happens in the lcm case.

\state Theorem \label describetight Let $S$ be a strongly finitely aligned $0$-left cancellative semigroup.  Let $\pi
\colon \hull \to A$ be a cover-to-join $\ast $-representation to an associative $\ast $-algebra $A$ such that $\pi
(\hull )$ spans $A$ and write $\pi _s$ for $\pi (\theta _s)$ for $s\in S$.  Then $A$ is spanned by elements of the form
$\pi _s\pi (f_{\Lambda })\pi _t^*$ with $\Lambda \subseteq S$ finite and $s,t\in \Lambda $.  If $S$ is right reducitive
and categorical at $0$, then $A$ is spanned by elements of the form $\pi _s\pi _t^*$ with $s,t\in S$.

\Proof Let $T$ be the set of elements of the form $\pi _s\pi (f_{\Lambda })\pi _t^*$ with $\Lambda \subseteq S$ finite
and $s,t\in \Lambda $.
 Trivially, $T\subseteq \pi (\hull )$.  Also, if $e\in E(S)$ with $se=s$, then $\pi _s =\pi _s\pi (f_s)\pi _{e}^*\in T$.
So it suffices to prove that $T$ is an inverse semigroup.
 Since $(\pi _s\pi (f_{\Lambda })\pi _t^*)^* = \pi _t\pi (f_{\Lambda })\pi _s^*$, it suffices to show that $T$ is a
semigroup.  So let $s_i,t_i\in \Lambda _i\subseteq S$, for $i=1,2$, with $\Lambda _i$ finite.  Let $\{w_1,\ldots ,
w_n\}$ be a basis for $t_1S\cap s_2S$ and write $w_i = t_1x_i=s_2y_i$.

 First observe that
  $$
  \pi _{s_1}\pi (f_{\Lambda _1})\pi _{x_i}\pi (f_{w_i})\pi _{y_i}^*\pi (f_{\Lambda _2})\pi _{t_2}^* = \pi _{s_1x_i}\pi
(f_{\Lambda _1x_i\cup \{w_i\}\cup \Lambda _2y_i})\pi _{t_2y_i}^*
  $$
  and $s_1x_i,t_2y_i\in \Lambda _1x_i\cup \{w_i\}\cup \Lambda _2y_i$.
   Then applying this and \ref {formulaalignedcase}, yields
   $$
  \pi _{s_1}\pi (f_{\Lambda _1})\pi _{t_1}^*\cdot \pi _{s_2}\pi (f_{\Lambda _2})\pi _{t_2}^* = \sum _{i=1}^n\pi
_{s_1}\pi (f_{\Lambda _1})\pi _{x_i}\pi (f_{w_i})\pi _{y_i}^*\pi (f_{\Lambda _2})\pi _{t_2}^* \quebra = \sum _{i=1}^n\pi
_{s_1x_i}\pi (f_{\Lambda _1x_i\cup \{w_i\}\cup \Lambda _2y_i})\pi _{t_2y_i}^*
  $$
  as required.

   The statement in the categorical at zero case is proved similarly, using the final formula in \ref
{formulaalignedcase}.  \endProof

Actually, if one follows the above proofs carefully one can weaken the assumption that $\pi $ is cover-to-join to just
ask that if a constructible set $X$ can be written as a finite union of constructible sets $X=\bigcup _{i=1}^nX_i$, then
$\pi (1_X)=\bigvee _{i=1}^n\pi (1_{X_i})$.  That is, Theorem~\ref {describetight} is true for $\theta $-tight $\ast
$-representations $\pi $ where $\theta \colon E(\hull )\to \P (S\setminus \{0\})$ is the restriction of the regular
representation of $\hull $.

\section The spectrum of the {\sla } of constructible sets

\label StringSection \fix Throughout this section we will fix a $0$-left cancellative semigroup $S$ admitting least
common multiples.

As mentioned in the preamble to Part (2), we are interested in the study of the {\sla } $\ehull $ of constructible sets,
with a special emphasis on its spectrum.  In order to exhibit examples of characters on $\ehull $ we shall make use of
the idea behind \ref {CharacterFromRep} so, considering the representation
  $$
  \varepsilon \colon \ehull \to \P (S^\star ),
  $$
  introduced in \ref {IntroEpsilon},
  and given $\sigma $ in the essential subset for $\varepsilon $, we may define
  $$
  \varphi ^\varepsilon _\sigma \colon X\in \ehull \mapsto [\sigma \in \varepsilon (X)] \ \in \ \{0, 1\}.
  $$
  \# $\varphi _\sigma $; Character associated with the string $\sigma $;
  which is a character on $\ehull $, as discussed near \ref {BigSubset}.

Observe that, for trivial reasons, the essential subset for $\varepsilon $ coincides with the essential subset $S^\star
_\ess $ for $\theta ^\star $, which we have seen in \ref {EssThetaStar} to consist of all strings except for the
singletons $\{s\}$, where $s$ is a
  degenerate\fn {Recall from \ref {DefDegenElt} that $s$ is degenerate if $s$ is irreducible and $Ss=\{0\}$.}
  element of $S$.

\definition \label IntroPhiSigmaNonDeg
  \iItemize
  \iItem A string $\sigma $ will be called \emph {degenerate} if $\sigma =\{s\}$, where $s$ is a degenerate element.
The set of all non-degenerate strings is therefore $S^\star _\ess $.
  \iItem For every non-degenerate string $\sigma $, we shall denote by $\varphi _\sigma $ the character of $\ehull $
given by
  $$
  \varphi _\sigma (X)= [\sigma \in \varepsilon (X)], \for X\in \ehull .
  $$

Employing the terminlogy introduced in \ref {BigSubset} we have that
  $$
  \spec _{S^\star } = \{\varphi _\sigma :\sigma \in S^\star _\ess \},
  $$
  which is thus a subset of $\spec $, allowing for a first glimpse of our main object of study.

\state Proposition One has that $\spec _{S^\star }$ is dense in the tight spectrum of $\ehull $.

\Proof We will prove this as an application of \ref {ApproxTightChar}, and hence our task consists in showing that
$\varepsilon (X)$ is nonempty for every nonempty $X$ in $\ehull $.  But this is immediate from Theorem \ref
{BiggerCovariance}, which implies that $r\in X$ if and only if $\delta _r\in \varepsilon (X)$.  \endProof

We should note that, since $\ehull $ is a subset of $\P (S')$, one may use \ref {TautoCharsDense} to produce another set
of characters which is also dense in the $\varspec {_\tight }$.  However, characters arising from the identity
representation of $\ehull $ on $S'$ have a very small chance of being tight, so we shall not be able to benefit much
from \ref {TautoCharsDense} in the present context.  As we shall see later (see \ref {BigNewTightResult}), characters
arising from strings are much more likely to be ultracharacters, and hence tight.

Suppose we are given $\varphi _\sigma $ and we want to recover $\sigma $ from $\varphi _\sigma $.  In the special case
in which $S$ has right local units, we have that
  $$
  \varphi _\sigma (E^\theta _s) = 1 \iff
  \sigma \in \varepsilon (E^\theta _s) = E^\star _s \explain {FRFstarR.iii} \iff
  s\in \sigma ,
  \equationmark ClueStringRLU
  $$
  so $\sigma $ is recovered as the set $\{s\in S: \varphi _\sigma (E^\theta _s)=1\}$.  Without assuming right local
units, the last part of \ref {ClueStringRLU} cannot be trusted, but it may be replaced with
  $$
  \phantom {\varphi _\sigma (E^\theta _s) = 1 \iff
  \sigma \in \varepsilon (E^\theta _s) = E^\star _s}
  \cdots \explain {FRFstarR.ii} \iff
  \sigma \cap E^\theta _s\neq \emptyset ,
  \equationmark ClueString
  $$
 so we at least know which $E^\theta _s$ have a nonempty intersection with $\sigma $.

\state Proposition \label PropsInterior
  Given any string $\sigma $, one has that the set
  $$
  \{s\in S: \varphi _\sigma (E^\theta _s)=1\},
  $$
  coincides with $\interior \sigma $, namely the interior of $\sigma $, as defined in \ref {DefOpenString}.

\Proof
  Observe that since $0\notin \sigma $, one has for every $s$ in $S$, that
  $$
  \sigma \cap sS = \sigma \cap (sS\setminus \{0\}) = \sigma \cap E^\theta _s.
  $$
  By \ref {ClueString} one then has that $s$ lies in the set displayed in the statement if and only if $\sigma \cap
sS\neq \emptyset $, so the statement follows.  \endProof

\medskip Given any character $\varphi $ of $\ehull $, regardless of whether or not it is of the form $\varphi _\sigma $
as above, we may still consider the set
  $$
  \sigma _\varphi :=\{s\in S: \varphi (E^\theta _s)=1\},
  \equationmark defineSigFi
  $$
  \# $\sigma _\varphi $; String associated with the character $\varphi $ (when nonempty);
  so that, when $\varphi =\varphi _\sigma $, we get $\sigma _{\varphi } =\interior \sigma $, by \ref {PropsInterior}.

\state Proposition \label StringFromChar
  If $\varphi $ is any character of\/ $\ehull $, and $\sigma _\varphi $ is nonempty, then $\sigma _\varphi $ is a string
which is moreover closed under least common multiples.

\Proof
  Since $E^\theta _0=\emptyset $, we have that $\varphi (E^\theta _0)=0$, whence $0\notin \sigma _\varphi $.  If $t\in
\sigma _\varphi $, and if $s$ divides $t$, then
  $tS\subseteq sS$, whence also $E^\theta _t\subseteq E^\theta _s$, and then $1=\varphi (E^\theta _t)\leq \varphi
(E^\theta _s)$, so $s\in \sigma _\varphi $.

To prove that $\sigma _\varphi $
  is closed under least common multiples, let $s,t\in \sigma $, and let $r$ be a least common multiple of $s$ and $t$.
Then $E^\theta _r= E^\theta _s\cap E^\theta _t$, whence
  $$
  \varphi _\sigma (E^\theta _r)= \varphi _\sigma (E^\theta _s)\varphi _\sigma (E^\theta _t) = 1,
  $$
  so $r\in \sigma _\varphi $.  This also implies that $\sigma _\varphi $ satisfies \ref {StringDirected}.
  \endProof

Based on \ref {IntroPhiSigmaNonDeg.ii} we may define a map from the set of all non-degenerate strings to
  $\spec $,
  \# $\textbackslash spec$; The spectrum of $\ehull $;
  the spectrum of $\ehull $, by
  $$
  \Phi : \sigma \in S^\star _\ess \mapsto \varphi _\sigma \in \spec ,
  \equationmark DefinePhi
  $$
  \# $\Phi $; Map sending non-degenerate strings to their associated characters;
  but if we want the dual correspondence suggested by \ref {defineSigFi}, namely
  $$
  \varphi \mapsto \sigma _\varphi ,
  \equationmark ProposeSigma
  $$
  to give a well defined map from $\spec $ to $S^\star $, we need to worry about its domain because we have not checked
that $\sigma _\varphi $ is always nonempty, and hence $\sigma _\varphi $ may fail to be a string.  The appropriate
domain is evidently given by the set of all characters $\varphi $ such that $\sigma _\varphi $ is nonempty but, before
we formalize this map, it is interesting to introduce a relevant sub{\sla } of $\ehull $.

\state Proposition \label IntroEOne
  The subset of\/ $\ehull $ given by\fn {This should be contrasted with \ref {FormOfHull}, where the general form of an
element of $\ehull $ is shown to be $uF^\theta _\Lambda $, where $u$ is in $\tilde S$, rather than $S$.}
  $$
  \eone = \{sF^\theta _\Lambda ,\
  \Lambda \subseteq S \hbox { is finite, and } s\in \Lambda \},
  $$
  \# $\eone $; A sub{\sla } of $\ehull $;
  is an ideal of\/ $\ehull $.  Moreover, for every $X$ in $\ehull $, one has that $X$ lies in $\eone $ if and only if
$X\subseteq E^\theta _s$, for some $s$ in $S$.

\Proof
  Let us first prove the last sentence of the statement.  The ``only if\kern 1pt" part is evident since
  $$
  sF^\theta _\Lambda \subseteq E^\theta _s,
  $$
  so let us focus on the ``if\kern 1pt" part.  We thus suppose that $X\in \ehull $ is such that $X\subseteq E^\theta
_s$, for some $s\in S$.  Observing that
  $$
  X\subseteq E^\theta _s=sF^\theta _s,
  $$
  we see that $X=sxF^\theta _\Delta $, where $x$ and $\Delta $ are as in \ref {FormOfSubSets}, whence
  $X\in \eone $.  The first part of the statement, namely that $\eone $ is an ideal, now follows easily.
  \endProof

Observe that if $S$ admits right local units, then $\eone =\ehull $, thanks to \ref {LocUnitsFormOfSlatice}.

By \ref {PandI} we may then view $\seone $ as an open subset of $\spec $.  The next result is intended to distinguish
the elements of $\seone $ within $\spec $.

\state Proposition \label DistinguishSeone
  Let $\varphi $ be a character on $\ehull $.  Then the following are equivalent:
  \iItemize
  \iItem $\varphi \in \seone $,
  \iItem $\varphi (E^\theta _s)=1$, for some $s$ in $S$,
  \iItem $\sigma _\varphi $ is nonempty, and hence it is a string by \ref {StringFromChar}.

\Proof
  The equivalence between (i) and (ii) follows from \ref {SpecIdealNoVanish} and \ref {IntroEOne}, while (ii) and (iii)
are obviously equivalent.
  \endProof

If $S$ admits right local units, we have seen that $\eone =\ehull $, so $\sigma _\varphi $ is a string for every
character $\varphi \in \spec $.

The vast majority of non-degenerate strings $\sigma $ lead to a character $\varphi _\sigma $ belonging to $\seone $, but
there are exceptions.

\state Proposition \label IrrStrings
  If $\sigma $ is a non-degenerate string in $S^\star _\ess $ then $\varphi _\sigma $ does not belong to $\seone $ if
and only if $\sigma =\{s\}$, where $s$ is an irreducible element of $S$.

\Proof By \ref {DistinguishSeone} to say that $\varphi _\sigma $ is not in $\seone $, is to say that for all $t$ in $S$
we have that $\varphi _\sigma (E^\theta _t)=0$, or, equivalently that $\sigma \notin E^\star _t$, by \ref {FRFstarR.ii}.
The conclusion then follows from \ref {ZeroChars.i}.
  \endProof

By \ref {DistinguishSeone} we have
 that the largest set of characters on which the correspondence described in \ref {ProposeSigma} produces a bona fide
string is precisely $\seone $, so we may now formally introduce the map suggested by that correspondence.

\definition \label IntroSigma
  We shall let
  $$
  \Sigma :\seone \mapsto S^\star ,
  $$
  \# $\Sigma $; Map sending characters to their associated strings;
  be the map
  given by
  $$
  \Sigma (\varphi ) = \sigma _\varphi =\{s\in S: \varphi (E^\theta _s)=1\}, \for \varphi \in \seone .
  $$

For every string $\sigma $, excluding the exceptional ones discussed in \ref {IrrStrings}, we then have that
  $$
  \Phi (\sigma )=\varphi _\sigma \in \seone ,
  $$
  and
  $$
  \Sigma \big (\Phi (\sigma )\big )=\interior \sigma ,
  \equationmark SigmaPhiIsInterior
  $$
  by \ref {PropsInterior}.  The nicest situation is for open strings:

\state Proposition \label SigmaPhiOnOpen
  If $\sigma $ is an open string, then
  \iItemize
  \iItem $\sigma $ is non-degenerate,
  \iItem $\Phi (\sigma )\in \seone $, and
  \iItem $\Sigma \big (\Phi (\sigma )\big )=\sigma $.

\Proof
  If $s$ is an irreducible element in $S$, then the string $\{s\}$ is certainly not open, so an open string cannot be
any of the exceptional strings discussed in \ref {IrrStrings}, much less a degenerate string.  Therefore $\Phi (\sigma
)\in \seone $.  The third point follows
  from \ref {SigmaPhiIsInterior}
  and the fact that $\sigma =\interior \sigma $.  \endProof

Given that the composition $\Sigma \circ \Phi $ is so well behaved for open strings, we will now study the reverse
composition $\Phi \circ \Sigma $ on a set of characters related to open strings.

\definition \label DefineOpenChar
  A character $\varphi $ in $\spec $ will be called an \emph {open} character if $\sigma _\varphi $ is a (nonempty) open
string.

We remark that every open character belongs to $\seone $ by \ref {DistinguishSeone}, although not all characters in
$\seone $ are open.

By \ref {SigmaPhiOnOpen} it is clear that $\varphi _\sigma $ is an open character for every open string $\sigma $.

If $S$ admits right local units, we have seen that every string in $S^\star $ is open, and also that $\sigma _\varphi $
is a string for every character.  Therefore every character in $\spec $ is open.

The composition $\Phi \circ \Sigma $ is not as well behaved as the one discussed in \ref {SigmaPhiOnOpen}, but there is
at least some relationship between a character $\varphi $ and its image under $\Phi \circ \Sigma $, as we shall now see.

\state Proposition \label RioMaior
  Given any open character $\varphi $, one has that
  $$
  \varphi \leq \Phi \big (\Sigma (\varphi )\big ).
  $$

\Proof
  In an effort to decongest notation, throughout this proof we will write $\sigma $ for $\sigma _\varphi $, so the
inequality in the statement reads $\varphi \leq \varphi _{\sigma }$.  In order to prove it, it is clearly enough to
argue that there is no $X$ in $\ehull $ such that
  $$
  \varphi (X)=1, \and \varphi _\sigma (X)=0.
  \equationmark AbsurdHyp
  $$
  Arguing by contradiction, we assume that such an $X$ exists.  Observing that $\varphi $ is in $\seone $, we may choose
$s$ in $S$ such that $\varphi (E^\theta _s)=1$.  Setting $X'=X\cap E^\theta _s$, we have that
  $$\def \quad { }
  \matrix {
  \varphi (X')& = & \varphi (X)\hfil \varphi (E^\theta _s) & = & 1, &\ \ \hbox {and} \cr \pilar {14pt}
  \varphi _\sigma (X')& = & \varphi _\sigma (X)\varphi _\sigma (E^\theta _s) & = & 0,}
  $$
  which means that we may suppose without loss of generality that the originally chosen $X$ is a subset of $E^\theta
_s$.  Using \ref {IntroEOne} we have that $X\in \eone $, so we may write $X=rF^\theta _\Lambda $, where $\Lambda $ is a
finite subset of $S$, and $r\in \Lambda $.  Noticing that $X\subseteq E^\theta _r$, we have
  $$
  1=\varphi (X)\leq \varphi (E^\theta _r),
  $$
  so we see that $r\in \sigma =\sigma _\varphi $.  Since $\varphi $ is an open character, $\sigma $ is an open string,
so there exists some
   $t$ in $S$, such that $rt\in \sigma $.  Therefore $\sigma \cap E^\theta _r\neq \emptyset $, and we deduce from \ref
{FRFstarR.ii} that $\sigma $ is in $E^\star _r$.

Notice that to say that $\varphi _\sigma (X)=0$ is the same as saying that
  $$
  \sigma \notin \varepsilon (X) = \theta ^\star _r(F^\star _\Lambda ),
  $$
  which implies that
  $$
  r\inv *\sigma = {\theta ^\star _r}\inv (\sigma ) \notin F^\star _\Lambda ,
  $$
  or, equivalently, that
  $$
  r\inv *\sigma \not \subseteq F^\theta _\Lambda ,
  $$
  by \ref {FRFstarR.i}.  This said, we may pick $y\in r\inv *\sigma $, such that $ty=0$, for some $t\in \Lambda $.  In
particular $ry\in \sigma $, so
  $$
  \varphi (E^\theta _{ry})=1.
  \equationmark XiRX
  $$

We next claim that $X$ and $E^\theta _{ry}$ are disjoint.  In fact, should this not be the case, we could find some
  $$
  s\in X\cap E^\theta _{ry} = rF^\theta _\Lambda \cap E^\theta _{ry},
  $$
  which may therefore be written as
  $$
  s = rp = ryq,
  $$
  for suitable $p$ in $F^\theta _\Lambda $, and $q$ in $S$.  Therefore $p=yq$, by $0$-left cancellativity, and then
  $$
  0\neq tp=tyq=0,
  $$
  a contradiction.  This proves that $X\cap E^\theta _{ry}=\emptyset $, so
  $$
  0 = \varphi (X\cap E^\theta _{ry}) = \varphi (X) \varphi (E^\theta _{ry}) \={XiRX}
  \varphi (X),
  $$
  contradicting \ref {AbsurdHyp}, and thus concluding the proof.
  \endProof

We have thus arrived at an important result.

\state Theorem \label BigNewTightResult
  Let $S$ be a $0$-left cancellative semigroup admitting least common multiples.  Then, for every open, maximal string
$\sigma $ over $S$, one has that $\varphi _\sigma $ is an ultracharacter.

\Proof
  Let $\psi $ be a character such that $\varphi _\sigma \leq \psi $.  For every $s$ in $S$ we then have that
  $$
  s\in \sigma \explain {SigmaPhiOnOpen.iii}\Imply \varphi _\sigma (E^\theta _s)=1 \Imply \psi (E^\theta _s)=1 \Imply
s\in \sigma _\psi ,
  $$
  which means that $\sigma \subseteq \sigma _\psi $, and hence that $\sigma =\sigma _\psi $, by maximality.  It follows
that $\psi $ is an open character, so \ref {RioMaior} applies for $\psi $, and we deduce that
  $$
  \varphi _\sigma \leq \psi \explain {RioMaior}\leq \Phi \big (\Sigma (\psi )\big ) = \Phi (\sigma ) = \varphi _\sigma ,
  $$
  so $\varphi _\sigma =\psi $, proving that $\varphi _\sigma $ is maximal.
  \endProof

The previous result raises the question as to whether $\sigma _\varphi $ is a maximal string for every ultracharacter
$\varphi $, but this is not true in general.  Consider for example the unital semigroup
  $$
  S=\{1,a,aa,0\},
  $$
  in which $a^3=0$.
  The $\theta $-constructible subsets of $S$ are precisely

  \vbox { \bigskip \begingroup \offinterlineskip
  $$
  \vcenter {\halign {
  \vrule height 14pt depth 8pt\ \ \hfill #\hfill \ \ &\vrule \ \ \hfill #\hfill \ \ &\vrule \ \ \hfill #\hfill \ \
\vrule \cr \tabrule
  $E^\theta _1 = F^\theta _1 = \{1,a,aa\}$ & & \cr \tabrule
  $\hfill F^\theta _a = \{1,a\}$ & $E^\theta _a = aF^\theta _a = \{a,aa\}$ & \cr \tabrule
  $\hfill F^\theta _{aa} = \{1\}$ & $\hfill aF^\theta _{aa} = \{a\}$ & $E^\theta _{aa} = aaF^\theta _{aa} = \{aa\}$ \cr
\tabrule
  }}
  $$
  \vskip -5pt\centerline {\eightrm List of $\scriptstyle \theta $-constructible sets}
  \endgroup
  \bigskip }

  \noindent and there are three strings over $S$, namely

  \bigskip \begingroup \offinterlineskip
  $$
  \vcenter {\halign {
  \vrule height 14pt depth 8pt\ \ \hfill #\hfill \ \ &\vrule \ \ \hfill #\hfill \ \ &\vrule \ \ \hfill #\hfill \ \
\vrule \cr \tabrule
  $\delta _1=\{1\}$ & $\delta _a=\{1,a\}$ & $\delta _{aa}=\{1,a,aa\}$ \cr \tabrule
  }}
  $$
  \endgroup \bigskip

  Since the correspondence $s\mapsto \delta _s$ is a bijection from $S'$ to $S^\star $, we see that $\theta ^\star $ is
isomorphic to $\theta $, and in particular the $\theta ^\star $-constructible subsets of $S^\star $, listed below,
mirror the $\theta $-constructible ones.

  \vbox {
  \bigskip \begingroup \offinterlineskip
  $$
  \vcenter {\halign {
  \vrule height 14pt depth 8pt\ \ \hfill #\hfill \ \ &\vrule \ \ \hfill #\hfill \ \ &\vrule \ \ \hfill #\hfill \ \
\vrule \cr \tabrule
  $E^\star _1 = F^\star _1 = \{\delta _1,\delta _a,\delta _{aa}\}$ & & \cr \tabrule
  $\hfill F^\star _a = \{\delta _1,\delta _a\}$ & $E^\star _a = aF^\star _a = \{\delta _a,\delta _{aa}\}$ & \cr \tabrule
  $\hfill F^\theta _{aa} = \{\delta _1\}$ & $\hfill aF^\star _{aa} = \{\delta _a\}$ & $E^\theta _{aa} = aaF^\star _{aa}
= \{\delta _{aa}\}$ \cr \tabrule
  }}
  $$
  \endgroup
  \vskip -5pt\centerline {\eightrm List of $\scriptstyle \theta ^\star $-constructible sets}
  }
  \bigskip

Observe that the string $\sigma :=\delta _a=\{1,a\}$ is a proper subset of the string $\{1,a,aa\}$, and hence $\sigma $
is not maximal.  But yet notice that $\varphi _\sigma $ is an ultracharacter, since $\{\delta _a\}$ is a
  minimal\fn {Whenever $e_0$ is a nonzero minimal element of a {\sla } $E$, the character $\varphi (e)=[e_0\leq e]$ is
an ultracharacter.}
  member of $\P (S^\star ,\theta ^\star )$.
  We thus get an example of

\medskip {\it ``A string $\sigma $ which is not maximal but such that $\varphi _\sigma $ is an ultracharacter.''}
\medskip

\noindent On the other hand, since $\sigma =\sigma _{\varphi _\sigma }$, this also provides an example of

\medskip {\it ``An ultracharacter $\varphi $ such that $\sigma _\varphi $ is not maximal.''}  \medskip

This suggests the need to single out the strings which give rise to ultracharacters:

\definition
  We will say that a string $\sigma $ is \emph {quasi-maximal} whenever $\varphi _\sigma $ is an ultracharacter.
  The set of all quasi-maximal strings will be denoted by $S^\propto $.
  \# $S^\propto $; The set of all quasi-maximal strings;

Adopting this terminology, the conclusion of \ref {BigNewTightResult} states that every open, maximal string is
quasi-maximal.

\state Theorem \label UltraIsPhiSigma
  Let $S$ be a $0$-left cancellative semigroup admitting least common multiples.  Then, every open ultracharacter on
$\ehull $ is of the form $\varphi _\sigma $ for some open, quasi-maximal string $\sigma $.

\Proof
  Let $\varphi $ be an open ultracharacter on $\ehull $.  Letting $\sigma =\sigma _\varphi $, we have that $\sigma $ is
open by definition, and by \ref {RioMaior} it follows that $\varphi \leq \varphi _\sigma $, and hence $\varphi =\varphi
_\sigma $, by maximality. That $\sigma $ is a quasi-maximal string is due to the fact that $\varphi _\sigma $ is an
ultracharacter.
  \endProof

The importance of quasi-maximal strings evidenced by the last result begs for a better understanding of such strings.
While we are unable to provide a complete characterization, we can at least exhibit some further examples beyond the
maximal ones.

To explain what we mean, recalll from \ref {FRFstarR.i} that a string $\sigma $ belongs to some $F^\star _\Lambda $ if
and only if $\sigma $ is contained in $F^\theta _\Lambda $.  It is therefore possible that $\sigma $ is maximal among
all strings contained in $F^\theta _\Lambda $, and still not a maximal string.  An example is the string $\{1,a\}$
mentioned above, which is maximal within $F^\theta _a$, but not maximal in the strict sense of the word.

\state Proposition \label RelativelyMaximalStrings
  Let $\Lambda $ be a nonempty finite subset of $S$ and suppose that $\sigma $ is an open string such that $\sigma
\subseteq F^\theta _\Lambda $.  Suppose moreover that $\sigma $ is maximal among the strings contained in $F^\theta
_\Lambda $, in the sense that for every string $\mu $, one has that
  $$
  \sigma \subseteq \mu \subseteq F^\theta _\Lambda \Imply \sigma =\mu .
  $$
  Then $\varphi _\sigma $ is an ultracharacter, and hence $\sigma $ is a quasi-maximal string.

\Proof
  We begin by following the first steps of the proof of \ref {BigNewTightResult}:
  let $\psi $ be a character such that $\varphi _\sigma \leq \psi $.  For every $s$ in $S$ we then have that
  $$
  s\in \sigma \explain {SigmaPhiOnOpen.iii}\Imply \varphi _\sigma (E^\theta _s)=1 \Imply \psi (E^\theta _s)=1 \Imply
s\in \sigma _\psi ,
  $$
  which means that $\sigma \subseteq \sigma _\psi $.  By hypothesis we have that $\sigma \subseteq F^\theta _\Lambda $,
hence $\sigma \in F^\star _\Lambda $, by \ref {FRFstarR.i}, so
  $$
  1 =
  [\sigma \in F^\star _\Lambda ] =
  [\sigma \in \varepsilon (F^\theta _\Lambda )] =
  \varphi _\sigma (F^\theta _\Lambda ) \leq \psi (F^\theta _\Lambda ),
  $$
  and we conclude that $\psi (F^\theta _\Lambda )=1$.  We claim that this entails that $\sigma _\psi \subseteq F^\theta
_\Lambda $.  In fact, given any $s\in \sigma _\psi $, we have
  $$
  1 =
  \psi (E^\theta _s)\psi (F^\theta _\Lambda ) =
  \psi (E^\theta _s\cap F^\theta _\Lambda ),
  $$
  so $E^\theta _s\cap F^\theta _\Lambda $ cannot possibly be the empty set.  Picking any $r$ in $E^\theta _s\cap
F^\theta _\Lambda $, we have that $r=sx$, for some $x$ in $S$, and for every $t$ in $\Lambda $, one has
  $$
  0\neq tr = tsx,
  $$
  so in particular $ts\neq 0$, showing that $s\in F^\theta _\Lambda $.  We have thus proved that $\sigma \subseteq
\sigma _\psi \subseteq F^\theta _\Lambda $, and we deduce from the relative maximality of $\sigma $ that $\sigma =\sigma
_\psi $.

It follows that $\psi $ is an open character, so \ref {RioMaior} applies for $\psi $, and we deduce that
  $$
  \varphi _\sigma \leq \psi \explain {RioMaior}\leq \Phi \big (\Sigma (\psi )\big ) = \Phi (\sigma ) = \varphi _\sigma ,
  $$
  so $\varphi _\sigma =\psi $, proving that $\varphi _\sigma $ is maximal.
  \endProof

Here are some further questions we have come across and which are still to be answered:

\state Questions \rm
  \iItemize
  \iItem Is there an intrinsic characterization of quasi-maximal strings?
  \iItem Under which assumptions on $S$ is every quasi-maximal string maximal?
  \iItem Is it possible to characterize the strings $\sigma $ for which $\varphi _\sigma $ is a tight character?  These
should be called \emph {tight} strings.

\section Ground characters

\def \dualrep {\hat \theta } \def \ac #1#2{\dualrep _{#1}(#2)} \def \acinv #1#2{\dualrep _{#1}\inv (#2)} \def \exist
{\exists \kern 0.8pt}

In the last section we were able to fruitfully study open characters using strings, culminating with Theorem \ref
{UltraIsPhiSigma}, stating that every open ultracharacter is given in terms of a string.  However nothing of interest
was said about an ultracharacter when it is not open.  The main purpose of this section is thus to obtain some useful
information about non-open ultracharacters.  The main result in this direction is Theorem \ref {NonOpenUltra}, below.

\fix Throughout this section we fix a $0$-left cancellative semigroup $S$ admitting least common multiples.  For each
$s$ in $S$ let
  $$
  \hat F_ s = \{\varphi \in \spec : \varphi (F^\theta _s)=1\}, \and
  \hat E_ s = \{\varphi \in \spec : \varphi (E^\theta _s)=1\},
  $$
  \# $\dualrep $; Dual representation of $S$ on $\textbackslash spec$;
  \# $\hat F_s$; Domain of $\dualrep _s$;
  \# $\hat E_s$; Range of $\dualrep _s$;
  and for every $\varphi $ in $\hat F_ s$, consider the character $\ac s\varphi $ given by
  $$
  \ac s\varphi (X) = \varphi \big (\theta _s\inv (E^\theta _s\cap X)\big ), \for X\in \ehull .
  $$ Observing that
  $$
  \ac s\varphi (E^\theta _s) = \varphi \big (\theta _s\inv (E^\theta _s)\big ) = \varphi (F^\theta _s) = 1,
  \equationmark ActionOnPhiIsChar
  $$
  we see that $\ac s\varphi $ is indeed a (nonzero) character, and that $\ac s\varphi $ belongs to $\hat E_ s$.  As a
consequence we get a map
  $$
  \dualrep _ s \colon \hat F_ s \to \hat E_ s,
  $$
  which is easily seen to be bijective, with inverse given by
  $$
  \acinv s\varphi (X) =
  \varphi \big (\theta _s(F^\theta _s\cap X)\big ), \for \varphi \in \hat E_ s, \for X\in \ehull .
  $$
  We may then see each $\dualrep _ s$ as an element of $\I \big (\spec \big )$, and it is not hard to see that the
correspondence
  $$
  \dualrep : s\in S \mapsto \dualrep _ s\in \I (\spec )
  $$
  is a representation of $S$ on $\spec $.

All of this may also be deduced from the fact that any inverse semigroup, such as $\hull $, admits a canonical
representation on the spectrum of its idempotent {\sla } (see \cite [Section 10]{actions}), and that $\dualrep $ may be
obtained as the composition
  $$
  S\arw \theta \hull \longrightarrow \I \big (\spec \big ),
  $$
  where the arrow in the right-hand-side is the canonical representation mentioned above.

\definition We shall refer to $\dualrep $ as the \emph {dual representation} of $S$.

In order to study the relationship between the dual representation and the representation
  $\rho $ of $\hull $ described in \ref {MapForISG}, let us prove the following technical result.

\state Lemma \label CovarPhiLemma
  Given $s$ in $S$, and $\sigma $ in $S^\star _\ess $, one has that
  \iItemize
  \iItem $\varphi _\sigma \in \hat F_ s \IFF \sigma \in F^\star _s$,
  \iItem if the equivalent conditions in (i) are satisfied, then $\ac s{\varphi _\sigma }=\varphi _{\theta ^\star
_s(\sigma )}$,
  \iItem $\varphi _\sigma \in \hat E_ s \IFF \sigma \in E^\star _s$,
  \iItem if the equivalent conditions in (iii) are satisfied, then $\acinv s{\varphi _\sigma }=\varphi _{\theta ^{\star
-1}_s(\sigma )}$.

\Proof
  (i)\enspace We have
  $$
  \varphi _\sigma \in \hat F_ s \IFF
  \varphi _\sigma (F^\theta _s)=1 \IFF
  \sigma \in \varepsilon (F^\theta _s) \={IntroEpsilon} F^\star _s.
  $$

\itmproof (iii) Follows as above by replacing the letter ``$F$" by the letter ``$E$".

\itmproof (ii) Assuming (i), one has for every $X\in \ehull $, that
  $$
  \ac s{\varphi _\sigma }(X) =
  \varphi _\sigma \big (\theta _s\inv (E^\theta _s\cap X)\big ) =
  \Bool {\sigma \in \varepsilon \big (\theta _s\inv (E^\theta _s\cap X)\big )} \={MapForISG} $$ $$ =
  \Bool {\sigma \in \theta ^{\star -1}_s \big (\varepsilon (E^\theta _s\cap X)\big )} =
  \Bool {\theta ^\star _s (\sigma )\in \varepsilon (E^\theta _s\cap X)} = \cdots
  \equationmark DualPhiSigmaX
  $$

  Observe that
  $$
  \varepsilon (E^\theta _s\cap X)=\varepsilon (E^\theta _s)\cap \varepsilon (X) = E^\star _ s\cap \varepsilon (X),
  $$
  and since $\theta ^\star _s (\sigma )$ is evidently in $E^\star _ s$, one has that \ref {DualPhiSigmaX} coincides with
  $$
  \Bool {\theta ^\star _s (\sigma )\in \varepsilon (X)} = \varphi _{\theta ^\star _s(\sigma )}(X),
  $$
  thus proving (ii).

\itmproof (iv) Assuming (iii), one has for every $X\in \ehull $, that
  $$
  \acinv s{\varphi _\sigma }(X) =
  \varphi _\sigma \big (\theta _s(F^\theta _s\cap X)\big ) =
  \Bool {\sigma \in \varepsilon \big (\theta _s(F^\theta _s\cap X)\big )} \quebra =
  \Bool {\sigma \in \theta ^\star _s\big (\varepsilon (F^\theta _s\cap X)\big )} =
  \Bool {\theta ^{\star -1}_s(\sigma )\in \varepsilon (F^\theta _s\cap X)},
  $$
  and conclusion follows as in the proof of (ii).  \endProof

Considering the representation $\theta ^\star $ of $S$ on $S^\star $, observe that $S^\star _\ess $ is an invariant\fn
  {The essential subset for a representation is evidently invariant!}
  subset of $S^\star $, and it is easy to see that it is also invariant under the representation
  $\rho $ of $\hull $ described in \ref {MapForISG}.  Together with the dual representation of $\hull $ on $\spec $
mentioned above, we thus have two natural representations of $\hull $, which are closeely related, as the following
immediate consequence of the above result asserts:

\state Proposition \label CovarPhi
  The mapping
  $$
  \Phi \colon S^\star _\ess \to \spec
  $$
  of \ref {DefinePhi} is covariant relative to the natural representations of $\hull $ referred to above.

\Proof
  Follows immediately from \ref {CovarPhiLemma}, and the fact that $\hull $ is generated by the $\theta _s$ and their
inverses.
  \endProof

The fact that the correspondence between strings and characters
  (see e.g.~\ref {SigmaPhiIsInterior} and \ref {RioMaior})
  is not a perfect one
  is partly responsible for the fact that expressing the covariance properties of the map $\Sigma $ of \ref {IntroSigma}
cannot be done in the same straightforward way as we did for $\Phi $ in \ref {CovarPhi}.

Let us first treat the question of covariance regarding $\acinv s\varphi $.  Of course, for this to be a well defined
character we need $\varphi $ to be in $\hat E_ s$, meaning that $\varphi (E^\theta _s)=1$, which is also equivalent to
saying that $s\in \sigma _\varphi $.  In particular characters with empty strings are immediately ruled out.

\state Lemma \label BackInvarString
  For every $s$ in $S$, and every character $\varphi $ in $\hat E_ s$, one has that
  $$
  \sigma _{\acinv s\varphi } = \{p\in S: sp\in \sigma _\varphi \}.
  $$

\Proof
  For $p$ in $S$, notice that $p\in \sigma _{\acinv s\varphi }$, iff
  $$
  1=\acinv s\varphi (E^\theta _p) =
  \varphi \big (\theta _s(F^\theta _s\cap E^\theta _p)\big ) =
  \varphi (E^\theta _{sp}),
  $$
  which in turn is equivalent to saying that $sp\in \sigma _\varphi $.  \endProof

The set appearing in the right hand side of the equation displayed in \ref {BackInvarString} is precisely the same set
mentioned in definition \ref {MappingStrings.ii} of $s\inv *\sigma _\varphi $, except that this notation is reserved for
the situation in which the intersection of $\sigma $ with $sS$ is nonempty, which precisely means that $s\in \interior
\sigma _\varphi $.

\state Proposition \label DualBackOnStrings
  Pick $s$ in $S$ and let $\varphi $ be any character in $\hat E_ s$.  Then $s\in \sigma _\varphi $, and moreover
  \iItemize
  \iItem if $s$ is in $\interior \sigma _\varphi $, then $\sigma _\varphi \in E^\star _s$, and $\sigma _{\acinv s\varphi
} = \theta ^{\star -1}_s(\sigma _\varphi )$,
  \iItem if $s$ is not in $\interior \sigma _\varphi $, then $\sigma _{\acinv s\varphi } = \emptyset $.

\Proof
    If $s$ is in the interior of $\sigma _\varphi $, there is some $p$ in $S$ such that $sp\in \sigma _\varphi $, whence
also $sp\in \sigma _\varphi \cap E^\theta _s$.  It follows that $\sigma _\varphi \cap E^\theta _s$ is nonempty, so \ref
{FRFstarR.ii} implies that $\sigma _\varphi $ lies in $E^\star _s$.  The last statement of (i) then follows at once from
\ref {BackInvarString}.

On the other hand, if $s$ is not in the interior of $\sigma _\varphi $, the conclusion again follows from \ref
{BackInvarString}.  \endProof

Regarding the behavior of strings associated to characters of the form $\ac s\varphi $, we have:

\state Lemma \label BirthOfString
  For every $s$ in $S$, and every character $\varphi $ in $\hat F_ s$, one has that $\ac s\varphi $ belongs to $\seone $
(and hence \ref {DistinguishSeone} implies that $\sigma _{\ac s\varphi }$ is a string), and moreover
  \iItemize
  \iItem if $\sigma _\varphi $ is nonempty, then $\sigma _\varphi \in F^\star _s$, and $\sigma _{\ac s\varphi } = \theta
^\star _s(\sigma _\varphi )$,
  \iItem if $\sigma _\varphi $ is empty, then $\sigma _{\ac s\varphi } = \delta _s$.

\Proof The fact that $\ac s\varphi $ belongs to $\seone $ follows from \ref {ActionOnPhiIsChar} and \ref
{DistinguishSeone}.

Given any $t$ in $\sigma _\varphi $, we have that $\varphi (E^\theta _t)=1$, so
  $$
  \varphi (E^\theta _t\cap F^\theta _s)= \varphi (E^\theta _t)\varphi (F^\theta _s) =1,
  $$
  whence $E^\theta _t\cap F^\theta _s\neq \emptyset $.  Analyzing any element in this nonempty intersection one quickly
realizes that $t\in F^\theta _s$, from whence it follows that $\sigma _\varphi \subseteq F^\theta _s$.  Under the
assumption of (i) one has that $\sigma _\varphi $ is a string by \ref {StringFromChar}, and the first conclusion of (i)
then follows from \ref {FRFstarR.i}.

Regardless of any assumption on $\sigma _\varphi $, pick any $t$ in $S$, and let $r$ be a least common multiple of $s$
and $t$, so that $E^\theta _r = E^\theta _s\cap E^\theta _t$, and $r = su=tv$, for suitable $u$ and $v$ in $\tS $.  We
then have that
  $$
  t\in \sigma _{\ac s\varphi } \iff
  \ac s\varphi (E^\theta _t) = 1 \quebra \iff
  \varphi \big (\theta _s\inv (E^\theta _s\cap E^\theta _t)\big )= 1 \iff
  \varphi \big (\theta _s\inv (E^\theta _r)\big )= 1 \iff \cdots
  $$
  Notice that
  $$
  E^\theta _r =
  E^\theta _{su} =
  \theta _s(E^\theta _u \cap F^\theta _s),
  $$
  so,
  $$
  \varphi \big (\theta _s\inv (E^\theta _r)\big ) =
  \varphi \big (\theta _s\inv \theta _s(E^\theta _u \cap F^\theta _s)\big ) =
  \varphi (E^\theta _u \cap F^\theta _s)
  $$
  We then deduce that
  $$
  t\in \sigma _{\ac s\varphi } \iff \varphi (E^\theta _u \cap F^\theta _s)=1.
  \equationmark CondFortInDualAct
  $$

In case $t\|s$, in which case one may take $r=s$, and $u=1$, above, one obviously has that $\varphi (E^\theta _u \cap
F^\theta _s)=1$, so the above argument shows that $t\in \sigma _{\ac s\varphi }$, thus proving that
  $$
  \delta _s\subseteq \sigma _{\ac s\varphi }.
  $$

Now suppose that $t {\not |}\, s$, and that we are under the conditions of (ii).  Then the element $u$ above is
necessarily different from $1$, so it must lie in $S$.  Moreover $\varphi (E^\theta _u)=0$, since otherwise $u\in \sigma
_\varphi $, whence
  $$
  \varphi (E^\theta _u \cap F^\theta _s) = \varphi (E^\theta _u)\varphi (F^\theta _s) = 0,
  $$
  and then $t\notin \sigma _{\ac s\varphi }$, by \ref {CondFortInDualAct}.  This concludes the proof of (ii), and it now
remains to prove the last assertion of (ii).  In order to do this we first observe that $\theta ^\star _s(\sigma
_\varphi )\in E^\star _s$, so $s\in \theta ^\star _s(\sigma _\varphi )$ by \ref {FRFstarR.iii}, and hence also
  $$
  \delta _s\subseteq \theta ^\star _s(\sigma _\varphi ).
  \equationmark DivisorIn
  $$
  Choosing any $t\in \sigma _{\ac s\varphi }$, let $r$, $u$, and $v$ be as above.  If $u=1$, then $t\|s$, so $t\in
\theta ^\star _s(\sigma _\varphi )$, by \ref {DivisorIn}.  Henceforth assuming that $u\in S$, we have by \ref
{CondFortInDualAct} that $\varphi (E^\theta _u) = 1$, so $u\in \sigma _\varphi $, and since $t\|su$, it follows that
$t\in \theta ^\star _s(\sigma _\varphi )$, thus showing that $\sigma _{\ac s\varphi }\subseteq \theta ^\star _s(\sigma
_\varphi )$.

In order to prove the reverse inclusion, pick any $t\in \theta ^\star _s(\sigma _\varphi )$.  Then $t\|sp$, for some $p$
in $\sigma _\varphi $, so we may write $tx=sp$, for some $x$ in $\tS $.  A moment's reflection will convince the reader
that
  $$
  \theta _s(E^\theta _p\cap F^\theta _s)\subseteq E^\theta _s \cap E^\theta _t,
  $$
  so
  $$
  \ac s\varphi (E^\theta _t) =
  \varphi \big (\theta _s\inv (E^\theta _s\cap E^\theta _t)\big ) \geq
  \varphi (E^\theta _p\cap F^\theta _s) =
  \varphi (E^\theta _p)\varphi (F^\theta _s) =
  \varphi (E^\theta _p) =1.
  $$
  This shows that $t\in \sigma _{\ac s\varphi }$, concluding the proof of (ii).  \endProof

We may interpret the above result, and more specifically the identity
  $$
  \sigma _{\ac s\varphi } = \theta ^\star _s(\sigma _\varphi ),
  $$
  as saying that the correspondence $\varphi \mapsto \sigma _\varphi $ is covariant with respect to the actions
$\dualrep $ and $\theta ^\star $, on $\spec $ and $S^\star $, respectively, except that the term ``$\sigma _\varphi $"
appearing is the right-hand-side above is not a well defined string since it may be empty, even though the
left-hand-side is always well defined.  In the problematic case of an empty string, \ref {BirthOfString.ii} then gives
the undefined right-hand-side the default value of $\delta _s$.

\definition \label DefineGroundChars
  A character $\varphi $ in $\spec $ will be called a \emph {ground} character if $\sigma _\varphi $ is empty.

By \ref {DistinguishSeone}, the ground characters are precisely the members of $\spec \setminus \seone $.

Besides the ground characters, a character $\varphi $ may fail to be open because $\sigma _\varphi $, while being a bona
fide string, is not an open string.  In this case we have by \ref {PropsOpenString.ii} that $\sigma _\varphi =\delta
_s$, for some $s$ in $S$ such that $s\notin sS$.

\state Lemma \label PrepOrbitGroundChar
  Let $\varphi $ be a character such that $\sigma _\varphi =\delta _s$, where $s$ is such that $s\notin sS$.  Then
$\varphi \in \hat E_ s$, and $\acinv s\varphi $ is a ground character.

\Proof Since $s\in \sigma _\varphi $, we have that $\varphi (E^\theta _s)=1$, so $\varphi \in \hat E_ s$, whence
  $$
  \psi :=\acinv s\varphi
  $$
  is a well defined character.  In order to prove that $\psi $ is a ground character, we argue by contradiction and
suppose instead that there exists some $t\in \sigma _\psi $.  Observing that $\psi \in \hat F_ s$, we then have that
  $$
  st\in s\sigma _\psi \subseteq s*\sigma _\psi =\theta ^\star _s(\psi )\={BirthOfString.i}
  \sigma _{\dualrep _s(\psi )}
  = \sigma _\varphi = \delta _s.
  $$
  It follows that $st\|s$, whence $stu=s$, for some $u$ in $\tS $, so
  $$
  s=stu\in sS,
  $$
  contradicting the hypothesis.  This shows that $\psi $ is indeed a ground character, concluding the proof.  \endProof

We may now give a precise characterization of non-open characters in terms of the ground characters:

\state Proposition \label NonOpenChars
  Denote by \
  $\sop $
  \# $\sop $; Set of all open characters on $\ehull $;
  the set of all open characters on $\ehull $.  Then
  $$
  \spec \setminus \sop = \big \{\ac u\varphi : u\in \tS , \ \varphi \hbox { is a ground character in } \hat F_ u\big \}.
  $$
  Moreover for each $\psi $ in the above set, there is a unique pair $(u, \varphi )$, with $u$ in $\tS $, and $\varphi $
a ground character, such that $\psi =\ac u\varphi $.

\Proof If $\varphi $ is a non-open character, then $\sigma _\varphi $ is either empty, in which case $\varphi $ is a
ground character, or $\sigma _\varphi $ is a (nonempty) non-open string.
  In the latter case we have by \ref {PropsOpenString.ii} that $\sigma _\varphi =\delta _s$, for some $s$ in $S$, such
that $s\notin sS$, so it follows from \ref {PrepOrbitGroundChar} that
  $$
  \psi := \acinv s\varphi
  $$
  is a ground character, necessarily in $\hat F_ s$.  Observing that $\varphi =\ac s\psi $, we see that $\varphi $ lies
in the set appearing in the right-hand-side in the statement.

Conversely, if $\varphi $ is a ground character in $\hat F_ u$, we must show that $\ac u\varphi $ is not open.  In case
$u=1$, there is nothing to do since ground characters are obviously not open, so we henceforth suppose that $u\in S$.
In keeping with our tradition of naming elements in $S$ by $s$, $t$, and $r$, while reserving $u$ and $v$ for elements
which, in principle, are allowed to range in all of $\tS $, we will write $s$ for $u$, so that $\varphi \in \hat F_ s$,
and our task consists in showing that $\ac s\varphi $ is not open.

Since $\varphi $ is a ground character, we have by \ref {BirthOfString.ii} that
  $$
  \sigma _{\ac s\varphi }=\delta _s.
  $$ Arguing by contradiction, suppose that $\ac s\varphi $ is open, whence $\delta _s$ is an open string.  Observing
that $s$ lies in $\delta _s$, there exists some $p$ in $S$ such that $sp\in \delta _s$, and then we may find $x$ in $\tS
$ such that $spx=s$.  Letting $e=px$, we then claim that $F^\theta _s\subseteq E^\theta _e$.  To see this, it is enough
to notice that if $t\in F^\theta _s$, then
  $$
  0\neq st=set,
  $$
  so $0$-left cancellativity applies giving $t=et\in E^\theta _e$, and proving our claim.  Recalling that $\varphi \in
\hat F_ s$, we then have
  $$
  1=\varphi (F^\theta _s)\leq \varphi (E^\theta _e),
  $$
  whence $e\in \sigma _\varphi =\emptyset $, a contradiction.  This shows that $\ac s\varphi $ is not open, as desired.

To prove the last assertion in the statement, suppose that
  $$
  \ac {u_1}{\varphi _1}=\ac {u_2}{\varphi _2},
  $$
  where $u_1, u_2\in \tS $, and $\varphi _1$ and $\varphi _2$ are ground characters.  We first observe that one cannot
have $u_1\in S$, and $u_2=1$ (or vice-versa), since otherwise $\ac {u_1}{\varphi _1}$ is not a ground character by \ref
{BirthOfString}, so it cannot possibly coincide with the ground character $\varphi _2$.  If both $u_1$ and $u_2$
coincide with $1$, there is nothing to do, so we suppose from now on that
  $$
  s_i:=u_i\in S, \for i=1,2,
  $$
  hence our hypothesis reads:
  $$
  \ac {s_1}{\varphi _1}=\ac {s_2}{\varphi _2}.
  $$

Using \ref {BirthOfString.ii} we have that
  $$
  \delta _{s_1} =
  \sigma _{ \ac {s_1}{\varphi _1}} =
  \sigma _{\ac {s_2}{\varphi _2}} =
  \delta _{s_2}.
  $$
  By the first part of the proof we have that $\ac {s_1}{\varphi _1}$ is not open, so the string displayed above is
likewise not open, and we deduce from the uniqueness in \ref {PropsOpenString.ii} that $s_1=s_2$.  The fact that
$\varphi _1=\varphi _2$ now follows easily.
  \endProof

We may now combine several of our earlier results to give a description of all ultracharacters on $\ehull $:

\state Theorem \label NonOpenUltra
  Let $S$ be a $0$-left cancellative semigroup admitting least common multiples.  Denote by
  $\sinf $
  \# $\sinf $; Set of all ultracharacters on $\ehull $;
  the set of all ultracharacters on $\ehull $, and by
  $$
  \sinfop = \sop \cap \sinf ,
  $$
  \# $\sinfop $; Set of all open ultracharacters on $\ehull $;
  namely the subset formed by all open ultracharacters.  Then
  \iItemize
  \iItem $\sinfop \phantom {\setminus \sinfop \ } = \big \{\varphi _\sigma : \sigma \hbox { is an open, quasi-maximal
string in }S\big \}$, and
  \iItem $\sinf \setminus \sinfop = \big \{\ac u\varphi : u\in \tS , \ \varphi \hbox { is a ground, ultracharacter in }
\hat F_ u\big \}$.

\Proof In order to prove (i) observe that
  if $\varphi $ is an open ultracharacter, then $\varphi =\varphi _\sigma $ for some open, quasi-maximal string $\sigma
$ by \ref {UltraIsPhiSigma}.  On the other hand, if $\sigma $ is an open, quasi-maximal string, then $\varphi _\sigma $
is an ultracharacter by definition, while
  $$
  \sigma =\Sigma \big (\Phi (\sigma )\big ),
  $$
  by \ref {SigmaPhiOnOpen.iii}, so $\varphi _\sigma = \Phi (\sigma )$ is an open character.

By \cite [Proposition 3.2]{ExelPardo}, one has that the dual representation of $\hull $ on $\spec $ leaves the set of
ultracharacters invariant, so a character $\varphi $ in $\hat F_ u$ is an ultracharacter if and only if the same holds
for $\ac u\varphi $.
  This said, point (ii) follows immediately from \ref {NonOpenChars}.  \endProof

The upshot of all this is that, in order to fully understand the ultracharacters of $\ehull $, one first needs to
describe the open, quasi-maximal strings.  The remaining ultracharacters are therefore obtained as the orbit under $\hat
\theta $ of the ground, ultracharacters.

Should the above program be brought to completion, one would therefore be able to describe all tight characters, since
$\varspec {_\tight }$ is well known to be the closure of the set of ultracharacters.

\immediate \closeout 3

  \IfFileExists {./symbols.aux}{
    \bigbreak
    \centerline {\tensc List of symbols}
    \nobreak \bigskip
    \catcode `\@=11
    \input symbols.aux
    \catcode `\@=12
    \vskip 1cm
  }{\bigskip \noindent \vtt *** File ``symbols.aux" is missing.  It might be generated after a rerun. *** \bigskip }

\writeIndex 2{contents}{}{References}

\references

%* AielloContiRossiStammeier
\Bibitem AielloContiRossiStammeier
  V. Aiello, R. Conti, S. Rossi and N. Stammeier;
  The inner structure of boundary quotients of right LCM semigroups;
  arXiv e-prints, 1709.08839, 2017

%* Bedos
\Bibitem BedosSpielberg
  E. B\'edos, S. Kaliszewski, J. Quigg and J. Spielberg;
  On finitely aligned left cancellative small categories, Zappa-Sz\'ep products and Exel-Pardo algebras;
  arXiv e-prints, 1712.09432, 2017

%* BrownloweLarsenStammeier
\Bibitem BrownloweLarsenStammeier
  N. Brownlowe, N. S. Larsen and N. Stammeier;
  C*-Algebras of algebraic dynamical systems and right LCM semigroups;
  to appear in Indiana Univ. Math. J.,
  arXiv e-prints, 1503.01599, 2017

%* CarlsenMatsumoto
\Article MatsuCarl
  T. M. Carlsen and K. Matsumoto;
  Some remarks on the C*-algebras associated with subshifts;
  Math. Scand., 95 (2004), 145-160

%* CarlsenSilvestrov
\Article CarlsenSilvestrov
  T. M. Carlsen and S. Silvestrov;
  C*-crossed products and shift spaces;
  Expo. Math., 25 (2007), 275-307

%* Cherubini
\Article Cherubini
  A. Cherubini and M. Petrich;
  The Inverse Hull of Right Cancellative Semigroups;
  J. Algebra, 111 (1987), 74-113

%* Clifford
\Bibitem CP
  A.~H. Clifford and G.~B. Preston;
  The algebraic theory of semigroups. Vol. I;
  Mathematical Surveys, No. 7.;
  American Mathematical Society,
  Providence, R.I., 1961

%* Coburn1
\Article CoOne
  L. A. Coburn;
  The C*-algebra generated by an isometry I;
  Bull. Amer. Math. Soc., 73 (1967), 722-726

%* Coburn2
\Article CoTwo
  L. A. Coburn;
  The C*-algebra generated by an isometry II;
  Trans. Amer. Math. Soc., 137 (1969), 211-217

%* Dokuchaev
\Article DokuchaExel
  M. Dokuchaev and R. Exel;
  Partial actions and subshifts;
  J. Funct. Analysis, 272 (2017), 5038-5106

%* Exela1
\Article actions
  R. Exel;
  Inverse semigroups and combinatorial C*-algebras;
  Bull. Braz. Math. Soc. (N.S.), 39 (2008), 191-313

%* Exela2
\Article semigpds
  R. Exel;
  Semigroupoid C*-Algebras;
  J. Math. Anal. Appl., 377 (2011), 303-318

%* ExelLaca
\Article infinoa
  R. Exel and M. Laca;
  Cuntz-Krieger algebras for infinite matrices;
  J. reine angew. Math., 512 (1999), 119-172

%* ExelPardo
\Article ExelPardo
  R. Exel and E. Pardo;
  The tight groupoid of an inverse semigroup;
  Semigroup Forum, 92 (2016), 274-303

%* ExelSteinberg
\Bibitem announce
R. Exel and B. Steinberg;
The inverse hull of 0-left cancellative semigroups;
arXiv e-prints,\hfill \break 1710.04722, 2017

%* Howson
\Article Howson
  A. G. Howson;
  On the intersection of finitely generated free groups;
  J. London Math. Soc., 29 (1954), 428-434

%* Keimel
\Article Keimel
  K. Keimel;
  Alg\`ebres commutatives engendr\'ees par leurs \'el\'ements idempotents;
  Canad. J. Math.,  22 (1970), 1071-1078

%* KPRR
\Article KPRR
  A. Kumjian, D. Pask, I. Raeburn and J.  Renault;
  Graphs, groupoids, and Cuntz-Krieger algebras;
  J. Funct. Anal., 144 (1997), 505-541

%* KumjianPask
\Article KumjianPask
  A. Kumjian and  D. Pask;
  Higher rank graph C*-algebras;
  New York J. Math., 6 (2000), 1-20

%* KwasniewskiLarsen
\Bibitem KwasniewskiLarsen
  B. K. Kwasniewski and N. S. Larsen;
  Nica-Toeplitz algebras associated with right tensor C*-precategories over right LCM semigroups: part II examples;
  arXiv e-prints, 1706.04951, 2017

%* Li
\Article Li
  X. Li;
  Semigroup C*-algebras and amenability of semigroups;
  J. Funct. Anal., 262 (2012), 4302-4340

%* Matsumoto
\Article MatsuOri
  K. Matsumoto;
  On C*-algebras associated with subshifts;
  Internat. J. Math., 8 (1997), 357-374

%* Munn
\Article Munn
  W. D. Munn;
  Brandt congruences on inverse semigroups;
  Proc. London Math. Soc., {\rm (3)} 14 (1964), 154-164.

%* Murphy1
\Article MurOne
  G. J. Murphy;
  Ordered groups and Toeplitz algebras;
  J. Operator Theory,   18 (1987), 303-326

%* Murphy2
\Article MurTwo
  G. J. Murphy;
  Ordered groups and crossed products of C*-algebras;
  Pacific J. Math., 2 (1991), 319-349

%* Murphy3
\Article MurThree
  G. J. Murphy;
  Crossed products of C*-algebras by semigroups of automorphisms;
  Proc.  London Math. Soc., 3 (1994), 423-448

%* Nica
\Article Nica
  A. Nica;
  C*-algebras generated by isometries and Wiener-Hopf operators;
  J. Operator Theory, 27 (1992), 17-52

%* Paterson
\Bibitem Paterson
  A. L. T. Paterson;
  Groupoids, inverse semigroups, and their operator algebras;
  Birkh\umlaut auser, 1999

%* SpielbergA
\Bibitem SpielbergA
  J. Spielberg;
  Groupoids and C*-algebras for categories of paths;
  arXiv e-prints, 1111.6924v4, 2014

%* SpielbergB
\Bibitem SpielbergB
  J. Spielberg;
  Groupoids and $C^*$-algebras for left cancellative small categories;
  arXiv e-prints, 1712.07720, 2017

%* Stammeier
\Bibitem Stammeier
  N. Stammeier;
  A boundary quotient diagram for right LCM semigroups;
  Semigroup Forum, pp. 1-16 (2017, first online),
  arXiv e-prints, 1604.03172

%* Starling
\Bibitem Starling
  C. Starling;
  Boundary quotients of C*-algebras of right LCM semigroups;
  arXiv e-prints, 1409.1549, 2017

%* Steinberg
\Article SteinbergPrimitive
B. Steinberg;
Simplicity, primitivity and semiprimitivity of \'etale groupoid algebras with applications to inverse semigroup algebras;
J. Pure Appl. Algebra, 220 (2016), 1035-1054

\endgroup

\closeout 2
\close
\bye